\documentclass[12pt]{article}

\usepackage{amssymb,latexsym,amsmath}

\begin{document}

\sloppy
\renewcommand{\theequation}{\arabic{section}.\arabic{equation}}
\thinmuskip = 0.5\thinmuskip
\medmuskip = 0.5\medmuskip
\thickmuskip = 0.5\thickmuskip
\arraycolsep = 0.3\arraycolsep

\newtheorem{theorem}{Theorem}[section]
\newtheorem{lemma}[theorem]{Lemma}
\newtheorem{prop}[theorem]{Proposition}
\renewcommand{\thetheorem}{\arabic{section}.\arabic{theorem}}

\newcommand{\prf}{\noindent{\bf Proof.}\ }
\def\prfe{\hspace*{\fill} $\Box$

\smallskip \noindent}

\def\be{\begin{equation}}
\def\ee{\end{equation}}
\def\bea{\begin{eqnarray}}
\def\eea{\end{eqnarray}}
\def\beas{\begin{eqnarray*}}
\def\eeas{\end{eqnarray*}}

\newcommand{\R}{\mathbb R}
\newcommand{\N}{\mathbb N}
\newcommand{\T}{\mathbb T}

\def\supp{\mathrm{supp}\,}
\def\vol{\mathrm{vol}\,}
\def\sign{\mathrm{sign}\,}
\def\ekin{E_\mathrm{kin}}
\def\epot{E_\mathrm{pot}}

\def\a{a_{\rho}}
\def\aa{a_{\psi}}
\def\B{B_{\rho}}
\def\c{c_{\rho}}
\def\as{a_{\sigma}}
\def\Bs{B_{\sigma}}
\def\cs{c_{\sigma}}

\def\am{a_{\rho^m}}
\def\Bm{B_{\rho^m}}
\def\cm{c_{\rho^m}}
\def\dm{d_{\rho^m}}
\def\Em{e_{\rho^m}}
\def\g{d_{\rho}}
\def\e{e_{\rho}}
\def\d{\partial^{\mu}_s}
\def\A{{\cal A}}
\def\E{{\cal E}}
\def\K{{\cal K}_{\nu}}
\def\H{{\cal H}}
\def\Hc{{{\cal H}_C}}
\def\D{{\cal D}}
\def\Db{\overline{{\cal D}}}
\def\Eb{\overline{E}}
\def\EE{\overline{{\cal E}}}
\def\U{{\cal U}}
\def\A{{\cal A}}
\def\f{\left\langle f\right\rangle}
\def\sm{\left\langle \rho^m\right\rangle}
\def\s{\left\langle\rho\right\rangle}
\def\p{\left\langle\psi\right\rangle}
\def\si{\left\langle\sigma\right\rangle}
\def\open#1{\setbox0=\hbox{$#1$}
\baselineskip = 0pt
\vbox{\hbox{\hspace*{0.4 \wd0}\tiny $\circ$}\hbox{$#1$}}
\baselineskip = 11pt\!}
\def\fn{\open{f}}
\def\Rn{\open{R}}
\def\Pn{\open{P}}

\title{Stability in the Stefan problem with surface tension (I)}
\author{ Mahir Had\v{z}i\'{c} and Yan Guo\\
         Division of Applied Mathematics \\
         Brown University, Providence, RI 02912, U.S.A.}
\maketitle

\begin{abstract}
We develop a high-order energy method to prove asymptotic stability of flat steady surfaces
for the Stefan problem with surface tension - also known as the Stefan problem with
Gibbs-Thomson correction.
\end{abstract}
\section{Introduction}
\setcounter{equation}{0}
The Stefan problem is one of the best known parabolic two-phase free boundary problems.
It is a simple model of phase transitions in liquid-solid systems.
\par
Let $\Omega\subset\R^n$ denote a domain that contains a liquid and a solid separated by an interface
$\Gamma$. As the melting or cooling take place the boundary moves and we are naturally led
to a free boundary problem.
The unknowns are the temperatures of the liquid and the solid denoted respectively by $v^+$ and $v^-$
and the location of the interface $\Gamma$ separating the two different phases.
\par

We shall assume that $\Omega=\T^{n-1}\times[-1,1]$ where $\T^{n-1}$ stands for an $(n-1)$-dimensional torus.
Let us assume that the moving interface $\Gamma(t)$ is a graph given by
$x_n=\rho(t,x')$. Here $\rho:[0,T]\times\T^{n-1}\to\R$ is some smooth function
such that $\bigcup_{0\leq t\leq T}\Gamma(t)\subset\Omega$ and $T>0$. Define the liquid/solid phase $\Omega^{\pm}(t)$ by
setting
$
\Omega^{\pm}(t)=\Big\{(x',x_n)\in\Omega\Big|\,\,x_n\gtrless\rho(x',t)\Big\}.
$
We note that $\Omega=\Omega^+(t)\cup\Omega^-(t)$. In order to formulate the problem we first specify the initial conditions.
Let $\Gamma_0=\textrm{graph}(\rho_0)$ be the initial position of the free boundary and $v_0:\Omega\to\R$ be the initial temperature.
The unknowns are the interface $\Big\{\Gamma(t);\,t\geq0\Big\}$ and the temperature function
$
v:\Omega\times[0,T]\to\R.
$
We denote the normal velocity of
$\Gamma$ by $V$ and normalize it to be positive if $\Gamma$ is locally expanding $\Omega^-(t)$.
The mean curvature of $\Gamma(t)$ is given by
\[
\kappa(t)=\nabla\cdot\big(\frac{\nabla\rho(t)}{\sqrt{1+|\nabla\rho(t)|^2}}\big).
\]
The Stefan problem with surface tension is now given by:
\bea
\label{eq:T}
&&\partial_tv-\Delta v=0\quad\quad\,\textrm{in}\,\,\Omega,\,\, t>0,\\
\label{eq:D}
&&v=\kappa(t)\quad\quad\quad\quad\,\,\textrm{on}\,\,\Gamma(t),\,\, t\geq0,\\
&&\partial_nv=0\quad\quad\quad\quad\,\,\,\,\,\textrm{on}\,\,\T^{n-1}\times\{x_n=\pm1\}\\
\label{eq:N}
&&V=[\partial_{\nu}v]^+_-\quad\quad\quad\textrm{on}\,\,\Gamma(t),\,\, t>0,\\
&&v(\cdot,0)=v_0\quad\quad\quad\,\textrm{in}\,\,\Omega,\\
\label{eq:rhoin}
&&\Gamma(0)=\Gamma_0.
\eea
Given $v$ we write
$v^+$ and $v^-$ for the restriction of $v$ to $\Omega^+(t)$ and $\Omega^-(t)$,
respectively. With this notation
$[\partial_{\nu}v]^+_-$ stands for the jump
of the normal derivatives across the interface $\Gamma(t)$, namely
\[
[\partial_{\nu}v]^+_-:=\partial_{\nu}v^+-\partial_{\nu}v^-,
\]
where $\nu$ stands for the unit normal on the hypersurface $\Gamma(t)$ with respect to $\Omega^-(t)$.
If we replace the boundary condition~(\ref{eq:D}) with
\be
v=0\quad\quad\,\,\textrm{on}\quad\Gamma(t),\,\,\,\,\,\, t\geq0,
\ee
then we are referring to the {\em classical Stefan problem}.
\par
The difficulties in dealing with the existence of solutions of the
problem~(\ref{eq:T}) -~(\ref{eq:rhoin}) arise from
the nonlinear coupling between the temperature $v$ and the boundary $\rho$.
This connection is expressed through the boundary conditions~(\ref{eq:D})
and~(\ref{eq:N}). The equation~(\ref{eq:N}) is a Neumann-type boundary
condition for $v$. It is hyperbolic in nature as opposed to the parabolic
diffusion process in the regions $\Omega^+$ and $\Omega^-$.
\par
From the technical point of view the first major obstacle for the analysis is the moving boundary. To deal with
this issue we shall first transform the problem to the fixed domain by applying
the so-called Hanzawa transform. To this end let us fix a small positive constant
$\alpha<\frac{1}{3}$ and choose a cut-off function
$\phi\in C^{\infty}(\R)$ with $\textrm{Im}(\phi)\subset[0,1]$ and
\[
\phi(z)=\left\{
\begin{array}{l}
1,\quad|z|\leq \alpha,\\
0,\quad|z|> 1-\alpha
\end{array}
\right.
\quad
||\phi'||_{L^{\infty}(\R)}<C.
\]
Define now a diffeomorphism
\[
\Theta(x',x_n,t)=(x',x_n+\phi(x_n)\rho(x',t),t)
\]
and the function $u(x',x_n,t)=v(\Theta(x',x_n,t))$. Observe that $u(x',0,t)=v(x',\rho(x',t),t)$ and at the
outer boundaries $\partial\Omega^{\pm}:=\T^{n-1}\times\{x_n=\pm1\}$, we have
$v|_{\partial\Omega^{\pm}}=u|_{\partial\Omega^{\pm}}$.
This is the version of the transform first introduced by Hanzawa (cf.~\cite{Ha}). In the new coordinates
the heat operator $\partial_t-\Delta$ transforms into a more complicated operator whose coefficients depend on
the interface function $\rho$ and the cut-off function $\phi$. Following the calculations from \cite{EsPrSi}, we
find that the
Laplace operator $\Delta$ in the new coordinates takes the form
\[
\Delta_{\Theta}u=\Delta_{x'}u+\a u_{nn}-\B\cdot\nabla_{x'}u_n-\g u_n,
\]
where
\be\label{eq:abro}
\a:=\frac{1+|\phi\nabla\rho|^2}{(1+\phi'\rho)^2},\quad\B:=\frac{2\phi}{1+\phi'\rho}\nabla\rho,
\ee
\be\label{eq:dro}
\g:=\frac{\phi\Delta\rho}{1+\phi'\rho}-
\frac{(\phi^2)'|\nabla\rho|^2}{(1+\phi'\rho)^2}+
\frac{\phi''\rho(1+|\phi\nabla\rho|^2)}{(1+\phi'\rho)^3}.
\ee
Furthermore, the operator $\partial_t$ in the new coordinates reads
\[
(\partial_t)_{\Theta}u=\partial_tu+\e u_n
\]
where
\be\label{eq:ero}
\e:=-\frac{\phi\rho_t}{1+\phi'\rho}.
\ee
Note that RHS of~(\ref{eq:D}) remains unchanged in the new coordinates.
In order to transform the boundary condition~(\ref{eq:N}) into the new coordinates we first observe
that
\[
(\partial_i)_{\Theta}u=\partial_iu-\frac{\phi\partial_i\rho}{1+\phi'\rho}\partial_nu,\,\,\,\,1\leq i\leq n-1,\quad\quad
(\partial_n)_{\Theta}u=\frac{1}{1+\phi'\rho}\partial_nu.
\]
We thus conclude that at the boundary $\T^{n-1}\times\{x_n=0\}$:
\[
\nabla v|_{\Gamma}=\nabla_{\Theta}u|_{x_n=0}=(\nabla_{x'}u,0)-\frac{\partial_n}{1+\phi'\rho}
(\phi\nabla\rho,1)=(\nabla_{x'}u,0)-\partial_nu(\nabla\rho,1).
\]
The outward unit normal is given by $\nu(x',t)=\frac{(-\nabla\rho(x',t),1)}{\sqrt{1+|\nabla\rho|^2}}$. Thus the normal
velocity $V$ takes the form
$
V=\frac{-\partial_t\rho}{\sqrt{1+|\nabla\rho|^2}}.
$
Using the above expressions we derive the formula for $[\partial_nv]^+_-|_{\Gamma}$. Namely on $\T^{n-1}\times\{x_n=0\}$
we have
\[
[\partial_nv]^+_-|_{\Gamma}=[(\nabla_{x'}u,0)-\partial_nu(\nabla\rho,1)]^+_-\cdot
\frac{(-\nabla\rho(x',t),1)}{\sqrt{1+|\nabla\rho|^2}}=[\partial_nu]^+_-\sqrt{1+|\nabla\rho|^2}.
\]
It is thus easy to see that the equation~(\ref{eq:N}) transforms into
\[
\partial_t\rho=(1+|\nabla\rho|^2)[\partial_nu]^-_+.
\]
For the sake of notational simplicity we also set
\[
\s:=\sqrt{1+|\nabla\rho|^2}.
\]
The Stefan problem~(\ref{eq:T}) -~(\ref{eq:rhoin}) now takes the following form:
\bea\label{eq:temp}
&&u_t-\Delta_{x'}u-\a u_{nn}+\B\cdot\nabla_{x'}u_n+\c u_n=0\\
\label{eq:gt}
&&u=\kappa\qquad\qquad\mbox{on}\ \T^{n-1}\times\{x_n=0\}\\
\label{eq:neumann}
&&\partial_nu=0\qquad\quad\mbox{on}\ \T^{n-1}\times\{x_n=\pm1\}\\
\label{eq:inu}
&&u(x,0)=u_0(x)\qquad x\in\Omega\\
\label{eq:inrho}
&&\rho(x',0)=\rho_0(x')\qquad x'\in\T^{n-1}\\
\label{eq:jump}
&&\rho_t=\s^2[u_n]^-_+\qquad\mbox{on}\ \T^{n-1}\times\{x_n=0\},
\eea
where we set $\c:=\g+\e$ with $\g$ and $\e$ given by~(\ref{eq:dro}) and~(\ref{eq:ero}) respectively.
Recall that $\a$ and $\B$ are given by~(\ref{eq:abro}).
In order to deal with the hyperbolic equation~(\ref{eq:jump}) we introduce
the regularization of the jump relation (\ref{eq:jump}):
\be\label{eq:jumpr}
\rho_t+\epsilon\Delta^2\rho_t=\s^2[u_n]^-_+\qquad\quad \mbox{on}\ \T^{n-1}\times\{x_n=0\}.
\ee
We shall refer to the problem of finding the solution to
(\ref{eq:temp})-(\ref{eq:inrho}) and (\ref{eq:jumpr}) as
to the regularized Stefan problem.\\
{\em Notation.}
For notational simplicity we define for any multi-index $\mu=(\mu_1,\dots,\mu_{n-1})$ and $s\in\N$
\be\label{eq:differential}
\d=\partial_t^s\partial_{x_1}^{\mu_1}\dots\partial_{x_{n-1}}^{\mu_{n-1}}.
\ee
Note that the operator $\d$ acts only in directions tangential to the
boundary $\T^{n-1}\times\{x_n=0\}$.
The Latin letters are always used to refer to the differentiation with respect to the time variable $t$ and
Greek letters to refer to the differentiation with respect to the first $n-1$ spatial variables $x_1,\dots,x_{n-1}$.
If each component of $\mu'$ is not greater then that of $\mu$ and $s'\leq s$, we write $(\mu',s')\leq (\mu,s)$. We write
$(\mu',s')<(\mu,s)$ if $(\mu',s')\leq(\mu,s)$ and $|\mu'|<|\mu|$ or $s'<s$. We also denote
$
C^{\mu'}_{s'}={(\mu,s)\choose (\mu',s')}.
$
For given functions $\omega:\T^{n-1}\to\R$ and $\U:\Omega\to\R$, we denote $\omega_i=\partial_{x_i}\omega$, $i=1,\dots, n-1$
and $\U_n=\partial_{x_n}\U$. With $x=(x_1,\dots,x_n)$ and $x'=(x_1,\dots,x_{n-1})$ we set
\[
\nabla_{x'}\U=(\partial_{x_1}\U,\dots,\partial_{x_{n-1}}\U),\quad
|\nabla_{x'}^2\U|^2=\sum_{i,j=1}^{n-1}(\partial_{x_ix_j}\U)^2,\quad
\Delta_{x'}\U=\sum_{i=1}^{n-1}\partial_{x_ix_i}\U.
\]
The Einstein summation convention is used throughout the paper when
dealing with repeated indices. The letter $C$ will stand for a generic
constant that may change from line to line.
\par
We define the following high-order energy norms:
\be
\begin{array}{l}
\displaystyle
\E(\U,\omega;\psi)(t):=\sum_{|\mu|+2s\leq 2k}\int_{\Omega}\Big\{(\d\U(t))^2+|\d \nabla_{x'}\U(t)|^2+a_{\psi(t)}(\d\U_n(t))^2\Big\}+\\
\displaystyle
+\sum_{|\mu|+2s\leq 2k}\int_{\T^{n-1}}\Big\{|\d\nabla\omega(t)|^2\left\langle\psi(t)\right\rangle^{-1}
+I_{\psi(t)}(\nabla^2\d\omega,\nabla^2\d\omega)(t)\Big\},
\end{array}
\label{eq:instantgen}
\ee
\be
\begin{array}{l}
\displaystyle
\D(\U,\omega;\psi)(t):=\sum_{|\mu|+2s\leq 2k}\int_{\Omega}\Big\{(\d\U_t(t))^2+|\d\nabla_{x'}\U(t)|^2
+a_{\psi(t)}(\d\U_n(t))^2+\\
\displaystyle
|\nabla_{x'}^2\d\U(t)|^2+2a_{\psi(t)}|\d\nabla_{x'}\U_n(t)|^2
+\big(a_{\psi(t)}\d\U_{nn}(t)\big)^2\Big\}\\
\displaystyle
+2\sum_{|\mu|+2s\leq 2k}\int_{\T^{n-1}}|\d\nabla\omega_t(t)|^2\left\langle\psi(t)\right\rangle^{-1},
\end{array}
\label{eq:dissipationgen}
\ee
where for given functions $\omega,\psi:\T^{n-1}\to\R$, we define
\be\label{eq:ipsi}
I_{\psi}(\nabla^2\omega,\nabla^2\omega):=|\nabla^2\omega|^2\p^{-1}
-\sum_{k=1}^{n-1}\big(\nabla\omega_k\cdot\nabla\psi\big)^2\p^{-3}.
\ee
Recall that $a_{\psi(t)}$ is given by~(\ref{eq:abro}).
It is crucial to observe that $I_{\psi}$ is
a positive definite bilinear form. Namely,
\be
\begin{array}{l}
\displaystyle
I_{\psi}(\nabla^2\omega,\nabla^2\omega)=
\int_{\T^{n-1}}\Big\{|\nabla^2\omega|^2\p^{-3}+
\Big(|\nabla^2\omega|^2|\nabla\psi|^2-
\sum_{k=1}^{n-1}\big(\nabla\omega_k\cdot\nabla\psi\big)^2\Big)\p^{-3}\Big\}\\
\displaystyle
\quad\quad\geq\int_{\T^{n-1}}|\nabla^2\omega|^2\p^{-3}.
\end{array}
\label{eq:posdef}
\ee
Note that we have used the Cauchy-Schwarz inequality in the last estimate.
The {\em instant energy} $\E$ and the {\em dissipation} $\D$ are respectively given by
\be\label{eq:instant}
\E\equiv \E(u,\rho):=\E(u,\rho;\rho),
\ee
\be\label{eq:dissipation}
\D\equiv\D(u,\rho):=\D(u,\rho;\rho).
\ee
In the rest of the paper we shall always assume $k\geq n$, where $n$ is the dimension of the space the domain
$\Omega$ belongs to.
Observe that the stationary solutions to the Stefan problem~(\ref{eq:temp}) -~(\ref{eq:jump}) are
given by $(u,\rho)\equiv(0,\bar{\rho})$, where $\bar{\rho}\in\R$ is a given constant.
Note that $\E(u,\rho-\bar{\rho})=\E(u,\rho)$.
The main result of the paper is the following theorem.
\begin{theorem}\label{th:global}
There exists a sufficiently small constant $M^{\ast}>0$ such that if the
initial data satisfy $(u_0,\rho_0)$ satisfy
\[
\E(u_0,\rho_0)+
\Big|\int_{\T^{n-1}}\rho_0-\int_{\Omega}u_0(1+\phi'\rho_0)\Big|
\leq M^{\ast},
\]
then there exists a unique global solution to the Stefan
problem~(\ref{eq:temp}) -~(\ref{eq:jump}) satisfying the global bound
\be\label{eq:star1}
\E(u,\rho)(t)+\frac{1}{2}\int_s^t\D(u,\rho)(\tau)\,d\tau\leq\E(u,\rho)(s),\quad t\geq s\geq0.
\ee
Moreover, given the stationary solution $(u,\rho)\equiv(0,\bar{\rho})$,
such that
\[
\int_{\T^{n-1}}\bar{\rho}=\int_{\T^{n-1}}\rho_0-\int_{\Omega}u_0(1+\phi'\rho_0),
\]
then for such small initial datum there exist constants $K_1,K_2>0$ such that
\[
\E(u,\rho)(t)+||\rho(t)-\bar{\rho}||_2^2\leq K_1e^{-K_2t}\quad\mbox{for all}\quad t\geq0.
\]
\end{theorem}
The proof of this theorem will strongly rely on careful examination
of the regularized Stefan problem~(\ref{eq:temp}) - (\ref{eq:inrho}) and~(\ref{eq:jumpr}).
For this purpose we introduce the appropriate energy norms incorporating the additional
viscosity coefficient $\epsilon$.
\be
\begin{array}{l}
\displaystyle
\E_{\epsilon}(\U,\omega;\psi):=\E(\U,\omega;\psi)
+\sum_{|\mu|+2s\leq 2k}\epsilon\int_{\T^{n-1}}
\Big\{|\d\Delta\nabla\omega|^2\left\langle\psi\right\rangle^{-1}\\
\displaystyle
\quad\qquad+I_{\psi}(\nabla^2\d\Delta\omega,\nabla^2\d\Delta\omega)\Big\}
\end{array}
\label{eq:energyeps}
\ee
\be\label{eq:dissipationeps}
\D_{\epsilon}(\U,\omega;\psi):=\D(\U,\omega;\psi)
+2\epsilon\sum_{|\mu|+2s\leq 2k}\int_{\T^{n-1}}|\nabla^2\d\Delta\omega|^2
\left\langle\psi\right\rangle^{-1}.
\ee
The above norms are the weighted versions of parabolic Sobolev norms given by
\be\label{eq:equivE}
||(\U,\omega)||_{\E_{\epsilon}}:=\sum_{|\mu|+2s\leq 2k}\Big\{||\d\U||_{H^1(\Omega)}^2+||\d\nabla\omega||_{H^1}^2+
\epsilon||\d\Delta\nabla\omega||_{H^1}^2\Big\}
\ee
and
\be\label{eq:equivD}
||(\U,\omega)||_{\D_{\epsilon}}:=\sum_{|\mu|+2s\leq 2k}\Big\{||\d\U_t||_{L^2(\Omega)}^2+
||\d\nabla\U||_{H^1(\Omega)}^2+||\d\nabla\omega_t||_2^2+\epsilon||\d\Delta\nabla\omega_t||_2^2\Big\}.
\ee
Given $||\psi||_{\infty}$ small enough so that
$(1+\phi'\psi)^2\geq\delta>0$ and $||\nabla\psi||_{\infty}$ bounded, we conclude
that there exists $C>0$ so that
\[
\frac{1}{C}||(\U,\omega)||_{\E_{\epsilon}}\leq\E_{\epsilon}\leq C||(\U,\omega)||_{\E_{\epsilon}},\quad
\frac{1}{C}||(\U,\omega)||_{\D_{\epsilon}}\leq \D_{\epsilon}\leq C||(\U,\omega)||_{\D_{\epsilon}}.
\]
In this sense the above norms are equivalent and this observation will be often implicitly used
throughout the paper.
The major part of the analysis will be concerned with proving the following result,
which states that the regularized Stefan problem has unique global solutions with
small initial data - independent of $\epsilon$.
\begin{theorem}\label{th:globalr}
There exists a sufficiently small constant $M>0$ independent of $\epsilon$,
such that the following
statement holds:
if for given initial data $(u^{\epsilon}_0,\rho^{\epsilon}_0)$ the inequality
\[
\E_{\epsilon}(u^{\epsilon}_0,\rho^{\epsilon}_0)+
\Big|\int_{\T^{n-1}}\rho^{\epsilon}_0-\int_{\Omega}u^{\epsilon}_0(1+\phi'\rho^{\epsilon}_0)\Big|
\leq M
\]
holds, then there exists a unique global solution $(u^{\epsilon},\rho^{\epsilon})$ to the regularized Stefan
problem~(\ref{eq:temp})- (\ref{eq:inrho}) and~(\ref{eq:jumpr}). Moreover,
\be\label{eq:star}
\E_{\epsilon}(u,\rho)(t)+\frac{1}{2}\int_0^t\D_{\epsilon}(u,\rho)(\tau)\,d\tau
\leq\E_{\epsilon}(u_0,\rho_0),\quad t\geq0.
\ee
\end{theorem}
The Stefan problem has been studied in a variety of mathematical literature over the past century
(see for instance \cite{Vi}).
It has been known that classical Stefan problem admits unique global classical solutions in $\R^1$
(\cite{Fri1}, \cite{Fri2} and \cite{Ka}).
Local classical solutions are established in~\cite{Ha} and~\cite{Me}.
\par
If the diffusion equation~(\ref{eq:T}) is
replaced by the elliptic equation $\Delta v=0$, then
the resulting problem is called the Hele-Shaw problem (or the quasi-stationary Stefan problem)
with surface tension. Global solutions for the
Hele-Shaw problem
in two dimensions with small initial data have been established in~\cite{CoPu}.
In~\cite{Ch}, stability of the solutions close to the steady state sphere
is established. Global stability for the one-phase Hele-Shaw problem is established
in~\cite{FriRe2}. Local-in-time
solutions in parabolic H\"older spaces in arbitrary dimensions are established in~\cite{ChHoYi}.

As to the Stefan problem with surface tension, global weak existence
theory (without uniqueness) is analyzed in \cite{Lu} and~\cite{Roe}.
An existence theory is also developed in~\cite{AlWa}. In
\cite{FriRe} the authors consider the Stefan problem with small
surface tension i.e. $\sigma\ll1$ if~(\ref{eq:D}) is substituted by
$v=\sigma\kappa(t)$. The local existence result for the Stefan
problem is studied in \cite{Ra}. In \cite{EsPrSi} the authors prove
a local existence and uniqueness result in suitable Besov spaces,
relying on the $L^p$-regularity theory.
\par
We establish a {\em global-in-time} existence, uniqueness and
exponential decay of classical solutions to the Stefan problem with
surface tension near a flat steady state (Theorem~\ref{th:global}).
The major difficulty consists of proving Theorem~\ref{th:globalr}
which establishes the existence and uniqueness result for the
regularized Stefan problem with the energy estimate
\be\label{eq:main}
\E_{\epsilon}(t)+\int_0^t\D_{\epsilon}(\tau)\,d\tau\leq\E_{\epsilon}(0)+
C\int_0^t\sqrt{\E_{\epsilon}(\tau)}\D_{\epsilon}(\tau)\,d\tau, \ee
where $C$ does not depend on $\epsilon$. Combined with the smallness
assumption on the instant energy $\E_{\epsilon}$ the
estimate~(\ref{eq:main}) gives~(\ref{eq:star}) and the
global-in-time existence. For a fixed $\epsilon$ we first construct
local-in-time solution for the regularized Stefan problem
~(\ref{eq:temp})- (\ref{eq:inrho}) and~(\ref{eq:jumpr}). The crux of
our method is the use of high order energy estimates, for the
differential operator $\d$ acts only in tangential directions with
respect to the boundary $\T^{n-1}\times\{x_n=0\}$. This is very
convenient when deriving the energy identities because the Neumann
boundary operator commutes with $\d$. The diffusion
equation~(\ref{eq:temp}) is then used to control high order
derivatives of $u$ with respect to the normal direction $x_n$, as it
is presented in Lemmas~\ref{lm:auximp} and~\ref{lm:auximp1}. We set
up an iteration scheme, which generates a sequence of iterates
$\{(u^m,\rho^m)\}_{m\in\N}$. Such iteration is well defined, but it
breaks the natural energy setting due to lack of exact cancelations
in the presence of the cross-terms. With fixed $\epsilon$, we
crucially use the regularization to prove that
$\{(u^m,\rho^m)\}_{m\in\N}$ is a Cauchy sequence in the energy
space. As $m\to\infty$ the unpleasant cross-terms disappear and we
recover~(\ref{eq:main}) in the limit. We conclude the proof of
Theorem~\ref{th:global} by letting $\epsilon\to0$.
\par
This work is the first step in our program of developing a robust
energy method to investigate
and characterize morphological stabilities/instabilities arising in numerous free boundary problems
in applied PDE. In particular, in a forth-coming paper we are going to establish stability
and instability(!) of steady spheres in the Stefan problem with surface tension.

The article is organized as follows: In Chapter~\ref{ch:basic}
we derive general energy identities for a
model Stefan problem. In Chapter~\ref{ch:local} the iteration
scheme for proving the local existence is set up and
the actual energy identities are derived, based on Chapter~\ref{ch:basic}.
Furthermore, some basic estimates
are established, which are then used in Chapter~\ref{ch:energye}
to prove the crucial energy estimates.
Chapter~\ref{ch:localwp}
is entirely devoted to the proof of the local-in-time existence nd uniqueness.
The main results, Theorems~\ref{th:global} and~\ref{th:globalr}
are proved in Chapter~\ref{ch:main}.
\section{Energy identities}\label{ch:basic}
Let $I=[0,q]$ for some $0<q\leq\infty$.
The derivation of the energy identities crucially depends on the following model
problem:
\bea
\label{eq:model1}
&&\U_t-\Delta_{x'}\U-\aa\U_{nn}=f\quad\textrm{on}\,\,\,\,\,\Omega\times I\\
\label{eq:model2}
&&\U=\Delta\chi\p^{-1}-\psi_i\psi_j\chi_{ij}\p^{-3}+G\quad\textrm{on}\,\,\,\,\,\T^{n-1}\times\{x_n=0\}\times I\\
&&\partial_n\U=0\quad\textrm{on}\,\,\,\,\,\T^{n-1}\times\{x_n=\pm1\}\times I\\
&&[\U_n]^-_+=\big(\omega_t+\epsilon\Delta^2\omega\big)\p^{-2}+h\quad\textrm{on}\,\,\,\,\,\T^{n-1}\times\{x_n=0\}\times I
\label{eq:model3}
\eea
We shall denote
\[
g:=-\psi_i\psi_j\chi_{ij}\p^{-3}+G.
\]
For most of the identities we shall derive, only the leading term
$\Delta\chi\p^{-1}$ in~(\ref{eq:model2}) will be relevant. We can thus write the equation~(\ref{eq:model2}) in the alternative
form
\be\label{eq:model2'}
\U=\Delta\chi\p^{-1}+g\quad\textrm{on}\,\,\,\,\,\T^{n-1}\times\{x_n=0\}\times I.
\ee
We define the energies $\bar{\E}_{\epsilon}$ and $\bar{\D}_{\epsilon}$ (for the model problem) by setting $k=0$ in the
definitions~(\ref{eq:energyeps}) and~(\ref{eq:dissipationeps}) of $\E_{\epsilon}$ and $\D_{\epsilon}$, respectively.
\begin{lemma}\label{lm:model}
Let each of the functions $\U$, $\omega$, $\chi$ and $\psi$ be five times continuously differentiable
with respect to the space variable and each of its spatial partial derivatives of order $\leq5$ once continuously
differentiable with respect to the time variable.
The following identity holds:
\be\label{eq:mainenergy}
\frac{d}{dt}\bar{\E}_{\epsilon}(\U,\omega;\psi)+\bar{\D}_{\epsilon}(\U,\omega;\psi)=
\int_{\Omega}\{P+R\}-\int_{\T^{n-1}}\{Q+S+T\},
\ee
where
\be
\label{eq:P}
P\equiv P(\U,\psi,f):=f\U-
(\aa)_n\U_n\U,\qquad\qquad\qquad\qquad\qquad\qquad\qquad\qquad\qquad
\ee
\be
\begin{array}{l}
\displaystyle
R\equiv R(\U,\psi,f):=f^2-2f\big(B_{\psi}\cdot\nabla_{x'}\U_n\big)
+\U_n^2(\aa)_t-2\U_t\U_n(\aa)_t\qquad\qquad\qquad\quad\\
\displaystyle
-2\nabla_{x'}\U_n\cdot\nabla_{x'}(\aa)\U_n+2\Delta_{x'}\U\U_n(\aa)_n,
\end{array}
\label{eq:R}
\ee
\be
\begin{array}{l}
\displaystyle
Q\equiv Q(\chi,\omega,\psi,g,h)=
\nabla\omega_t\cdot\big(\nabla\chi-\nabla\omega\big)\p^{-1}+
\epsilon\Delta\nabla\omega_t\cdot\big(\Delta\nabla\chi-\Delta\nabla\omega\big)\p^{-1}\qquad\\
\displaystyle
-\frac{1}{2}\Big\{|\nabla\omega|^2\p^{-1}_t+\epsilon|\Delta\nabla\omega|^2\p^{-1}_t\Big\}+
\omega_t\nabla\chi\cdot\nabla(\p^{-1})\\
\displaystyle
+\epsilon\Delta\nabla\omega_t
\cdot\nabla(\p^{-1})\Delta\chi
-(\omega_t+\epsilon\Delta^2\omega_t)g
-\p^2h\U|_{\T^{n-1}},
\end{array}
\label{eq:Q}
\ee
\be
\begin{array}{l}
\displaystyle
S\equiv S(\chi,\omega,\psi,g,h):=2\nabla\omega_t\cdot\big(\nabla\chi_t-\nabla\omega_t\big)\p^{-1}+
2\epsilon\Delta\nabla\omega_t\cdot\big(\Delta\nabla\chi_t-\Delta\nabla\omega_t\big)\p^{-1}\\
\displaystyle
+2\nabla\chi_t\cdot\nabla(\p^{-1})\omega_t+
2\epsilon\Delta\nabla\omega_t\cdot\nabla(\p^{-1})\Delta\chi_t
-2\big(\omega_t+\epsilon\Delta^2\omega_t\big)\big(\Delta\chi\p^{-1}_t
+g_t\big)\\
\displaystyle
-2h\p^2\U_t|_{\T^{n-1}},
\end{array}
\label{eq:S}
\ee
\be
\begin{array}{l}
\displaystyle
T\equiv T(\chi,\omega,\psi,G,h)=A(\chi,\omega,\psi,g,h)+B(\chi,\omega,\psi,g,h)+
2\Delta G\big(\omega_t+\epsilon\Delta^2\omega_t\big)\qquad\\
\displaystyle
\quad+2\Delta\U|_{\T^{n-1}}h\p^2.
\end{array}
\label{eq:bigT}
\ee
Here, the functions $A$ and $B$ are given by:
\be
\begin{array}{l}
\displaystyle
A:=2\Delta\omega_t
\big(\Delta\chi-\Delta\omega\big)\p^{-1}-
2\Delta\omega_t
\psi_i\psi_j(\chi_{ij}-\omega_{ij})\p^{-3}
-|\nabla^2\omega|^2\p^{-1}_t\\
\displaystyle
+2\omega_{it}
\nabla\omega_i\cdot\nabla(\p^{-1})
-2\nabla\omega_t\cdot\nabla(\p^{-1})\Delta\omega
+\omega_{jk}\omega_{ik}\psi_j\psi_{it}\p^{-3}
+\big(\nabla\omega_k\cdot\nabla\psi\big)^2\p^{-3}_t\\
\displaystyle
-2\omega_{jk}\Big\{\big[\psi_i\psi_j\omega_{kt}\p^{-3}\big]_i-
\psi_i\psi_j\omega_{ikt}\p^{-3}\Big\}
+2\omega_{kt}\Big\{\big[\psi_i\psi_j\omega_{ij}\p^{-3}\big]_k-
\psi_i\psi_j\omega_{ijk}\p^{-3}\Big\}
\end{array}
\label{eq:A}
\ee
and
\be
\begin{array}{l}
\displaystyle
B:=2\epsilon\Delta^2\omega_t
\big(\Delta^2\chi-\Delta^2\omega\big)\p^{-1}
-2\epsilon\Delta^2\omega_t
\psi_i\psi_j(\Delta\chi_{ij}-\Delta\omega_{ij})\p^{-3}\\
\displaystyle
-\epsilon|\nabla^2\Delta\omega|^2\p^{-1}_t
+2\epsilon\Delta\omega_{it}
\nabla\Delta\omega_i\cdot\nabla(\p^{-1})\\
\displaystyle
-2\epsilon\Delta\nabla\omega_t\cdot\nabla(\p^{-1})\Delta^2\omega
+\epsilon\Delta\omega_{jk}\Delta\omega_{ik}\psi_j\psi_{it}\p^{-3}
+\epsilon\big(\Delta\nabla\omega_k\cdot\nabla\psi\big)^2\p^{-3}_t\\
\displaystyle
-2\epsilon\Delta\omega_{jk}\Big\{\big[\psi_i\psi_j\Delta\omega_{kt}\p^{-3}\big]_i-
\psi_i\psi_j\Delta\omega_{ikt}\p^{-3}\Big\}\\
\displaystyle
+2\epsilon\Delta\omega_{kt}\Big\{\big[\psi_i\psi_j\Delta\omega_{ij}\p^{-3}\big]_k-
\psi_i\psi_j\Delta\omega_{ijk}\p^{-3}\Big\}\\
\displaystyle
+2\epsilon\Big\{\Delta\big(\Delta\chi\p^{-1}-\psi_i\psi_j\chi_{ij}\p^{-3}\big)-
\big(\Delta^2\chi\p^{-1}-\psi_i\psi_j\Delta\chi_{ij}\p^{-3}\big)\Big\}\Delta^2\omega_t.
\end{array}
\label{eq:B}
\ee
\end{lemma}
\prf
We start by multiplying the equation~(\ref{eq:model1}) by $\U$ and integrating
over $\Omega$. By a direct computation,
\be\label{eq:firstidentity}
\frac{1}{2}\partial_t\int_{\Omega}\U^2+\int_{\Omega}\Big\{|\nabla_{x'}\U|^2+
\aa(\U_n)^2\Big\}-\int_{\T^{n-1}}\p^2[\U_n]^-_+\U=\int_{\Omega}P(\U,\psi,f),
\ee
where $P(\U,\psi,f)$ is given by~(\ref{eq:P}).
Using the boundary conditions~(\ref{eq:model2'}) and~(\ref{eq:model3}), we obtain
\be
\begin{array}{l}
-\p^2[\U_n]^-_+\U=
-\big\{\omega_t\Delta\chi\p^{-1}+\epsilon\Delta^2\omega_t\Delta\chi\p^{-1}\big\}\\
-\p^2h\big[\Delta\chi\p^{-1}+g\big]-(\omega_t+\epsilon\Delta^2\omega_t)g.
\end{array}
\label{eq:hilfe1}
\ee
Integrating by parts in the first term on the RHS of~(\ref{eq:hilfe1}) above we arrive at
\be
\begin{array}{l}
-\int_{\T^{n-1}}\Big\{\omega_t\Delta\chi\p^{-1}+\epsilon\Delta^2\omega_t\Delta\chi\p^{-1}\Big\}=\\
\int_{\T^{n-1}}\nabla(\omega_t\p^{-1})\cdot\nabla\chi+\epsilon\Delta\nabla\omega_t\cdot\nabla(\Delta\chi\p^{-1}).
\end{array}
\label{eq:hilfe2}
\ee
By the product rule, the integrand on the RHS of~(\ref{eq:hilfe2}) can be written as
\[
\nabla\omega_t\p^{-1}\cdot\nabla\chi+\epsilon\Delta\nabla\omega_t\cdot\Delta\nabla\chi\p^{-1}+
\omega_t\nabla(\p^{-1})\cdot\nabla\chi+\epsilon\Delta\omega_t\nabla(\p^{-1})\cdot\Delta\nabla\chi.
\]
In each of the terms $\nabla\omega\p^{-1}\cdot\nabla\chi$ and $\Delta\nabla\omega_t\cdot\Delta\nabla\chi\p^{-1}$ we
set $\chi=\omega+(\chi-\omega)$ to obtain
\be
\begin{array}{l}
\displaystyle
\frac{1}{2}\partial_t\Big\{|\nabla\omega|^2\p^{-1}+\epsilon|\Delta\nabla\omega|^2\p^{-1}\Big\}
+\nabla\omega_t\cdot\big(\nabla\chi-\nabla\omega\big)\p^{-1}+\\
\displaystyle
+\epsilon\Delta\nabla\omega_t\cdot\big(\Delta\nabla\chi-\Delta\nabla\omega\big)\p^{-1}
-\frac{1}{2}\Big\{|\nabla\omega|^2\p^{-1}_t+\epsilon|\Delta\nabla\omega|^2\p^{-1}_t\Big\}+\\
\displaystyle
+\omega_t\nabla\chi\cdot\nabla(\p)^{-1}+\epsilon\Delta\nabla\omega_t
\cdot\nabla(\p^{-1})\Delta\chi
\end{array}
\label{eq:chitrick}
\ee
Thus plugging~(\ref{eq:chitrick}) into~(\ref{eq:hilfe2}) yields
\be\label{eq:firstboundary}
-\int_{\T^{n-1}}\p^2[\U_n]^-_+\U=
\frac{1}{2}\partial_t\int_{\T^{n-1}}\Big\{|\nabla\omega|^2\p^{-1}+\epsilon|\Delta\nabla\omega|^2\p^{-1}\Big\}
+\int_{\T^{n-1}}Q(\chi,\omega,\psi,g,h),
\ee
where $Q(\chi,\omega,\psi,g,h)$ is given by~(\ref{eq:Q}).
To complete the derivation of~(\ref{eq:mainenergy}), we take the square of the equation~(\ref{eq:model1})
and integrate over $\Omega$:
\be
\begin{array}{l}
\displaystyle
\partial_t\Big\{\int_{\Omega}|\nabla_{x'}\U|^2+\aa\U_n^2\Big\}+
\int_{\Omega}\Big\{\U_t^2+|\nabla_{x'}^2\U|^2+2a_{\psi}|\nabla_{x'}\U_n|^2
+\aa^2\U_{nn}^2\Big\}\\
\displaystyle
-2\int_{\T^{n-1}}\p^2[\U_n]^-_+\U_t+2\int_{\T^{n-1}}\p^2[\U_n]^-_+\Delta_{x'}\U=\int_{\Omega}R(\U,\psi,f),
\end{array}
\label{eq:omega2}
\ee
where $R(\U,\psi,f)$ is given by~(\ref{eq:R}).
The goal is to evaluate the
two integrals over $\T^{n-1}$ on LHS of~(\ref{eq:omega2})
using the boundary conditions~(\ref{eq:model2'}) and~(\ref{eq:model3}).
We first treat the integral $\int_{\T^{n-1}}\p^2[\U_n]^-_+\U_t$. Integrating by
parts in the leading order term, we obtain
\be
\begin{array}{l}
\displaystyle
-2\int_{\T^{n-1}}\p^2[\U_n]^-_+\U_t=
2\int_{\T^{n-1}}\Big\{|\nabla\omega_t|^2\p^{-1}+\epsilon|\Delta\nabla\omega_t|^2\p^{-1}\Big\}\\
\displaystyle
\quad\quad+\int_{\T^{n-1}}
S(\chi,\omega,\psi,g,h),
\end{array}
\label{eq:secondboundary}
\ee
where $S(\chi,\omega,\psi,g,h)$ is given by~(\ref{eq:S}).
Note that the expression~(\ref{eq:S}) is obtained similarly to~(\ref{eq:chitrick}) by setting
$\chi=\omega+(\chi-\omega)$ in the leading order terms.
\par
The second integral over $\T^{n-1}$ in the identity~(\ref{eq:omega2}) is more delicate.
We shall make use of the boundary conditions~(\ref{eq:model2}) and~(\ref{eq:model3}) to evaluate it.
The relation~(\ref{eq:model2}) is used is to exploit the full algebraic structure of the
curvature-type term $\Delta\chi\p^{-1}-\psi_i\psi_j\chi_{ij}\p^{-3}$, which is important
in the energy estimates later on.
We have:
\be
\begin{array}{l}
\displaystyle
\p^2[\U_n]^-_+\Delta_{x'}\U=
\Delta\big(\Delta\chi\p^{-1}-\psi_i\psi_j\chi_{ij}\p^{-3}\big)(\omega_t+\epsilon\Delta^2\omega_t)\\
\displaystyle
+\Delta G\big(\omega_t+\epsilon\Delta^2\omega_t\big)
+\p^2h\Delta_{x'}\U.
\end{array}
\label{eq:hilfe3}
\ee
Observe that
\be
\begin{array}{l}
\displaystyle
2\int_{\T^{n-1}}\Delta\big(\Delta\chi\p^{-1}-\psi_i\psi_j\chi_{ij}\p^{-3}\big)\omega_t=
2\int_{\T^{n-1}}\big(\Delta\chi\p^{-1}-\psi_i\psi_j\chi_{ij}\p^{-3}\big)\Delta\omega_t=\\
\displaystyle
\partial_t\int_{\T^{n-1}}|\nabla^2\omega|^2\p^{-1}-\int_{\T^{n-1}}|\nabla^2\omega|^2\p^{-1}_t
+2\int_{\T^{n-1}}\omega_{it}\nabla\omega_i\cdot\nabla(\p^{-1})-\\
\displaystyle
2\int_{\T^{n-1}}\nabla\omega_t\cdot\nabla(\p^{-1})\Delta\omega
+2\int_{\T^{n-1}}\Delta\omega_t
\big(\Delta\chi-\Delta\omega\big)\p^{-1}-\\
\displaystyle
2\int_{\T^{n-1}}\omega_{kkt}\psi_i\psi_j\omega_{ij}\p^{-3}-
2\int_{\T^{n-1}}\Delta\omega_t\big(\psi_i\psi_j(\chi_{ij}-\omega_{ij})\p^{-3}\big).
\end{array}
\label{eq:hilfe4}
\ee
Here, just like in~(\ref{eq:chitrick}) we substituted $\chi=\omega+(\chi-\omega)$ in the
leading order terms.
Note that we have repeatedly used integration by parts.
Integrating by parts twice, we obtain
\be
\begin{array}{l}
\displaystyle
-2\int_{\T^{n-1}}\omega_{kkt}\psi_i\psi_j\omega_{ij}\p^{-3}=
2\int_{\T^{n-1}}\omega_{kt}\psi_i\psi_j\omega_{ijk}\p^{-3}+\\
\displaystyle
2\int_{\T^{n-1}}\omega_{kt}\Big\{\big[\psi_i\psi_j\omega_{ij}\p^{-3}\big]_k-
\psi_i\psi_j\omega_{ijk}\p^{-3}\Big\}=\\
\displaystyle
-2\int_{\T^{n-1}}\omega_{ikt}\psi_i\psi_j\omega_{jk}\p^{-3}
-2\int_{\T^{n-1}}\omega_{jk}\Big\{\big[\psi_i\psi_j\omega_{kt}\p^{-3}\big]_i-
\psi_i\psi_j\omega_{ikt}\p^{-3}\Big\}\\
\displaystyle
+2\int_{\T^{n-1}}\omega_{kt}\Big\{\big[\psi_i\psi_j\omega_{ij}\p^{-3}\big]_k-
\psi_i\psi_j\omega_{ijk}\p^{-3}\Big\}.
\end{array}
\label{eq:hilfe5}
\ee
We now single out the $t$-derivative in the first term on RHS of~(\ref{eq:hilfe5}) to obtain
\be
\begin{array}{l}
\displaystyle
-2\omega_{ikt}\psi_i\psi_j\omega_{jk}\p^{-3}=
-\partial_t\Big\{\sum_{k=1}^{n-1}\big(\nabla\omega_k\cdot\nabla\psi\big)^2\p^{-3}\Big\}
+\omega_{jk}\omega_{ik}\psi_j\psi_{it}\p^{-3}\\
\displaystyle
\quad\quad+\big(\nabla\omega_k\cdot\nabla\psi\big)^2\p^{-3}_t.
\end{array}
\label{eq:hilfe6}
\ee
We combine the identities~(\ref{eq:hilfe4}), ~(\ref{eq:hilfe5}) and~(\ref{eq:hilfe6}) to conclude
\beas
&&2\int_{\T^{n-1}}\Delta\big(\Delta\chi\p^{-1}-\psi_i\psi_j\chi_{ij}\p^{-3}\big)\omega_t=\\
&&\partial_t\int_{\T^{n-1}}|\nabla^2\omega|^2\p^{-1}
-\partial_t\Big\{\sum_{k=1}^{n-1}\int_{\T^{n-1}}\big(\nabla\omega_k\cdot\nabla\psi\big)^2\p^{-3}\Big\}
+\int_{\T^{n-1}}A\\
&&=\partial_t\int_{\T^{n-1}}I_{\psi}(\nabla^2\omega,\nabla^2\omega)+\int_{\T^{n-1}}A,
\eeas
where $I_{\psi}$ and $A$ are given by~(\ref{eq:ipsi}) and~(\ref{eq:A}) respectively.
\par
In the $\epsilon$-dependent part on RHS of~(\ref{eq:hilfe3})
we set $\omega^{\ast}=\Delta\omega$ and $\chi^{\ast}=\Delta\chi$. We can write
\beas
&&\Delta\big(\Delta\chi\p^{-1}-\psi_i\psi_j\chi_{ij}\p^{-3}\big)\Delta^2\omega_t=
\big(\Delta\chi^{\ast}\p^{-1}-\psi_i\psi_j\chi^{\ast}_{ij}\psi^{-3}\big)\Delta\omega^{\ast}_t\\
&&+\Big\{\Delta\big(\chi^{\ast}\p^{-1}-\psi_i\psi_j\chi_{ij}\p^{-3}\big)-
\big(\Delta\chi^{\ast}\p^{-1}-\psi_i\psi_j\chi^{\ast}_{ij}\p^{-3}\big)\Big\}\Delta\omega^{\ast}_t.
\eeas
We may now apply the same computation as in~(\ref{eq:hilfe4}) to conclude
\beas
2\epsilon\int_{\T^{n-1}}\Delta\big(\Delta\chi\p^{-1}-\psi_i\psi_j\chi_{ij}\p^{-3}\big)\Delta^2\omega_t=
\partial_t\int_{\T^{n-1}}\epsilon I_{\psi}(\nabla^2\Delta\omega,\nabla^2\Delta\omega)
+\int_{\T^{n-1}}B,
\eeas
where $B$ is given by~(\ref{eq:B}).
We combine the above identities to write the final form of the
second integral over $\T^{n-1}$ in the identity~(\ref{eq:omega2}):
\be
\begin{array}{l}
\displaystyle
2\int_{\T^{n-1}}\p^2[\U_n]^-_+\Delta_{x'}\U=
\partial_t\Big\{\int_{\T^{n-1}}I_{\psi}(\nabla^2\omega,\nabla^2\omega)+
\epsilon I_{\psi}(\nabla^2\Delta\omega,\nabla^2\Delta\omega)\Big\}\\
\displaystyle
\quad\quad+\int_{\T^{n-1}}T(\chi,\omega,\psi,G,h),
\end{array}
\label{eq:thirdboundary}
\ee
where $T$ is given by~(\ref{eq:bigT}). By summing the
identities~(\ref{eq:firstidentity}) and~(\ref{eq:omega2}),
plugging~(\ref{eq:firstboundary}) in~(\ref{eq:firstidentity})
and~(\ref{eq:secondboundary}) and~(\ref{eq:thirdboundary}) into~(\ref{eq:omega2})
and collecting terms, we conclude the proof of the lemma.
\prfe
\section{Iteration scheme and the basic estimates}\label{ch:local}
We shall set up an iterative scheme in order to solve the regularized Stefan problem locally-in-time.
For given $\rho^m$ and Cauchy data $u^{\epsilon}_0\in C^{\infty}(\Omega)$, $\rho^{\epsilon}_0\in C^{\infty}(\T^{n-1})$,
we solve the following problem:
\bea
&&\label{eq:tempm}
u^{m+1}_t-\Delta_{x'}u^{m+1}-\am u^{m+1}_{nn}+\Bm\cdot\nabla_{x'}u^{m+1}_n+\cm u^{m+1}_n=0\\
&&\label{eq:gtm}
u^{m+1}=\kappa^m\qquad\quad    \mbox{on}\ \T^{n-1}\times\{x_n=0\}\\
&&\label{eq:neumannm}
\partial_nu^{m+1}=0\qquad\qquad \mbox{on}\ \T^{n-1}\times\{x_n=\pm1\}\\
&&\label{eq:initialm}
u^{m+1}(x,0)=u^{\epsilon}_0(x),\quad\rho^{m+1}(x',0)=\rho^{\epsilon}_0(x').
\eea
Here
\be\label{eq:kappam}
\kappa^m:=\nabla\cdot\big(\frac{\nabla\rho^m}{\sm}\big).
\ee
The solution to the problem~(\ref{eq:tempm}) -~(\ref{eq:initialm}) exists and is smooth
(see Chapter 4 of~\cite{Li}). Having
obtained $u^{m+1}$, we solve the equation
\be\label{eq:jumprm}
\rho^{m+1}_t+\epsilon\Delta^2\rho^{m+1}_t=\sm^2[u^{m+1}_n]^-_+\qquad\quad \mbox{on}\ \T^{n-1}\times\{x_n=0\}
\ee
for $\rho^{m+1}$. We aim for proving the convergence of the sequence $(u^m,\rho^m)$
to the solution of the regularized Stefan problem in
the energy space.
Applying the tangential differential operator $\d$ (recall~(\ref{eq:differential})) to the equations~(\ref{eq:tempm}), (\ref{eq:gtm})
and~(\ref{eq:jumprm}), we obtain
\[
\d u^{m+1}_t-\d\Delta_{x'}u^{m+1}-\am \d u^{m+1}_{nn}=f^m_{\mu,s}
\]
\beas
&&\d u^{m+1}=\d\kappa^m=\Delta\d\rho^m\sm^{-1}+g^m_{\mu,s}\\
&&\quad\quad=\Delta\d\rho^m\sm^{-1}-\rho^m_i\rho^m_j\Delta\d\rho_{ij}^m\sm^{-1}+G^m_{\mu,s}
\eeas
\[
[\d u^{m+1}_n]^-_+\sm^2=\d\rho^{m+1}_t+\epsilon\Delta^2\rho^{m+1}_t+h^m_{\mu,s}\sm^2,
\]
where
\be
\begin{array}{l}
\displaystyle
f^m_{\mu,s}=\Big\{\d\Big(\am u^{m+1}_{nn}\Big)-\am\d u^{m+1}_{nn}\Big\}-
\d\Big(\Bm\cdot\nabla_{x'}u^{m+1}_n\Big)
-\d\Big(\cm u^{m+1}_n\Big),
\end{array}
\label{eq:ef}
\ee
\be\label{eq:gsmall}
g^m_{\mu,s}=\sum_{|\mu'|+s'\atop<|\mu|+s}C^{\mu'}_{s'}
\partial^{\mu'}_{s'}\Delta\rho^m\partial^{\mu-\mu'}_{s-s'}(\sm^{-1})+
\d(\nabla\rho^m\cdot\nabla(\sm^{-1})),
\ee
\be
\begin{array}{l}
\displaystyle
G^m_{\mu,s}=\sum_{|\mu'|+s'\atop<|\mu|+s}C^{\mu'}_{s'}
\partial^{\mu'}_{s'}\Delta\rho^m\partial^{\mu-\mu'}_{s-s'}(\sm^{-1})-\\
\displaystyle
\quad-\big\{\d\big(\rho^m_i\rho^m_j\rho_{ij}^m\sm^{-1}\big)-
\rho^m_i\rho^m_j\d\rho_{ij}^m\sm^{-1}\big\},
\end{array}
\label{eq:gbig}
\ee
\be\label{eq:ha}
h^m_{\mu,s}=\sum_{|\mu'|+s'\atop<|\mu|+s}C^{\mu'}_{s'}\partial^{\mu'}_{s'}(\rho^{m+1}_t
+\epsilon\Delta^2\rho^{m+1}_t)\partial^{\mu-\mu'}_{s-s'}(\sm^{-2}).
\ee
For any $l\in\N$ let us define
\be\label{eq:energies}
\E^l:=\E_{\epsilon}(u^l,\rho^l;\rho^{l-1}),\quad\D^l:=\D_{\epsilon}(u^l,\rho^l;\rho^{l-1}),
\ee
where $\E_{\epsilon}$ and $\D_{\epsilon}$
are defined by~(\ref{eq:energyeps}) and~(\ref{eq:dissipationeps}) respectively.
Setting $\U=\d u^{m+1}$, $\omega=\d\rho^{m+1}$, $\chi=\d\rho^m$, $\psi=\rho^m$,
$f=f^m_{\mu,s}$, $g=g^m_{\mu,s}$, $G=G^m_{\mu,s}$ and $h=h^m_{\mu,s}$,
the identity~(\ref{eq:mainenergy}) implies
\be\label{eq:energymm}
\frac{d}{dt}\E^{m+1}(t)+\D^{m+1}(t)=\int_{\Omega}\big\{P^m+R^m\big\}
-\int_{\T^{n-1}}\big\{Q^m+S^m+T^m\big\}.
\ee
Here $P^m=\sum_{|\mu|+2s\leq2k}P^m_{\mu,s}$ and $R^m$, $Q^m$, $S^m$ and $T^m$ are defined analogously, whereby
\be\label{eq:PQ}
P^m_{\mu,s}=P(\d u^{m+1},\,\rho^m,f^m_{\mu,s}),\quad Q^m_{\mu,s}=Q(\d\rho^m,\,\d\rho^{m+1},\,\rho^m,\,g^m_{\mu,s},\,h^m_{\mu,s})
\ee
and
\be
\begin{array}{l}
\displaystyle
R^m_{\mu,s}:=R(\d u^{m+1},\,\rho^m,f^m_{\mu,s}),\quad
S^m_{\mu,s}:=S(\d\rho^m,\d\rho^{m+1},\rho^m,g^m_{\mu,s},h^m_{\mu,s}),\\
\displaystyle
T^m_{\mu,s}:=T(\d\rho^m,\d\rho^{m+1},\rho^m,g^m_{\mu,s},h^m_{\mu,s}).
\end{array}
\label{eq:RST}
\ee
The inequality~(\ref{eq:posdef}) implies that the instant energy $\E^l$ is positive definite.
In order to estimate $P^m$, $R^m$, $Q^m$, $S^m$ and $T^m$ we first need to establish some basic
auxiliary estimates.
\begin{lemma}
The following identity holds
\be\label{eq:consm}
\partial_t\Big\{\int_{\Omega} u^{m+1}(1+\phi'\rho^m)\Big\}=\partial_t\Big\{\int_{\T^{n-1}}\rho^{m+1}\Big\}.
\ee
\end{lemma}
\prf
We multiply the equation (\ref{eq:tempm}) with $(1+\phi'\rho^m)$ and integrate over
$\Omega$. We thus obtain
\beas
&&\int_{\Omega}(1+\phi'\rho^m)(u^{m+1}_t-\Delta_{x'}u^{m+1})-
\int_{\Omega}\frac{1+|\phi\nabla\rho^m|^2}{1+\phi'\rho^m}u^{m+1}_{nn}+
2\int_{\Omega}\phi\nabla\rho^m\cdot\nabla_{x'}u^{m+1}_n\\
&&+\int_{\Omega}\phi\Delta\rho^mu^{m+1}_n-
\int_{\Omega}\Big(\frac{(\phi^2)'|\nabla\rho^m|^2}{1+\phi'\rho^m}-
\frac{\phi''\rho^m(1+|\phi\nabla\rho^m|^2)}{(1+\phi'\rho^m)^2}\Big)u^{m+1}_n\\
&&-\int_{\Omega}\phi\rho^m_tu^{m+1}_n=0.
\eeas

Integrating by parts we have
$
\int_{\Omega}-\phi\rho^m_tu^{m+1}_n=
\int_{\Omega}\phi'\rho^m_tu^{m+1}.
$
Using this identity, we obtain
\[
\int_{\Omega}(1+\phi'\rho^m)u^{m+1}_t+\int_{\Omega}(1+\phi'\rho^m)\e u^{m+1}_n=
\partial_t\Big\{\int_{\Omega}(1+\phi'\rho^m)u^{m+1}\Big\}.
\]
Observe that the integration by parts implies
\[
\int_{\Omega}-(1+\phi'\rho^m)\Delta_{x'}u^{m+1}=
\int_{\Omega}\phi'\nabla_{x'}\rho^m\cdot\nabla_{x'}u^{m+1}
\]
and
$
\int_{\Omega}\phi\Delta\rho^mu^{m+1}_n=-\int_{\Omega}\phi\nabla\rho^m\cdot\nabla_{x'}u^{m+1}_n.
$
Thus
\[
-\int_{\Omega}(1+\phi'\rho^m)\Delta_{x'}u^{m+1}+\int_{\Omega}2\phi\nabla\rho^m\cdot\nabla_{x'}u^{m+1}_n+
\int_{\Omega}\phi\Delta\rho^mu^{m+1}_n=0
\]
Note further that
\[
\partial_n\Big(\frac{1+|\phi\nabla\rho^m|^2}{1+\phi'\rho^m}\Big)=
\frac{(\phi^2)'|\nabla\rho^m|^2}{1+\phi'\rho^m}-
\frac{\phi''\rho^m(1+|\phi\nabla\rho^m|^2)}{(1+\phi'\rho^m)^2}
\]
 Using integration by parts again we have
\beas
&&-\int_{\Omega}\frac{1+|\phi\nabla\rho^m|^2}{1+\phi'\rho^m}u^{m+1}_{nn}
-\int_{\Omega}\Big(\frac{(\phi^2)'|\nabla\rho^m|^2}{1+\phi'\rho^m}-
\frac{\phi''\rho^m(1+|\phi\nabla\rho^m|^2)}{(1+\phi'\rho^m)^2}\Big)u^{m+1}_n=\\
&&-\int_{\T^{n-1}}\sm^2[u^{m+1}_n]^-_+=-\partial_t\Big\{\int_{\T^{n-1}}\rho^{m+1}\Big\}.
\eeas
This finishes the proof of the lemma.
\prfe
\par
The importance of this identity is reflected in the fact that it allows to control terms
with purely temporal derivatives of $\rho^{m+1}$:
\begin{lemma}\label{lm:tau}
There exist positive constants $K$ and $\theta_0<1$ such that for any
$\theta\leq\theta_0$ such that if
\[
\Big|\int_{\T^{n-1}}\rho^{\epsilon}_0-\int_{\Omega}u^{\epsilon}_0(1+\phi'\rho^{\epsilon}_0)\Big|<\theta,
\]
\[
\E^m\leq\theta,\quad||\nabla\rho^{m-1}||_{\infty}\leq 1\quad\textrm{and}
\quad\sum_{p=0}^{k}||\partial_p\rho^m||_2\leq K\sqrt{\E^m}+\theta,
\]
then\quad
$
||\nabla\rho^m||_{\infty}\leq1\quad\textrm{and}
\quad\sum_{p=0}^{k}||\partial_p\rho^{m+1}||_2\leq K\sqrt{\E^{m+1}}+\theta.
$
\end{lemma}
\prf
Observe first that the assumption on $\rho^{m-1}$ implies that
$\left\langle \rho^{m-1}\right\rangle\leq\sqrt{2}$. Using the Sobolev inequality,
we obtain
\[
||\nabla\rho^m||_{\infty}\leq C_{\ast}||\nabla\rho^m||_{H^{(n+1)/2}}\leq C_{\ast}\sqrt{2}\E^m
\leq C_{\ast}\sqrt{2}\theta.
\]
Thus, choosing $\theta_0\leq\frac{1}{C_{\ast}\sqrt{2}}$ guarantees $||\nabla\rho^m||_{\infty}\leq1$.
By~(\ref{eq:consm}) we have
\[
\int_{\T^{n-1}}\partial_s\rho^{m+1}_t=\int_{\Omega}\partial_su^{m+1}_t+
\int_{\Omega}\phi'\sum_{p=0}^{s+1}{s+1\choose p}\partial_pu^{m+1}
\partial_{s+1-p}\rho^m.
\]
Also, $\int_{\T^{n-1}}\rho^{m+1}=\int_{\Omega}(u^{m+1}+\phi'u^{m+1}\rho^m)+
\int_{\T^{n-1}}\rho^{\epsilon}_0-\int_{\Omega}u^{\epsilon}_0(1+\phi'\rho^{\epsilon}_0).
$
Thus, for $0\leq s\leq k-1$, using the Cauchy-Schwarz inequality, definition~(\ref{eq:energies}) of $\E^{m+1}$ and
the main assumption in the statement of the lemma, we obtain
\beas
\Big|\int_{\T^{n-1}}\partial_s\rho^{m+1}_t\Big|&\leq&
||\partial_{s+1}u^{m+1}||_2+C\sum_{p=1}^{s+1}||\partial_pu^{m+1}||_{L^2(\Omega)}
||\partial_{s+1-p}\rho^m||_2\\
&\leq&\sqrt{\E^{m+1}}+C\sqrt{\E^{m+1}}\sum_{p=0}^{k}||\partial_p\rho^m||_2.
\eeas
Furthermore,
\beas
\Big|\int_{\T^{n-1}}\rho^{m+1}\Big|&\leq& C||u^{m+1}||_2+
C||u^{m+1}||_{L^2(\Omega)}||\rho^m||_2+\Big|\int_{\T^{n-1}}\rho^{\epsilon}_0-\int_{\Omega}u^{\epsilon}_0(1+\phi'\rho^{\epsilon}_0)\Big|\\
&\leq&C\sqrt{\E^{m+1}}+C\sqrt{\E^{m+1}}\sum_{p=0}^{k}||\partial_p\rho^m||_2+\theta.
\eeas
By the Poincar\'e inequality we get
\be\label{eq:poincare}
\sum_{p=0}^k||\partial_p\rho^{m+1}||_2\leq C\sum_{p=0}^k||\partial_p\nabla\rho^{m+1}||_2+
C\sum_{p=0}^k||\int_{\T^{n-1}}
\partial_p\rho^{m+1}||_2.
\ee
The first term on the right-hand side is estimated by $C\sqrt{\E^{m+1}}$,
by the definition~(\ref{eq:energies}) of $\E^{m+1}$.
By the previous
two inequalities and the assumptions of the lemma, we can estimate the second sum on RHS
of~(\ref{eq:poincare}) by
$C\sqrt{\E^{m+1}}+C\sqrt{\E^{m+1}}\big(K\sqrt{\E^m}+\theta\big)+\theta$.
Keeping in mind that $\E^m\leq\theta$, we choose $1\leq K$ large enough
and $\theta<\frac{1}{C_{\ast}\sqrt{2}}$ small enough so that
\[
\sum_{p=0}^k||\partial_p\rho^{m+1}||_2\leq K\sqrt{\E^{m+1}}+\theta.
\]
\prfe
In the following we shall work under the standing assumption
\be\label{eq:standing}
\E^m\leq\theta,\quad||\nabla\rho^{m-1}||_{\infty}\leq1\quad\sum_{p=0}^{k}||\partial_p\rho^m||_2\leq K\sqrt{\E^m}+\theta,
\ee
with $\theta\leq\theta_0$ where $\theta_0$ is given as in Lemma~\ref{lm:tau}.
\begin{lemma}\label{lm:aux}
Let $\xi\in C^{\infty}(J,\R)$ and $J\subseteq\R$ an interval such that
every derivative of $\xi$ is uniformly bounded on $J$.
\begin{enumerate}
\item[\rm{(a)}]
Let $(\mu,r)$ be
a pair of indices such that $|\mu|+2r\leq2k$ and $(|\mu|,r)\neq(0,0)$.
Then there exists a positive constant $C$ such that
\be\label{eq:aux1}
||\partial^{\mu}_r\big[\xi(|\nabla\rho^m|^2)\big]||_{H^1}\leq C\sqrt{\E^m}
\ee
and
\be\label{eq:aux11}
||\partial^{\mu}_r\big[\xi(|\nabla\rho^m|^2)\big]||_{H^3}\leq \frac{C}{\sqrt{\epsilon}}\sqrt{\E^m}.
\ee
\item[\rm{(b)}]
There exists a positive constant $C$ such that
\be\label{aux2}
||\partial^{\mu}_s\nabla\rho^m||_{H^2}\leq C\sqrt{\D^{m+1}},
\ee
where $|\mu|+2s\leq 2k$.
Furthermore,
\be\label{eq:aux3}
||\partial^{\mu}_{s+1}\xi(|\nabla\rho^m|^2)||_2\leq C\sqrt{\D^m}
\ee
for all $(\mu,s)$ satisfying $|\mu|+2s\leq 2k$.
\item[\rm{(c)}]
For any pair of indices $(\mu,s)$ such that $|\mu|+2s\leq 2k$ there exists
a positive constant C and a small parameter $\lambda$ such that
\be\label{eq:higher}
||\d\rho^{m+1}_t||_2^2+
2\epsilon||\d\Delta\rho^{m+1}_t||_2^2+
\epsilon^2||\d\Delta^2\rho^{m+1}_t||_2^2\leq
\frac{C}{\lambda}\E^{m+1}+C\lambda\D^{m+1}.
\ee
\end{enumerate}
\end{lemma}

\prf
{\em Part (a):}
Let $\alpha=\mu+\tau$ for any given multi-index of length $n-1$ satisfying $|\tau|\leq1$.
Let first $r=0$. By assumption $|\alpha|\geq1$. Using Moser's inequality (cf. \cite{GuHaSp}, Lemma 5) and Leibniz'
rule (cf. \cite{GuHaSp}, Lemma 4), we have
\beas
&&||\partial^{\alpha}\big[\xi(|\nabla\rho^m|^2)\big]||_2\leq
C\max_{1\leq d\leq|\alpha|}\Big(||\xi^{(d)}(|\nabla\rho^m|^2))||_{\infty}
||\nabla\rho^m||_{\infty}^{|\alpha|-1}\Big)||\partial^{\alpha}(|\nabla\rho^m|^2)||_2\\
&&\quad\quad\leq C||\partial^{\alpha}(|\nabla\rho^m|^2)||_2
\leq C||\partial^{\alpha}\nabla\rho^m||_2||\nabla\rho^m||_{\infty}
\leq C\sqrt{\E^m}.
\eeas
Let $r\geq1$.
\beas
&&\partial^{\alpha}_r\big[\xi(|\nabla\rho^m|^2)\big]=
\partial^{\alpha}\Big\{\sum_{d=1}^r\sum_{s_1+\dots+s_d=r \atop s_i>0}
C_d\xi^{(d)}(|\nabla\rho^m|^2)\partial_{s_1}(|\nabla\rho^m|^2)
\dots\partial_{s_d}(|\nabla\rho^m|^2)\Big\}\\
&&=\sum_{d=1}^r\sum_{s_1+\dots+s_d=r \atop s_i>0}
\sum_{\gamma_1+\dots+\gamma_{d+1}\atop=\alpha}
C_dC_{\gamma_1\dots\gamma_{d+1}}
\partial^{\gamma_1}_{s_1}(|\nabla\rho^m|^2)
\dots\partial^{\gamma_d}_{s_d}(|\nabla\rho^m|^2)\partial^{\gamma_{d+1}}
\big[\xi^{(d)}(|\nabla\rho^m|^2)\big].
\eeas
For any $i=1,\dots,d$, we have
$
\partial^{\gamma_i}_{s_i}(|\nabla\rho^m|^2)=
\sum_{l=0}^{s_i}\sum_{\delta\leq\gamma_i}C_{l,\delta}
\partial^{\delta}_l\nabla\rho^m\partial^{\gamma_i-\delta}_{s_i-l}\nabla\rho^m.
$
Thus if $|\gamma_i|+2s_i\leq k+1$, then by the Sobolev inequality
\[
||\partial^{\delta}_l\nabla\rho^m||_{\infty}\leq
||\partial^{\delta}_l\nabla\rho^m||_{H^{\frac{n+1}{2}}}\leq C\sqrt{\E^m},
\]
and analogously $||\partial^{\gamma_i-\delta}_{s_i-l}\nabla\rho^m||_{\infty}\leq C\sqrt{\E^m}$,
implying $||\partial^{\gamma_i}_{s_i}(|\nabla\rho^m|^2)||_{\infty}\leq C\E^m$.
If $|\gamma_{d+1}|\leq k+1$, we use the Sobolev and Moser's inequality to conclude
\[
||\partial^{\gamma_{d+1}}
\big[\xi^{(d)}(|\nabla\rho^m|^2)\big]||_{\infty}\leq C\sqrt{\E^m}.
\]
If there exists $1\leq j\leq d$ such that
$|\gamma_j|+2s_j>k+1$, then $|\gamma_i|+2s_i\leq k+1$ for $1\leq i\leq d$, $i\neq j$ and additionally,
$|\gamma_{d+1}|\leq k+1$. Thus we can estimate the term containing $\gamma_j$ in superscript in $L^2$-norm and
the remaining terms in $L^{\infty}$-norm. If on the other hand $|\gamma_i|+2s_i\leq k+1$ for every
$1\leq i\leq d$, we estimate the term $\partial^{\gamma_{d+1}}
\big[\xi^{(d)}(|\nabla\rho^m|^2)\big]$ in $L^2$-norm and the remaining terms in $L^{\infty}$-norm.
We conclude that
\[
||\partial^{\alpha}_r\big[\xi(|\nabla\rho^m|^2)\big]||_2\leq C\sqrt{\E^m},
\]
for the specified range of $\alpha$-s and $r$-s.
The inequality~(\ref{eq:aux11})is proved similarly.

{\em Part (b):}
By~(\ref{eq:gtm}),
$
\Delta\rho^m=u^{m+1}\sm+\rho^m_i\rho^m_j\rho^m_{ij}\sm^{-2}.
$
Let $\gamma=\mu+\tau$ where $|\tau|=1$.
Applying $\partial^{\gamma}_s$ to the above identity, we get
\beas
\partial^{\gamma}_s\Delta\rho^m&=&\sum_{\gamma',s'}C^{\gamma'}_{s'}
\partial^{\gamma'}_{s'}(u^{m+1})\partial^{\gamma-\gamma'}_{s-s'}(\sm)+\\
&&\sum_{\sum (\gamma_l+s_l)\atop=\gamma+s}C^{\gamma_1,\gamma_2,\gamma_3,\gamma_4}_{s_1,s_2,s_3,s_4}
\partial^{\gamma_1}_{s_1}\rho^m_i\partial^{\gamma_2}_{s_2}\rho^m_j\partial^{\gamma_3}_{s_3}\rho^m_{ij}
\partial^{\gamma_4}_{s_4}(\sm^{-2})
\eeas
Observe that $\int_{\T^{n-1}}u^{m+1}=0$ since $u^{m+1}=\kappa^m=\nabla\cdot\big(\nabla\rho^m\sm^{-1}\big)$ on $\T^{n-1}$.
Let us fix $(\gamma',s')\leq(\gamma,s)$. If $|\gamma'|>0$ note that $\int_{\Omega}\partial^{\gamma'}_{s'}u^{m+1}=0$.
We use the trace inequality and then the Poincar\'e
inequality on $\Omega$ to deduce
\[
||\partial^{\gamma'}_{s'}u^{m+1}||_{L^2(\T^{n-1})}\leq
C||\partial^{\gamma'}_{s'}u^{m+1}||_{H^1(\Omega)}\leq
C||\partial^{\gamma'}_{s'}\nabla u^{m+1}||_{L^2(\Omega)}\leq C\sqrt{\D^{m+1}}.
\]
If $|\gamma'|=0$ by the Poincar\'e inequality and the
trace inequality:
\beas
||\partial_{s'}u^{m+1}||_{L^2(\T^{n-1})}&\leq&
C||\partial_{s'}\nabla_{x'}u^{m+1}||_{L^2(\T^{n-1})}\\
&\leq&C||\partial_{s'}\nabla_{x'}u^{m+1}||_{H^1(\Omega)}
\leq C\sqrt{\D^{m+1}}.
\eeas
By part (a),
$
||\partial^{\gamma-\gamma'}_{s-s'}(\sm)||_2\leq C+C\sqrt{\E^m}.
$
Furthermore, if for some $1\leq l\leq4$ we have $|\gamma_l|+2s_l\geq k$, we estimate the
term containing $\gamma_l$ in superscript in $L^2$-norm and the remaining terms in
$L^{\infty}$-norm. Using part (a) we deduce
\[
||\partial^{\gamma_1}_{s_1}\rho^m_i\partial^{\gamma_2}_{s_2}\rho^m_j\partial^{\gamma_3}_{s_3}\rho^m_{ij}
\partial^{\gamma_4}_{s_4}(\sm^{-2})||_2\leq C(\E^m)^{3/2}\sum_{|\mu|+2s\leq 2k}||\partial^{\mu}_s\nabla\rho^m||_{H^2}.
\]
Thus, summing over all pairs $(\mu,s)$ yields
\be\label{eq:previous1}
\sum_{|\mu|+2s\leq 2k\atop\gamma=\mu+\tau, |\tau|=1}||\partial^{\gamma}_s\Delta\rho^m||_2\leq (C+C\sqrt{\E^m})\sqrt{\D^{m+1}}+
C(\E^m)^{3/2}\sum_{|\mu|+2s\leq 2k}||\partial^{\mu}_s\nabla\rho^m||_{H^2}.
\ee
Since $\int_{\T^{n-1}}\nabla\rho^m=0$ and $|\gamma|\geq1$ elliptic regularity
implies
\be\label{eq:previous2}
\sum_{|\mu|+2s\leq 2k}||\partial^{\mu}_s\nabla\rho^m||_{H^2}\leq C\sum_{|\mu|+2s\leq 2k\atop\gamma=\mu+\tau, |\tau|=1}
||\partial^{\gamma}_s\Delta\rho^m||_2.
\ee
Combining~(\ref{eq:previous1}), ~(\ref{eq:previous2}) and choosing $\E^m$ sufficiently small we deduce the claim.\\
{\em Part (c):}
We apply the differential operator $\d$ to the  'jump relation'~(\ref{eq:jumprm})
and take squares on both sides to obtain
\[
||\d\rho^{m+1}_t||_2^2+
2\epsilon||\d\Delta\rho^{m+1}_t||_2^2+
\epsilon^2||\d\Delta^2\rho^{m+1}_t||_2^2=
||\d\big(\sm^2[u^{m+1}_n]^-_+\big)||_2^2.
\]
Next
\beas
&&||\d\Big(\sm^2[u^{m+1}_n]^-_+\Big)||_2=||\sum_{\mu',s'}C^{\mu'}_{s'}
\partial^{\mu'}_{s'}(\sm^2)\partial^{\mu-\mu'}_{s-s'}[u^{m+1}_n]^-_+||_2\\
&&\leq C\sum_{|\mu'|+2s'\atop\leq k}||\partial^{\mu'}_{s'}(\sm^2)||_{\infty}
||\partial^{\mu-\mu'}_{s-s'}[u^{m+1}_n]^-_+||_2+\\
&&\quad+ C\sum_{|\mu'|+2s'\atop>k}||\partial^{\mu'}_{s'}(\sm^2)||_2
||\partial^{\mu-\mu'}_{s-s'}[u^{m+1}_n]^-_+||_{L^{\infty}(\T^{n-1})}\\
&&\leq(C\sqrt{\E^m}+1)\big(\frac{C}{\lambda}\E^{m+1}+C\lambda\D^{m+1}\big)
\leq\frac{C}{\lambda}\E^{m+1}+C\lambda\D^{m+1}.
\eeas
Here we assume $\E^m\leq1$. Terms involving $L^{\infty}$-norm are first estimated by the Sobolev inequality.
Then we use the standard trace inequality
$
||v||_{L^2(\partial\Omega)}\leq \frac{C}{\lambda}||v||_{L^2(\Omega)}+\lambda
||\nabla v||_{L^2(\Omega)}
$
to bound the terms involving $[u^{m+1}_n]^-_+$.
Observe that the same proof is easily adapted to yield the bound
\be
||\partial_s\rho^{m+1}_t||_2\leq C\sqrt{\D^{m+1}},\qquad\textrm{for}\quad s\leq k.
\ee
\prfe
{\bf Remark.}
The estimate~(\ref{eq:higher}) will play a crucial role in the energy estimates for
the problem of local existence.
For any pair  $(\mu,s)$, $|\mu|+2s\leq2k$ and $|\mu|\geq1$, the elliptic regularity and
the estimate~(\ref{eq:higher}) imply
\be\label{eq:elliptic}
||\nabla^4\partial^{\mu-1}_s\nabla\rho^{m+1}_t||_2\leq
\frac{C}{\lambda\epsilon^2}\sqrt{\E^{m+1}}+\lambda\sqrt{\D^{m+1}}.
\ee
\begin{lemma}\label{lm:aux2}
There exists a positive constant $C$ such that
\begin{enumerate}
\item[\rm (a)]
For $r\geq0$, $|\mu|+2s\leq2k$
and $(|\mu|,s)\neq(0,0)$
\be
||\d\partial_{n^r}(\am)||_{L^2(\Omega)}\leq C\sqrt{\E^m}
\ee
\item[\rm (b)]
For $r\geq0$ and $|\mu|+2s\leq2k$
\be
||\d\partial_{n^r}(\Bm)||_{L^2(\Omega)}\leq C\sqrt{\E^m}
\ee
\item[\rm (c)]
For $r\geq0$ and $|\mu|+2s\leq2k$
\be\label{eq:aux6}
||\d\partial_{n^r}(\cm)||_{L^2(\Omega)}\leq C\sqrt{\E^m}+C\sqrt{\D^m}.
\ee
Furthermore, for $r\geq0$,
$|\mu|+2s\leq 2k-1$
\be\label{eq:aux7}
||\d\partial_{n^r}(\cm)||_{L^2(\Omega)}\leq C\sqrt{\E^m}.
\ee
\end{enumerate}
\end{lemma}
\prf
{\em Part (a):}
We note:
\beas
\d\partial_{n^r}(\am)&=&\d\partial_{n^r}\big((1+|\phi\nabla\rho^m|^2)(1+\phi'\rho^m)^{-2}\big)\\
&=&\d\Big(\sum_{p=0}^r\partial_{n^p}F_1(|\phi\nabla\rho^m|)\partial_{n^{r-p}}F_2(\phi'\rho^m)\Big),
\eeas
where $F_1(x):=1+x^2$ and $F_2(x):=(1+x)^{-2}$.
Observe that
\[
\partial_{n^{r-p}}F_2(\phi'\rho^m)=
\sum_{d=1}^{r-p}\sum_{s_1+\dots s_d\atop =r-p}C_dF_2^{(d)}(\phi'\rho^m)(\rho^m)^d
\phi^{(s_1)}\dots\phi^{(s_d)},
\]
and thus
\beas
||\d\partial_{n^{r-p}}F_2(\phi'\rho^m)||_{L^2(\Omega)}&\leq&
C\sum_{d=1}^{r-p}||\d\Big[F_2^{(d)}(\phi'\rho^m)(\rho^m)^d\Big]||_{L^2(\Omega)}\\
&\leq&C\sum_{d=1}^{r-p}\sum_{\mu',s'}C^{\mu'}_{s'}
||\partial^{\mu'}_{s'}\big[[F_2^{(d)}(\phi'\rho^m)\big]
\partial^{\mu-\mu'}_{s-s'}\big[G_d(\rho^m)\big]||_{L^2(\Omega)},
\eeas
where $G_d(x):=x^d$. By Moser's inequality,
for any pair of indices $(\nu,q)$ such that $|\nu|+2q\leq 2k$, $(|\nu|,q)\neq(0,0)$ we have
$||\partial^{\nu}_q\big[[F_2^{(d)}(\phi'\rho^m)\big]||_{L^2(\Omega)}\leq
C\sqrt{\E^m}$ if $|\nu|>0$, and $||\partial^{\nu}_q\big[[F_2^{(d)}(\phi'\rho^m)\big]||_{L^2(\Omega)}\leq
C\sqrt{\E^m}+\theta$ if $|\nu|=0$. We have used Lemma~\ref{lm:tau}.
Similarly, if $|\nu|+2q\leq 2k$, we have
$||\partial^{\nu}_q\big[G_d(\rho^m)\big]||_{L^2(\Omega)}\leq C\sqrt{\E^m}$
when $|\nu|>0$, and
$||\partial^{\nu}_q\big[G_d(\rho^m)\big]||_{L^2(\Omega)}\leq C\sqrt{\E^m}+\theta$ when $|\nu|=0$.
Thus
\be\label{eq:funny}
||\partial^{\mu'}_{s'}\big[[F_2^{(d)}(\phi'\rho^m)\big]
\partial^{\mu-\mu'}_{s-s'}\big[G_d(\rho^m)\big]||_{L^2(\Omega)}\leq C\sqrt{\E^m},
\ee
where we hit the terms with lower order derivatives with
$L^{\infty}$-norms, depending on whether $|\mu'|+2s'\leq k$ or
$|\mu'|+2s'\geq k$. Additionally,
we use the assumption that $\theta<1$ and $\E^m<1$.
This implies
$
||\d\partial_{n^{r-p}}F_2(\phi'\rho^m)||_{L^2(\Omega)}\leq C\sqrt{\E^m}
$
for $0\leq p\leq r$ and $|\mu|+2s\leq2k$.
In the same way we prove
$
||\d\partial_{n^p}F_1(|\phi\nabla\rho^m|)||_{L^2(\Omega)}\leq C\sqrt{\E^m}
$
for $0\leq p\leq r$, $|\mu|+2s\leq2k$
and $(|\mu|,s)\neq(0,0)$.
Using the same idea of estimating lower order terms with $L^{\infty}$-norms
as in the proof of (\ref{eq:funny}), we conclude the proof of part (a). The proof
of part (b) follows in a completely analogous way.\\
{\em Part (c):}
To prove part (c),
recall that $\cm=\dm+\Em$, where $\dm$ and $\Em$ are
given by
\[
\dm:=\frac{\phi\Delta\rho^m}{1+\phi'\rho^m}-
\frac{(\phi^2)'|\nabla\rho^m|^2}{(1+\phi'\rho^m)^2}+
\frac{\phi''\rho^m(1+|\phi\nabla\rho^m|^2)}{(1+\phi'\rho^m)^3},\quad
\Em:=\frac{\phi\rho^m_t}{1+\phi'\rho^m}.
\]
The analogous proof as in the part (a) implies that
$||\d\partial_{n^r}(\dm)||_{L^2(\Omega)}\leq C\sqrt{\E^m}$. Furthermore,
since $||\d\rho^m_t||\leq C\sqrt{\D^m}$ if $|\mu|+2s\leq2k$,
we use the same method as in part (a) to prove
$||\d\partial_{n^r}(\Em)||_{L^2(\Omega)}\leq C\sqrt{\D^m}$. We thus conclude (\ref{eq:aux6}).
From the definition of $\E^m$ (cf.~(\ref{eq:energies})) and the assumption of Lemma~{\ref{lm:tau}}, we have
$||\d\rho^m_t||\leq C\sqrt{\E^m}$ for $|\mu|+2s\leq 2k-1$, $|\mu|>0$; also
$||\partial_s\rho^m_t||\leq C\sqrt{\E^m}+\theta$ for all $s\leq k-1$.
Now we use the same method as in the proof of part (a) to deduce (\ref{eq:aux7}).
\prfe
\begin{lemma}\label{lm:auximp}
For all pairs of indices $(\mu,s)$ such that $|\mu|+2s+r\leq k+\frac{n}{2}+3$,
the  inequality
$
||\d\partial_{n^r}u^{m+1}||_{L^2(\Omega)}\leq C\sqrt{\E^{m+1}}
$
holds.
\end{lemma}
\prf
We prove the claim by induction in $r$. In case $r=1$ the claim is obvious
from the definition of $\E^m$.
Let the claim be true for all $r\leq\vartheta$ for some $1<\vartheta<k+n/2+3$.
We have to prove the claim for $r=\vartheta+1$, i.e.
\be
||\d\partial_{n^{\vartheta+1}}u^{m+1}||_{L^2(\Omega)}\leq C\sqrt{\E^{m+1}},
\ee
if $|\mu|+2s\leq k+\frac{n}{2}-\vartheta+2$.
Let $(\mu,s)$ be such a pair of indices. Then
\be
\begin{array}{l}
\displaystyle
\d\partial_{n^{\vartheta+1}}u^{m+1}=\d\partial_{n^{\vartheta-1}}u^{m+1}_{nn}\\
\displaystyle
=\d\partial_{n^{\vartheta-1}}\Big\{\big(u^{m+1}_t-\Delta_{x'}u^{m+1}+\Bm\cdot\nabla_{x'}u^{m+1}_n+
\cm u^{m+1}_n\big)\am^{-1}\Big\}\\
\displaystyle
=\sum_{\mu',s'}\sum_{w=0}^{\vartheta-1}C(\mu',s',w)\Big\{\partial^{\mu'}_{s'}\partial_{n^w}u^{m+1}_t
-\partial^{\mu'}_{s'}\partial_{n^w}\Delta_{x'}u^{m+1}+\\
\displaystyle
+\sum_{\mu'',s''}\sum_{p=0}^wC(\mu'',s'',p)\partial^{\mu''}_{s''}\partial_{n^p}\big(
\Bm\big)\partial^{\mu'-\mu''}_{s'-s''}\partial_{n^{w-p+1}}\nabla_{x'}u^{m+1}+\\
\displaystyle
+\sum_{\mu'',s''}\sum_{p=0}^{\vartheta-1}C(\mu'',s'',p)\partial^{\mu''}_{s''}\partial_{n^p}\big(
\cm\big)\partial^{\mu'-\mu''}_{s'-s''}\partial_{n^{w-p+1}}u^{m+1}\Big\}
\partial^{\mu-\mu'}_{s-s'}\partial_{n^{\vartheta-1-w}}\big(\am^{-1}\big).
\end{array}
\label{eq:uen1}
\ee
Observe that
$
||\partial^{\mu'}_{s'}\partial_{n^w}u^{m+1}_t||_{L^2(\Omega)}\leq C\sqrt{\E^{m+1}}
$
and
$
||\partial^{\mu'}_{s'}\partial_{n^w}\Delta_{x'}u^{m+1}||_{L^2(\Omega)}\leq C\sqrt{\E^{m+1}}
$
for any triple of indices $(\mu',s',w)\leq(\mu,s,\vartheta-1)$ by inductive hypothesis.
Further, by Lemma~\ref{lm:aux2} and the Sobolev inequality
$
||\partial^{\mu''}_{s''}\partial_{n^p}\Bm||_{L^{\infty}(\Omega)}\leq C\sqrt{\E^m},
$
$
||\partial^{\mu''}_{s''}\partial_{n^p}\cm||_{L^{\infty}(\Omega)}\leq C\sqrt{\E^m}
$
and
$
||\partial^{\mu-\mu'}_{s-s'}\partial_{n^{\vartheta-1-w}}\big(\am^{-1}\big)||_{L^{\infty}(\Omega)}\leq C\sqrt{\E^m}+C,
$
for $(\mu'',s'',p)\leq(\mu',s',w)\leq(\mu,s,\vartheta-1)$.
Applying these estimates to the above identity and using the inductive assumption,
we obtain
\[
||\d\partial_{n^{\vartheta+1}}u^{m+1}||_{L^2(\Omega)}\leq C\sqrt{\E^m}\sqrt{\E^{m+1}}+C\E^m\sqrt{\E^{m+1}}
\leq C\sqrt{\E^{m+1}},
\]
where we recall the smallness assumption on $\E^m$, specially $\E^m\leq1$. This finishes the proof of the lemma.
\prfe
\begin{lemma}\label{lm:auximp1}
If $|\mu|+r+2s\leq 2k+1$,
\be
||\d\partial_{n^r}u^{m+1}||_{L^2(\Omega)}\leq C\sqrt{\E^{m+1}}.
\ee
\end{lemma}
\prf
We prove the claim by induction in $r$. In case $r=1$ the claim is obvious
from the definition of $\E^m$.
Let the claim be true for all $r\leq\vartheta$ for some $1<\vartheta<2k+1$.
We have to prove the claim for $r=\vartheta+1$, i.e.
\[
||\d\partial_{n^{\vartheta+1}}u^{m+1}||_{L^2(\Omega)}\leq C\sqrt{\E^{m+1}},
\]
if $|\mu|+2s+\vartheta\leq 2k$. Let $(\mu,s)$ be such a pair of indices.
Then
\be
\begin{array}{l}
\displaystyle
\d\partial_{n^{\vartheta+1}}u^{m+1}=\d\partial_{n^{\vartheta-1}}u^{m+1}_{nn}\\
\displaystyle
=\d\partial_{n^{\vartheta-1}}\Big\{\big(u^{m+1}_t-\Delta_{x'}u^{m+1}+\Bm\cdot\nabla_{x'}u^{m+1}_n+
\cm u^{m+1}_n\big)\am^{-1}\Big\}\\
\displaystyle
=\sum_{\mu',s'}\sum_{w=0}^{\vartheta-1}C(\mu',s',w)\Big\{\partial^{\mu'}_{s'}\partial_{n^w}u^{m+1}_t
-\partial^{\mu'}_{s'}\partial_{n^w}\Delta_{x'}u^{m+1}+\\
\displaystyle
+\sum_{\mu'',s''}\sum_{p=0}^wC(\mu'',s'',p)\partial^{\mu''}_{s''}\partial_{n^p}\big(
\Bm\big)\cdot\partial^{\mu'-\mu''}_{s'-s''}\partial_{n^{w-p+1}}\nabla_{x'}u^{m+1}+\\
\displaystyle
+\sum_{\mu'',s''}\sum_{p=0}^wC(\mu'',s'',p)\partial^{\mu''}_{s''}\partial_{n^p}\big(
\cm\big)\partial^{\mu'-\mu''}_{s'-s''}\partial_{n^{w-p+1}}u^{m+1}\Big\}
\partial^{\mu-\mu'}_{s-s'}\partial_{n^{\vartheta-1-w}}\big(\am^{-1}\big).
\end{array}
\label{eq:uen2}
\ee
We analyze
separately the case when $|\mu'|+w+2s'\leq k$ and $|\mu'|+w+2s'\geq k$.
\textbf{Case 1.}
In the case $|\mu'|+w+2s'\geq k$ note that
\[
|\mu-\mu'|+2(s-s')+\vartheta-w\leq2k-k=k.
\]
By the Sobolev inequality and Lemma~\ref{lm:aux2}, we have
\[
||\partial^{\mu-\mu'}_{s-s'}\partial_{n^{\vartheta-1-w}}\big(\am^{-1}\big)||_{L^{\infty}(\Omega)}\leq
C\sqrt{\E^m}.
\]
If $(\mu'',s'',p)\leq(\mu',s',w)$,
$|\mu''|+p+2s''\leq k$ and
\[
||\partial^{\mu''}_{s''}\partial_{n^p}\big(\cm\big)||_{L^{\infty}(\Omega)}
\leq||\partial^{\mu''}_{s''}\partial_{n^p}\big(\cm\big)||_{H^{n/2+1}(\Omega)}
\leq C\sqrt{\E^m},
\]
where we have used the Sobolev inequality and Lemma~\ref{lm:aux2} respectively.
Similarly
$
||\partial^{\mu''}_{s''}\partial_{n^p}\big(
\Bm\big)||_{L^{\infty}(\Omega)}\leq
C\sqrt{\E^m}.
$
This implies
\beas
&&\int_{\Omega}|\partial^{\mu''}_{s''}\partial_{n^p}\big(
\Bm\big)|^2|(\partial^{\mu'-\mu''}_{s'-s''}\partial_{n^{w-p+1}}\nabla_{x'}u^{m+1})|^2
(\partial^{\mu-\mu'}_{s-s'}\partial_{n^{\vartheta-1-w}}\big(\am^{-1}\big))^2\\
&&\leq C||\partial^{\mu''}_{s''}\partial_{n^p}\big(
\Bm\big)||_{L^{\infty}(\Omega)}^2
||\partial^{\mu'-\mu''}_{s'-s''}\partial_{n^{w-p+1}}\nabla_{x'}u^{m+1}||_{L^2(\Omega)}^2\times\\
&&\quad\times||\partial^{\mu-\mu'}_{s-s'}\partial_{n^{\vartheta-1-w}}\big(\am^{-1}\big)||_{L^{\infty}(\Omega)}^2\leq
C(\E^m)^2\E^{m+1},
\eeas
where we have used the inductive hypothesis to deduce
\[
||\partial^{\mu'-\mu''}_{s'-s''}\partial_{n^{w-p+1}}\nabla_{x'}u^{m+1}||_{L^2(\Omega)}^2\leq C\E^{m+1}.
\]
Analogously
\beas
&&\int_{\Omega}(\partial^{\mu''}_{s''}\partial_{n^p}\big(
\cm\big))^2(\partial^{\mu'-\mu''}_{s'-s''}\partial_{n^{w-p+1}}u^{m+1})^2
(\partial^{\mu-\mu'}_{s-s'}\partial_{n^{\vartheta-1-w}}\big(\am^{-1}\big))^2\\
&&\leq C(\E^m)^2\E^{m+1}
\eeas
If on the other hand $|\mu''|+p+2s''> k$, we use the Sobolev inequality and Lemma~\ref{lm:auximp}
to get
\[
||\partial^{\mu'-\mu''}_{s'-s''}\partial_{n^{w-p+1}}u^{m+1}||_{L^{\infty}(\Omega)}\leq
||\partial^{\mu'-\mu''}_{s'-s''}\partial_{n^{w-p+1}}u^{m+1}||_{H^{n/2+1}(\Omega)}\leq
C\sqrt{\E^{m+1}},
\]
where we note that
\[
|\mu'-\mu''|+(w-p+1)+n/2+1+2(s'-s'')\leq k+n/2+1,
\]
so that Lemma~\ref{lm:auximp} is
applicable. In analogous fashion it follows
\[
||\partial^{\mu'-\mu''}_{s'-s''}\partial_{n^{w-p+1}}\nabla_{x'}u^{m+1}||_{L^{\infty}(\Omega)}\leq C\sqrt{\E^{m+1}}.
\]
We also note that
\beas
&&\int_{\Omega}\Big(\partial^{\mu'}_{s'}\partial_{n^w}\partial u^{m+1}_t
-\partial^{\mu'}_{s'}\partial_{n^w}\Delta_{x'}u^{m+1}\Big)^2(\partial^{\mu-\mu'}_{s-s'}\partial_{n^{\vartheta-1-w}}\big(\am^{-1}\big))^2\\
&&\leq||\partial^{\mu'}_{s'}\partial_{n^w}u^{m+1}_t
-\partial^{\mu'}_{s'}\partial_{n^w}\Delta_{x'}u^{m+1}||_{L^2(\Omega)}^2
||\partial^{\mu-\mu'}_{s-s'}\partial_{n^{\vartheta-1-w}}\big(\am^{-1}\big)||_{L^{\infty}(\Omega)}^2\\
&&\leq C\E^{m+1}\E^m\leq C\E^{m+1}.
\eeas
Observe that we have used the inductive hypothesis in the last inequality above. This completes the
first case.

\textbf{Case 2.}
In the case $|\mu'|+w+2s'\leq k$,
by the Sobolev inequality and Lemmas~\ref{lm:aux2} and~\ref{lm:auximp},
\[
||\partial^{\mu''}_{s''}\partial_{n^p}\big(
\Bm\big)\partial^{\mu'-\mu''}_{s'-s''}\partial_{n^{w-p+1}}\nabla_{x'}u^{m+1}||_{L^{\infty}(\Omega)}\leq C\sqrt{\E^m}\sqrt{\E^{m+1}},
\]
and
$
||\partial^{\mu''}_{s''}\partial_{n^p}\big(
\cm\big)\partial^{\mu'-\mu''}_{s'-s''}\partial_{n^{w-p+1}}u^{m+1}||_{L^{\infty}(\Omega)}\leq C\sqrt{\E^{m+1}}
$
for $(\mu'',s'',p)\leq(\mu',s',w)$. By Lemma~\ref{lm:aux2},
$
||\partial^{\mu-\mu'}_{s-s'}\partial_{n^{\vartheta-1-w}}\big(\am^{-1}\big)||_{L^2(\Omega)}\leq C\sqrt{\E^m}+C.
$
We also note
\beas
&&\int_{\Omega}\Big(\partial^{\mu'}_{s'}\partial_{n^w}u^{m+1}_t
-\partial^{\mu'}_{s'}\partial_{n^w}\Delta_{x'}u^{m+1}\Big)^2(\partial^{\mu-\mu'}_{s-s'}\partial_{n^{\vartheta-1-w}}\big(\am^{-1}\big))^2\\
&&\leq||\partial^{\mu'}_{s'}\partial_{n^w}u^{m+1}_t
-\partial^{\mu'}_{s'}\partial_{n^w}\Delta_{x'}u^{m+1}||_{L^{\infty}(\Omega)}^2
||\partial^{\mu-\mu'}_{s-s'}\partial_{n^{\vartheta-1-w}}\big(\am^{-1}\big)||_{L^2(\Omega)}^2\\
&&\leq C\E^{m+1}\E^m+C\E^{m+1}\leq C\E^{m+1}.
\eeas
By the Sobolev inequality and Lemma~\ref{lm:auximp},
\[
||\partial^{\mu'}_{s'}\partial_{n^w}u^{m+1}_t
-\partial^{\mu'}_{s'}\partial_{n^w}\Delta_{x'}u^{m+1}||_{L^{\infty}(\Omega)}^2\leq C\E^{m+1}.
\]
We combine the above estimates to conclude
$
||\d\partial_{n^{\vartheta+1}}u^{m+1}||_{L^2(\Omega)}\leq C\sqrt{\E^{m+1}}
$
and this completes the second case and finishes the proof of the lemma.
\prfe
\section{Energy estimates}\label{ch:energye}
\begin{lemma}\label{lm:local3}
Let $K$ and $\theta\leq\theta_0$ be given as in Lemma~\ref{lm:tau}.
There exists $0<L\leq\theta$ and ${\cal T}^{\epsilon}$ such that
if
\[
\E_{\epsilon}(u^{\epsilon}_0,\rho^{\epsilon}_0)+
\Big|\int_{\T^{n-1}}\rho^{\epsilon}_0-\int_{\Omega}u^{\epsilon}_0(1+\phi'\rho^{\epsilon}_0)\Big|
\leq\frac{L}{2}
\]
and for some $m\in\N$
\[
\sup_{0\leq t\leq {\cal T}^{\epsilon}}\E^m(t)+
\int_0^{{\cal T}^{\epsilon}}\D^m(\tau)\,d\tau\leq L,
\]
\[
||\nabla\rho^{m-1}||_{\infty}\leq1,
\quad\sum_{p=0}^{k}||\partial_p\rho^m||_2\leq K\sqrt{\E^m}+\theta,
\]
then\quad
\[
\sup_{0\leq t\leq {\cal T}^{\epsilon}}\E^{m+1}(t)+
\int_0^{{\cal T}^{\epsilon}}\D^{m+1}(\tau)\,d\tau\leq L.
\]
\end{lemma}
\prf
With the preparation from Chapter~\ref{ch:basic}, we are ready to
estimate RHS of~(\ref{eq:energymm}) term by term.
Note that the assumptions of Lemma~\ref{lm:tau} are fulfilled and we
are thus able to use Lemmas~{\ref{lm:aux}} - \ref{lm:auximp1} in the forthcoming estimates.
Let $(\mu,s)$ be an arbitrary, but fixed pair of indices satisfying $|\mu|+2s\leq2k$.
\par
{\bf Term $\int_{\Omega}P^m_{\mu,s}$:}
Recall that $P^m_{\mu,s}=P(\d u^{m+1},\rho^m,f^m)$ is given by~(\ref{eq:PQ}), where $P$ is given by~(\ref{eq:P}) and
$f^m$ by~(\ref{eq:ef}).
Thus, combining~(\ref{eq:P}) and~(\ref{eq:ef}) we can estimate the first
term on RHS of $\int_{\Omega}P^m_{\mu,s}$:
\be
\begin{array}{l}
\displaystyle
\Big|\int_{\Omega}\Big\{\d(\am u^{m+1}_{nn})-\am\d u^{m+1}_{nn}\Big\}\d u^{m+1}\Big|\\
\displaystyle
\leq C\sum_{(|\mu'|,s')\neq(0,0)}\Big|\int_{\Omega}\partial^{\mu'}_{s'}\am\partial^{\mu-\mu'}_{s-s'}
u^{m+1}_{nn}\d u^{m+1}\Big|=C\sum_{|\mu'|+2s'\leq k\atop(|\mu'|,s')\neq(0,0)}+C\sum_{|\mu'|+2s'> k}\\
\displaystyle
\leq C\sqrt{\E^m}\sqrt{\D^{m+1}}\sqrt{\D^{m+1}}\leq C\sqrt{\E^m}\D^{m+1}.
\end{array}
\label{eq:z1}
\ee
In the first sum, observe that $||\partial^{\mu'}_{s'}\am||_{L^{\infty}(\Omega)}\leq C\sqrt{\E^m}$ for
$|\mu'|+2s'\leq k$, by
the Sobolev inequality and Lemma~\ref{lm:aux}. In the second sum observe that
$||\partial^{\mu-\mu'}_{s-s'}u^{m+1}_{nn}||_{L^{\infty}(\Omega)}\leq C\sqrt{\D^{m+1}}$,
by the Sobolev inequality and Lemma~\ref{lm:auximp1}.
By glancing at~(\ref{eq:P}) and~(\ref{eq:ef})
the second term in the expression $\int_{\Omega}P^m_{\mu,s}$ is given by
$
\int_{\Omega}\d(\Bm\cdot\nabla_{x'}u^{m+1}_n)\d u^{m+1}.
$
It is
estimated in a completely analogous way and is
bounded by $ C\sqrt{\E^m}\D^{m+1}$ again.
By~(\ref{eq:PQ}), the third term on RHS~(\ref{eq:P})
renders the third term in $\int_{\Omega}P^m_{\mu,s}$. We have
\be
\begin{array}{l}
\displaystyle
\Big|\int_{\Omega}\d\big(\cm u^{m+1}_n\big)\d u^{m+1}\Big|\leq
C\sum_{\mu',s'}\Big|\int_{\Omega}\partial^{\mu'}_{s'}\cm\partial^{\mu-\mu'}_{s-s'}
u^{m+1}_n\d u^{m+1}\Big|\\
\displaystyle
=C\sum_{|\mu'|+2s'\leq k}+C\sum_{|\mu'|+2s'> k}\\
\displaystyle
\leq C\sqrt{\E^m}\D^{m+1}+C\sqrt{\D^m}\sqrt{\D^{m+1}}\sqrt{\E^{m+1}}.
\end{array}
\label{eq:z2}
\ee
In the first sum we estimate $||\d\cm||_{\infty}$ like above ( since $|\mu'|+2s'\leq k$).
In the second sum, for $|\mu'|+2s'\geq k$, $||\partial^{\mu-\mu'}_{s-s'}u^{m+1}_n||_{L^{\infty}(\Omega)}\leq C\sqrt{\E^{m+1}}$
by the Sobolev inequality and Lemma~\ref{lm:auximp1}. By Lemma~\ref{lm:aux},
$||\partial^{\mu'}_{s'}\cm||_{L^2}\leq C\sqrt{\D^m}$.
Finally, the fourth term of $\int_{\Omega}P^m_{\mu,s}$ (again use~(\ref{eq:P}) to identify the fourth term
and the equations~(\ref{eq:PQ}) and~(\ref{eq:ef}) to plug in the appropriate values),
is estimated by using the Cauchy-Schwarz
inequality
\be
\begin{array}{l}
\displaystyle
\Big|\int_{\Omega}(\am)_n\d u^{m+1}_n\d u^{m+1}\Big|\leq
||(\am)_n||_{L^{\infty}(\Omega)}||\d u^{m+1}_n||_{L^2(\Omega)}||\d u^{m+1}||_{L^2(\Omega)}\\
\displaystyle
\quad\quad\leq C\sqrt{D^m}\sqrt{\D^{m+1}}\sqrt{\E^{m+1}}.
\end{array}
\label{eq:z3}
\ee
\par
{\bf Term $\int_{\T^{n-1}}Q^m_{\mu,s}$:}
We now proceed with the estimates for the expression $\int_{\T^{n-1}}Q^m_{\mu,s}$,
where $Q^m_{\mu,s}$ is given by~(\ref{eq:PQ}),
where $Q$ is defined by~(\ref{eq:Q}), $g^m_{\mu,s}$ is given
by~(\ref{eq:gsmall}) and $h^m_{\mu,s}$ is given by~(\ref{eq:ha}).
The first two terms of $\int_{\T^{n-1}}Q^m_{\mu,s}$ are the cross-terms
and they deserve special attention.
For any $\eta>0$,
\be
\begin{array}{l}
\displaystyle
\Big|\int_{\T^{n-1}}\d\nabla\rho^{m+1}_t\big\{\d\nabla\rho^m\sm^{-1}-\d\nabla\rho^{m+1}\sm^{-1}\big\}\Big|\\
\displaystyle
\leq\eta||\d\nabla\rho^{m+1}||_2^2+\frac{C}{\eta}\Big(||\d\nabla\rho^m||_2^2+||\d\nabla\rho^{m+1}||_2^2\Big)
\leq\eta\D^{m+1}+\frac{C}{\eta}(\E^m+\E^{m+1}).
\end{array}
\label{eq:z4}
\ee
\be
\begin{array}{l}
\displaystyle
\Big|\epsilon\int_{\T^{n-1}}\d\nabla\Delta\rho^{m+1}_t\cdot
\Big\{\Delta\d\nabla\rho^m\sm^{-1}-\Delta\d\nabla\rho^{m+1}\sm^{-1}\Big\}\Big|\\
\displaystyle
\leq\eta\epsilon||\d\Delta\nabla\rho^{m+1}_t||_2^2+\frac{\epsilon C}{\eta}
\Big(||\d\Delta\nabla\rho^m||_2^2
+||\d\Delta\nabla\rho^{m+1}||_2^2\Big)\\
\displaystyle
\leq\eta\D^{m+1}+\frac{C}{\eta}(\E^m+\E^{m+1}).
\end{array}
\label{eq:z5}
\ee
Observe that the constant $C$ does not depend on $\epsilon$.
The third term in $\int_{\T^{n-1}}Q^m_{\mu,s}$ is given by~(\ref{eq:Q}) and~(\ref{eq:PQ}). Note that
$
\partial_t\big(\sm^{-1}\big)=\frac{2\nabla\rho^m\nabla\rho^m_t}{\sm^3}.
$
By Lemma~\ref{lm:aux}, we
conclude $||\sm^{-1}_t||_{\infty}\leq C\D^m$. We
then obtain
\be
\begin{array}{l}
\displaystyle
\Big|\int_{\T^{n-1}}\Big\{|\d\nabla\rho^{m+1}|^2\sm^{-1}_t+\epsilon
|\d\Delta\nabla\rho^{m+1}|^2\sm^{-1}_t\Big\}\Big|\\
\displaystyle
\leq C\D^m\big(||\d\nabla\rho^{m+1}_t||_2^2+\epsilon||\d\Delta\nabla\rho^{m+1}||_2^2\big)\leq
C\E^{m+1}\D^m.
\end{array}
\label{eq:z6}
\ee
To estimate the fourth term of $\int_{\T^{n-1}}Q^m_{\mu,s}$ (which is obtained as the fourth term of~(\ref{eq:Q}) together
with the definition~(\ref{eq:PQ})), we
use the Cauchy-Schwarz inequality to get
\be
\begin{array}{l}
\displaystyle
\Big|\int_{\T^{n-1}}\d\rho^{m+1}_t\d\nabla\rho^m\cdot\nabla(\sm^{-1})\Big|\leq
||\d\rho^{m+1}_t||_2||\d\nabla\rho^m||_2||\nabla(\sm^{-1})||_{\infty}\leq\\
\displaystyle
C\sqrt{\D^{m+1}}\sqrt{\E^m}\sqrt{\D^m}
\leq C\sqrt{\E^m}(\D^m+\D^{m+1}).
\end{array}
\label{eq:z7}
\ee
Analogously, the fifth term in $\int_{\T^{n-1}}Q^m_{\mu,s}$ is estimated as follows:
\be
\begin{array}{l}
\displaystyle
\epsilon\Big|\int_{\T^{n-1}}\Delta\nabla\d\rho^{m+1}_t\cdot\nabla(\sm^{-1})\Delta\d\rho^m\Big|\leq\\
\displaystyle
||\sqrt{\epsilon}\Delta\nabla\d\rho^{m+1}_t||_2||\nabla(\sm^{-1})||_{\infty}||\sqrt{\epsilon}\Delta\d\rho^m||_2\leq
C\sqrt{\D^{m+1}}\sqrt{\D^m}\sqrt{\E^m}\\
\displaystyle
\leq C\sqrt{\E^m}(\D^m+\D^{m+1}).
\end{array}
\label{eq:z8}
\ee
We first note that the sixth term
of~(\ref{eq:Q}) contains $g$. Note that in $\int_{\T^{n-1}}Q^m_{\mu,s}$  $g=g^m_{\mu,s}$, where
$g^m_{\mu,s}$ is defined by~(\ref{eq:gsmall}). We shall first estimate $||g^m_{\mu,s}||_2$ and
$||\sqrt{\epsilon}\nabla g^m_{\mu,s}||_2$ and then
use the Cauchy-Schwarz inequality.
For any $|\mu'|+s'<|\mu|+s$, we have
\[
||\partial^{\mu'}_{s'}\Delta\rho^m\partial^{\mu-\mu'}_{s-s'}(\sm^{-1})||_2+
||\d(\nabla\rho^m\cdot\nabla(\sm^{-1}))||_2
\leq C\sqrt{\E^m}\sqrt{\D^m}.
\]
The inequality follows by estimating the term with smaller order space-derivatives
in $L^{\infty}$-norm, which can then be estimated by the Sobolev inequality and Lemma~\ref{lm:aux}.
Similarly, recalling~(\ref{eq:energies}):
\[
||\sqrt{\epsilon}\nabla\big(\partial^{\mu'}_{s'}\Delta\rho^m\partial^{\mu-\mu'}_{s-s'}(\sm^{-1})\big)||_2+
||\sqrt{\epsilon}\nabla\big(\d(\nabla\rho^m\cdot\nabla(\sm^{-1}))\big)||_2
\leq C\sqrt{\E^m}\sqrt{\D^m}.
\]
Therefore
$
||g^m_{\mu,s}||_2+||\sqrt{\epsilon}\nabla g^m_{\mu,s}||_2\leq C\sqrt{\E^m}\sqrt{\D^m}
$
and we can bound the sixth term in $\int_{\T^{n-1}}Q^m_{\mu,s}$ by
\be
\begin{array}{l}
\displaystyle
\Big|\int_{\T^{n-1}}\big(\d\rho^{m+1}_t+\epsilon\Delta^2\rho^{m+1}_t\big)g^m_{\mu,s}\Big|
\leq||\d\rho^{m+1}_t||_2||g^m_{\mu,s}||_2\\
\displaystyle
+||\sqrt{\epsilon}\Delta\nabla\rho^{m+1}_t||_2
||\sqrt{\epsilon}\nabla g^m_{\mu,s}||_2
\leq
C\sqrt{\D^{m+1}}\sqrt{\E^m}\sqrt{\D^m}
\leq C\sqrt{\E^m}(\D^m+\D^{m+1}).
\end{array}
\label{eq:z9}
\ee
The last term in $\int_{\T^{n-1}}Q^m_{\mu,s}$ is extracted from the
last term of~(\ref{eq:Q}), which contains $h$. By~(\ref{eq:PQ}),
$h=h^m_{\mu,s}$ where $h^m_{\mu,s}$ is given
by~(\ref{eq:ha}).
For the notational simplicity, we set $h^m_{\mu,s}=:h_1+\epsilon h_2$, where
\be\label{eq:ha12}
h_1:=\sum_{|\mu'|+s'\atop<|\mu|+s}C^{\mu'}_{s'}\partial^{\mu'}_{s'}\rho^{m+1}_t
\partial^{\mu-\mu'}_{s-s'}(\sm^{-2}),\,\,
h_2:=\sum_{|\mu'|+s'\atop<|\mu|+s}C^{\mu'}_{s'}\partial^{\mu'}_{s'}
\Delta^2\rho^{m+1}_t\partial^{\mu-\mu'}_{s-s'}(\sm^{-2}).
\ee
Note
that for $|\mu'|+s'<|\mu|+s$, we have
\be
\begin{array}{l}
\displaystyle
||\sm^2\partial^{\mu'}_{s'}\rho^{m+1}_t
\partial^{\mu-\mu'}_{s-s'}(\sm^{-2})||_2\\
\displaystyle
\leq||\sm^2||_{\infty}||\partial^{\mu'}_{s'}\rho^{m+1}_t
\partial^{\mu-\mu'}_{s-s'}(\sm^{-2})||_2
\leq C\sqrt{\E^m}\sqrt{\D^{m+1}}.
\end{array}
\label{eq:aid1}
\ee
Here,
if $|\mu'|+2s'\leq k$ then $||\partial^{\mu'}_{s'}\rho^{m+1}_t||_{\infty}\leq C\sqrt{\E^{m+1}}$ and
if $|\mu'|+2s'\geq k$ then $||\partial^{\mu-\mu'}_{s-s'}(\sm^{-2})||_{\infty}\leq C\sqrt{\E^m}$.
We use the Sobolev inequality and Lemma~{\ref{lm:aux}} to conclude the estimate. This implies
$||h_1||_2\leq C\sqrt{\E^m}\sqrt{\D^{m+1}}$. Note that in similar fashion,
for $|\mu'|+s'<|\mu|+s$,
\[
||\sqrt{\epsilon}\partial^{\mu'}_{s'}\Delta^2\rho^{m+1}_t\partial^{\mu-\mu'}_{s-s'}(\sm^{-2})||_2
\leq C\sqrt{\E^m}\sqrt{\D^{m+1}}.
\]
In other words $||\sqrt{\epsilon}h_2||_2\leq C\sqrt{\E^m}\sqrt{\D^{m+1}}$.
From the proof of part (b) of Lemma~\ref{lm:aux}, we deduce
$
||\d\kappa^m||_2\leq C\sqrt{\D^{m+1}}
$
and also
$
||\sqrt{\epsilon}\d\kappa^m||_2\leq C\sqrt{\D^{m+1}}
$
(recall here~(\ref{eq:kappam})).
Thus the last term of $\int_{\T^{n-1}}Q^m_{\mu,s}$ is bounded by
\be
\begin{array}{l}
\displaystyle
\Big|\int_{\T^{n-1}}\sm^2h^m_{\mu,s}\d\kappa^m\Big|\leq C||h_1||_2||\d\kappa^m||_2+C||\sqrt{\epsilon}h_2||_2||\sqrt{\epsilon}\d\kappa^m||_2\\
\displaystyle
\leq C\sqrt{\E^m}\D^{m+1}.
\end{array}
\label{eq:z10}
\ee
\par
{\bf Term $\int_{\Omega}R^m_{\mu,s}$:}
Note that $R^m_{\mu,s}$ is given by~(\ref{eq:RST}) where $f^m_{\mu,s}$, is defined by~(\ref{eq:ef}).
Our first task is to estimate the first term of $\int_{\Omega}R^m_{\mu,s}$,
namely $\int_{\Omega}(f^m_{\mu,s})^2$.
Observe that
\be
\begin{array}{l}
\displaystyle
\int_{\Omega}(f^m_{\mu,s})^2
\leq C\sum_{(\mu',s')\atop\neq(0,0)}\int_{\Omega}
(\partial^{\mu'}_{s'}\am)^2(\partial^{\mu-\mu'}_{s-s'}u^{m+1}_{nn})^2\\
\displaystyle
+C\sum_{(\mu',s')}\int_{\Omega}|\partial^{\mu'}_{s'}\Bm|^2|\partial^{\mu-\mu'}_{s-s'}\nabla_{x'}u^{m+1}_n|^2
+C\sum_{\mu',s'}\int_{\Omega}(\partial^{\mu'}_{s'}\cm)^2(\partial^{\mu-\mu'}_{s-s'}u^{m+1}_n)^2
\end{array}
\label{eq:aid2}
\ee
If $|\mu'|+2s'\leq k$, by Lemma~\ref{lm:aux2}
\[
||\partial^{\mu'}_{s'}\am||_{L^{\infty}(\Omega)}+
||\partial^{\mu'}_{s'}\Bm||_{L^{\infty}(\Omega)}+||\partial^{\mu'}_{s'}\cm||_{L^{\infty}(\Omega)}\leq C\sqrt{\E^m}.
\]
Thus for $|\mu'|+2s'\leq k$ and $(\mu',s')\neq(0,0)$, RHS of~(\ref{eq:aid2}) is bounded
by~$C\E^m\D^{m+1}$.
If $|\mu'|+2s'>k$, then $|\mu-\mu'|+2(s-s')\leq k-1$ and from Lemma~\ref{lm:auximp1}
\[
||\partial^{\mu-\mu'}_{s-s'}u^{m+1}_{nn}||_{L^{\infty}(\Omega)}+
||\partial^{\mu-\mu'}_{s-s'}\nabla_{x'}u^{m+1}_n||_{L^{\infty}(\Omega)}+
||\partial^{\mu-\mu'}_{s-s'}u^{m+1}_n||_{L^{\infty}(\Omega)}\,\leq
C\sqrt{\E^{m+1}}.
\]
In addition to this, for such $(\mu',s')$, we use Lemma~\ref{lm:aux2} to conclude
\[
||\partial^{\mu'}_{s'}\am||_{L^2(\Omega)}+||\partial^{\mu'}_{s'}\Bm||_{L^2(\Omega)}+
||\partial^{\mu'}_{s'}\cm||_{L^2(\Omega)}\,\leq C\sqrt{\D^m}.
\]
Thus, for every $\lambda>0$
\beas
&&\int_{\Omega}
(\partial^{\mu'}_{s'}\am)^2(\partial^{\mu-\mu'}_{s-s'}u^{m+1}_{nn})^2\leq
||\partial^{\mu'}_{s'}\am||_{L^2(\Omega)}^2||\partial^{\mu-\mu'}_{s-s'}u^{m+1}_{nn}||_{L^{\infty}(\Omega)}^2
\leq CD^m\E^{m+1}.
\eeas
Analogously,
\beas
&&\int_{\Omega}\Big\{|\partial^{\mu'}_{s'}\Bm|^2|\partial^{\mu-\mu'}_{s-s'}\nabla_{x'}u^{m+1}_n|^2+
(\partial^{\mu'}_{s'}\cm)^2(\partial^{\mu-\mu'}_{s-s'}u^{m+1}_n)^2\Big\}
\leq
CD^m\E^{m+1},
\eeas
for all $(\mu',s')\leq(\mu,s)$ satisfying $|\mu'|+2s'>k$.
Combining the estimates for $|\mu'|+2s'\leq k$ and $|\mu'|+2s'> k$, we obtain
\be
\int_{\Omega}(f^m_{\mu,s})^2\leq C\E^m\D^{m+1}+CD^m\E^{m+1}
\label{eq:z11}
\ee
The second term of the expression $\int_{\Omega}R^m_{\mu,s}$ (the
second term on RHS of~(\ref{eq:R}) and~(\ref{eq:RST})) is bounded by:
\be\label{eq:z13}
\Big|\int_{\Omega}(\d u^{m+1}_n)^2(\am)_t\Big|\leq
||(\am)_t||_{L^{\infty}(\Omega)}||\d u^{m+1}_n||_{L^2(\Omega)}^2\\
\leq C\sqrt{\E^m}\D^{m+1}
\ee
and similarly the third term
$
\Big|\int_{\Omega}\d u^{m+1}_t\d u^{m+1}_n(\am)_n\Big|\leq C\sqrt{\E^m}\D^{m+1}.
$
Note further that  for the fourth term in $\int_{\Omega}R^m_{\mu,s}$ (use~(\ref{eq:R}) and~(\ref{eq:RST})),
by Lemma~\ref{lm:aux},
\be
\begin{array}{l}
\displaystyle
\Big|\int_{\Omega}\nabla_{x'}\d u^{m+1}_n\cdot\nabla_{x'}\am\d u^{m+1}_n\Big|\\
\displaystyle
\leq
C||\nabla_{x'}\d u^{m+1}_n||_{L^2(\Omega)}||\nabla_{x'}\am||_{L^{\infty}(\Omega)}||\d u^{m+1}_n||_{L^2(\Omega)}
\leq C\sqrt{\D^{m+1}}\sqrt{\E^m}\sqrt{\D^{m+1}}.
\end{array}
\label{eq:z14}
\ee
By Lemma~\ref{lm:aux}, the last term in $\int_{\Omega}R^m_{\mu,s}$ (last term on RHS of~(\ref{eq:R})
and~(\ref{eq:RST})) is bounded by:
\be
\begin{array}{l}
\displaystyle
\Big|\int_{\Omega}\Delta_{x'}\d u^{m+1}(\am)_n\d u^{m+1}_n\Big|\\
\displaystyle
\leq||\Delta_{x'}\d u^{m+1}||_{L^2(\Omega)}||(\am)_n||_{L^{\infty}(\Omega)}
||\d u^{m+1}_n||_{L^2(\Omega)}\leq C\sqrt{\E^m}\D^{m+1}.
\end{array}
\label{eq:z15}
\ee
\par
{\bf Term $\int_{\T^{n-1}}S^m_{\mu,s}$:}
Note that $S^m_{\mu,s}$ is given by~(\ref{eq:RST}) where $S$ is given by~(\ref{eq:S}), $g^m_{\mu,s}$
by~(\ref{eq:gsmall}) and $h^m_{\mu,s}$ by~(\ref{eq:ha}).
The first two terms on RHS of~(\ref{eq:S}) are the cross-terms and in order to estimate them
we shall exploit the part (c) of Lemma~\ref{lm:aux}. It turns out that
the constants on the right-hand side will depend on $\epsilon$.
Note that for any $\eta,\lambda>0$
\be
\begin{array}{l}
\displaystyle
\Big|\int_{\T^{n-1}}\d\nabla\rho^{m+1}_t\Big(\d\nabla\rho^m_t\sm^{-1}-
\d\nabla\rho^{m+1}_t\sm^{-1}\Big)\Big|\\
\displaystyle
\leq
\frac{C}{\eta}||\d\nabla\rho^{m+1}_t||_2^2+\eta||\d\nabla\rho^m_t||_2^2
\leq\frac{C}{\eta\lambda\epsilon^4}\E^{m+1}+\frac{C\lambda}{\eta}\D^{m+1}+\eta\D^m.
\end{array}
\label{eq:z16}
\ee
In the last estimate we have used the estimate (\ref{eq:elliptic}).
Similarly,
\be
\begin{array}{l}
\displaystyle
\Big|\epsilon\int_{\T^{n-1}}\d\Delta\nabla\rho^{m+1}_t\Big(
\d\Delta\nabla\rho^m_t\sm^{-1}-\d\Delta\nabla\rho^{m+1}_t\sm^{-1}\Big)\Big|\\
\displaystyle
\leq\frac{\epsilon C}{\eta}||\d\Delta\nabla\rho^{m+1}_t||_2^2+\eta\epsilon
||\d\Delta\nabla\rho^m_t||_2^2
\leq
\frac{C}{\eta\lambda\epsilon^3}\E^{m+1}+\frac{C\epsilon\lambda}{\eta}\D^{m+1}+\eta\D^m.
\end{array}
\label{eq:z17}
\ee
The third term of $\int_{\T^{n-1}}S^m_{\mu,s}$  is given by the third term
on RHS of~(\ref{eq:S}) and~(\ref{eq:RST}):
\be
\begin{array}{l}
\displaystyle
\Big|\int_{\T^{n-1}}\d\rho^{m+1}_t\d\nabla\rho^m_t\cdot\nabla(\sm^{-1})\Big|\leq\\
\displaystyle
||\d\rho^{m+1}_t||_2||\d\nabla\rho^m_t||_2||\nabla(\sm^{-1})||_{\infty}\leq
C\sqrt{\D^{m+1}}\sqrt{\D^m}\sqrt{\E^m}\\
\displaystyle
C\sqrt{\E^m}(\D^m+\D^{m+1}).
\end{array}
\label{eq:z18}
\ee
Similarly, the fourth term in $\int_{\T^{n-1}}S^m_{\mu,s}$ (given by the fourth term on RHS of~(\ref{eq:S})
and~(\ref{eq:RST})) is bounded by:
\be
\begin{array}{l}
\displaystyle
\epsilon\Big|\int_{\T^{n-1}}\Delta\d\rho^m_t\Delta\nabla\d\rho^{m+1}_t\cdot\nabla(\sm^{-1})\Big|\leq\\
\displaystyle
||\sqrt{\epsilon}\Delta\d\rho^m_t||_2||\sqrt{\epsilon}\Delta\nabla\d\rho^{m+1}_t||_2||\nabla(\sm^{-1})||_{\infty}\leq
C\sqrt{\D^m}\sqrt{\D^{m+1}}\sqrt{\E^m}\\
\displaystyle
\leq C\sqrt{\E^m}(\D^m+\D^{m+1}).
\end{array}
\label{eq:z19}
\ee
Note that the fifth term in $\int_{\T^{n-1}}S^m_{\mu,s}$,
by~(\ref{eq:S}) and~(\ref{eq:RST})) involves the function $g^m_{\mu,s}$,
where $g^m_{\mu,s}$ is given by~(\ref{eq:gsmall}).
The crucial step in estimating this term in $\int_{\T^{n-1}}S^m$, is to
observe that $||(g^m_{\mu,s})_t||_2\leq C\sqrt{\E^m}\sqrt{\D^m}$. This is proved
by first differentiating $g^m_{\mu,s}$ with respect to $t$, and then in each product
estimating the
terms with lower order space derivatives in $L^{\infty}$-norm and the other one
in $L^2$-norm. The same method applies to show $||\sqrt{\epsilon}\nabla(g^m_{\mu,s})_t||_2\leq
C\sqrt{\E^m}\sqrt{\D^m}$. Also observe that
\[
||\Delta\d\rho^m\sm^{-1}_t||_2\leq||\Delta\d\rho^m||_2||\sm^{-1}_t||_{\infty}\leq
C\sqrt{\E^m}\sqrt{\D^m}.
\]
Analogous proof shows that
\[
||\sqrt{\epsilon}\nabla\big(\Delta\d\rho^m\sm^{-1}_t\big)||_2
\leq C\sqrt{\E^m}\sqrt{\D^m}.
\]
Using the above inequalities and the Cauchy-Schwarz inequality we
establish \be
\begin{array}{l}
\displaystyle
\Big|\int_{\T^{n-1}}\d\rho^{m+1}_t
\Big(\Delta\d\rho^m\sm^{-1}_t+(g^m_{\mu,s})_t\Big)\Big|
\leq||\d\rho^{m+1}_t||_2||\Delta\d\rho^m\sm^{-1}_t+(g^m_{\mu,s})_t||_2\\
\displaystyle
\quad\quad\leq C\sqrt{\D^{m+1}}\sqrt{\E^m}\sqrt{\D^m}
\leq C\sqrt{\E^m}\big(\D^m+\D^{m+1}\big).
\end{array}
\label{eq:z20}
\ee
Integrating by parts and using the analogous argument as in~(\ref{eq:z20})
we get
\be
\begin{array}{l}
\displaystyle
\epsilon\Big|\int_{\T^{n-1}}\d\Delta^2\rho^{m+1}_t
\Big(\Delta\d\rho^m\sm^{-1}_t+(g^m_{\mu,s})_t\Big)\Big|=\epsilon\Big|\int_{\T^{n-1}}\d\Delta\nabla\rho^{m+1}_t\cdot
\nabla\Big(\Delta\d\rho^m\sm^{-1}_t+(g^m_{\mu,s})_t\Big)\Big|\\
\displaystyle
\leq||\sqrt{\epsilon}\d\Delta\nabla\rho^{m+1}_t||_2
||\sqrt{\epsilon}\nabla\Big(\Delta\d\rho^m\sm^{-1}_t+(g^m_{\mu,s})_t\Big)||_2\\
\displaystyle
\leq C\sqrt{\D^{m+1}}\sqrt{\E^m}\sqrt{\D^m}
\leq C\sqrt{\E^m}\big(\D^m+\D^{m+1}\big)
\end{array}
\label{eq:z21}
\ee
The sixth and the last term in $\int_{\T^{n-1}}S^m$ ( given through the last term
on RHS of~(\ref{eq:S}) and~(\ref{eq:RST})) involves the function $h^m_{\mu,s}$,
where $h^m_{\mu,s}$ is given by~(\ref{eq:ha}). Since
$\d\kappa^m_t=\nabla\cdot\d\big(\nabla\rho^m\sm^{-1}_t\big)$, we can integrate by parts to obtain
\[
\int_{\T^{n-1}}\sm^2h^m_{\mu,s}\d\kappa^m_t=-\int_{\T^{n-1}}\nabla\big(\sm^2h^m_{\mu,s}\big)\cdot\big(\nabla\rho^m\sm^{-1}\big)_t.
\]
We split $h^m_{\mu,s}=h_1+\epsilon h_2$ as in~(\ref{eq:ha12}).
For $|\mu'|+s'<|\mu|+s$
\be
\begin{array}{l}
\displaystyle
||\nabla\big(\sm^2h_1\big)||_2
\leq||\nabla(\sm^2)||_{\infty}||\partial^{\mu'}_{s'}\rho^{m+1}_t
\partial^{\mu-\mu'}_{s-s'}(\sm^{-2})||_2+\\
\displaystyle
+||\sm^2||_{\infty}\Big(||\partial^{\mu'}_{s'}\nabla\rho^{m+1}_t
\partial^{\mu-\mu'}_{s-s'}(\sm^{-2})||_2+||\partial^{\mu'}_{s'}\rho^{m+1}_t
\partial^{\mu-\mu'}_{s-s'}\nabla(\sm^{-2})||_2\Big)\\
\displaystyle
\leq C\sqrt{\E^m}\sqrt{\E^m}\sqrt{\D^{m+1}}+ C\sqrt{\E^m}\sqrt{\D^{m+1}}
\leq C\sqrt{\E^m}\sqrt{\D^{m+1}}.
\end{array}
\label{eq:bosna}
\ee
Here we have used $L^2-L^{\infty}$ estimated by separating the cases
$|\mu'|+2s'\leq k$ and $|\mu'|+2s'> k$. Since
\[
||\sqrt{\epsilon}\nabla\big(\sm^2\partial^{\mu'}_{s'}\Delta^2\rho^{m+1}_t
\partial^{\mu-\mu'}_{s-s'}(\sm^{-2})\big)||_2\leq C\sqrt{\E^m}\sqrt{\D^{m+1}},
\]
we deduce
\be\label{eq:hercegovina}
||\sqrt{\epsilon}\nabla(\sm^2h_2)||_2\leq C\sqrt{\E^m}\sqrt{\D^{m+1}},
\ee
where $h_2$ is given by~(\ref{eq:ha12}). Similarly,
$
||\partial^{\mu}_{s+1}\big(\nabla\rho^m\sm^{-1}\big)||_2\leq C\sqrt{\E^m}\sqrt{\D^m}
$
and also
$
||\sqrt{\epsilon}\partial^{\mu}_{s+1}\big(\nabla\rho^m\sm^{-1}\big)||_2\leq C\sqrt{\E^m}\sqrt{\D^m}.
$
\par
In summary,
\be
\begin{array}{l}
\displaystyle
\Big|\int_{\T^{n-1}}\sm^2h^m_{\mu,s}\d\kappa^m_t\Big|\leq||\nabla\big(\sm^2h_1\big)||_2||\d\big(\nabla\rho^m\sm^{-1}\big)||_2+\\
\displaystyle
+||\sqrt{\epsilon}\big(\sm^2h_2\big)||_2||\sqrt{\epsilon}\d\big(\nabla\rho^m\sm^{-1}\big)||_2\\
\displaystyle
\leq C\E^m\sqrt{\D^m}\sqrt{\D^{m+1}}\leq C\E^m\D^m+C\E^m\D^{m+1}.
\end{array}
\label{eq:z22}
\ee
\par
{\bf Term $\int_{\T^{n-1}}T^m_{\mu,s}$:}
Recall that $T^m_{\mu,s}$ is defined by~(\ref{eq:RST}) where $T$ is given by~(\ref{eq:bigT}), $G^m_{\mu,s}$
by~(\ref{eq:gbig}) and $h^m_{\mu,s}$ by~(\ref{eq:ha}). In particular
the term $A$- the first term on RHS of~(\ref{eq:bigT}) is given by~(\ref{eq:A}).
The first two terms of the expression $A$ are the cross-terms.
Using integration by parts, we obtain
\be
\begin{array}{l}
\displaystyle
\Big|\int_{\T^{n-1}}\d\Delta\rho^{m+1}_t\big(\d\Delta\rho^m-\d\Delta\rho^{m+1}\big)\sm^{-1}\Big|=\\
\displaystyle
\Big|\int_{\T^{n-1}}\d\nabla\rho^{m+1}_t\cdot\nabla\Big(\big(\d\Delta\rho^m-\d\Delta\rho^{m+1}\big)\sm^{-1}\Big)\Big|\\
\displaystyle
\leq\lambda||\d\nabla\rho^{m+1}_t||_2^2+
\frac{C}{\lambda}||\nabla\big(\d\Delta\rho^m\sm^{-1}\big)||_2^2+\\
\displaystyle
\frac{C}{\lambda}||\nabla\big(\d\Delta\rho^{m+1}\sm^{-1}\big)||_2^2
\leq\lambda\D^{m+1}+\frac{C}{\epsilon\lambda}(\E^m+\E^{m+1}),
\end{array}
\label{eq:z23}
\ee
where we note that by the definition of $\E^m$,
$
||\partial^{\mu}_s(\nabla\rho^m\sm^{-1})||_{H^3}\leq\frac{C}{\sqrt{\epsilon}}\sqrt{\E^m}
$
for $|\mu|+2s\leq 2k$.
Similarly, for the second cross term:
\be
\begin{array}{l}
\displaystyle
\Big|\int_{\T^{n-1}}\d\Delta\rho^{m+1}_t\Big(\rho^m_i\rho^m_j\big(\d\rho^m_{ij}-\d\rho^{m+1}_{ij}\big)\Big)\sm^{-3}\Big|\\
\displaystyle
\quad\quad\leq\lambda\D^{m+1}+\frac{C}{\epsilon\lambda}(\E^m+\E^{m+1}).
\end{array}
\label{eq:z24}
\ee
The proof of~(\ref{eq:z24}) relies on the same idea as above; we first integrate by parts and then establish
the estimate
\[
||\nabla\Big(\rho^m_i\rho^m_j\big(\d\rho^m_{ij}-\d\rho^{m+1}_{ij}\big)\sm^{-3}\Big)||_2^2\leq
\frac{C}{\epsilon}(\E^m+\E^{m+1}).
\]
By $A-\mbox{(crossterms)}$ we denote the sum of all the remaining terms
in the expression $A$ (recall~(\ref{eq:A}) and the fact
that $\psi=\rho^m$ in our case).
Terms of the form $\sm^{-1}_i$, $\sm^{-1}_t$, $\sm^{-3}_t$, $\rho^m_i$ and $\rho^m_{ij}$ for $1\leq i,j\leq n-1$
are bounded in $L^{\infty}$-norm by $C\sqrt{\E^m}$, by
Lemma~\ref{lm:aux}. Terms of the form $\sm^{-3}$ are estimated by $1$ in $L^{\infty}$-norm. Note that
in the last two terms in the expression $A$ (\ref{eq:A}) the leading order derivatives cancel out after the
the product rule has been applied within the parentheses. Using these observations and applying the Cauchy-Schwarz
inequality, we conclude
\be\label{eq:z25}
|A-\mbox{(crossterms)}|\leq C\sqrt{\E^m}(\D^m+\D^{m+1}).
\ee
Recall now that the term $B$ (the second expression on RHS of~(\ref{eq:bigT})) is given by~(\ref{eq:B}).
The first two terms in the expression $B$ are again the cross-terms.
By part (c) of Lemma~\ref{lm:aux}, we obtain
\be
\begin{array}{l}
\displaystyle
\Big|\epsilon\int_{\T^{n-1}}\d\Delta^2\rho^{m+1}_t\big(\d\Delta^2\rho^m-\d\Delta^2\rho^{m+1}\big)\sm^{-1}\Big|\leq\\
\displaystyle
\epsilon^2||\d\Delta^2\rho^{m+1}_t||_2^2+C||\d\Delta^2\rho^m||_2^2+C||\d\Delta^2\rho^{m+1}||_2^2\leq\\
\displaystyle
\leq
\frac{C}{\lambda}\E^{m+1}+(\lambda+C\E^m)\D^{m+1}+ \frac{C}{\epsilon}\big(\E^m+\E^{m+1}\big).
\end{array}
\label{eq:z26}
\ee
Analogously, we establish
\be
\begin{array}{l}
\displaystyle
\Big|\epsilon\int_{\T^{n-1}}\d\Delta^2\rho^{m+1}_t
\Big(\rho^m_i\rho^m_j\big(\d\Delta\rho^m_{ij}-\d\Delta\rho^{m+1}_{ij}\big)\Big)\sm^{-3}\Big|\leq\\
\displaystyle
\quad\quad\frac{C}{\lambda}\E^{m+1}+(\lambda+C\E^m)\D^{m+1}+ \frac{C}{\epsilon}\big(\E^m+\E^{m+1}\big).
\end{array}
\label{eq:z27}
\ee
We denote the sum of the remaining terms in the expression $B$ by $B-\mbox{(crossterms)}$. The
same idea as in the estimates for $|A-\mbox{(crossterms)}|$ works.
 It is important to note that we have canceling of the
highest order derivatives within the parentheses in the last three expressions on RHS
of~(\ref{eq:B}).
In addition to that, we factorize $\epsilon=\sqrt{\epsilon}\times\sqrt{\epsilon}$.
By the Cauchy-Schwarz inequality,
\be\label{eq:z28}
|B-\mbox{(crossterms)}|\leq C\sqrt{\E^m}(\D^m+\D^{m+1}).
\ee
The third term of $\int_{\T^{n-1}}T^m_{\mu,s}$ is given by the third term on RHS of~(\ref{eq:bigT}) together
with~(\ref{eq:RST}). Recall that $G^m_{\mu,s}$ is given by~(\ref{eq:gbig}).
In order to estimate it, we first integrate by parts.
\beas
&&\int_{\T^{n-1}}\Delta G^m_{\mu,s}\big(\d\rho^{m+1}_t+\epsilon\d\Delta^2\rho^{m+1}_t\big)=\\
&&\quad\quad-\int_{\T^{n-1}}\nabla G^m_{\mu,s}\cdot\d\nabla\rho^{m+1}_t
-\epsilon\int_{\T^{n-1}}\Delta\nabla G^m_{\mu,s}\cdot\d\Delta\nabla\rho^{m+1}_t.
\eeas
The crucial observation is
$
||\nabla G^m_{\mu,s}||_2\leq C\sqrt{\E^m}\sqrt{\D^m},
$
$
||\sqrt{\epsilon}\Delta\nabla G^m_{\mu,s}||_2\leq C\sqrt{\E^m}\sqrt{\D^m}.
$
Both inequalities follow in the standard way, by using $L^{\infty}-L^2$ type estimates and~(\ref{eq:energies}).
The third term of $\int_{\T^{n-1}}T^m_{\mu,s}$ is then bounded by:
\be
\begin{array}{l}
\displaystyle
\Big|\int_{\T^{n-1}}\Delta G^m_{\mu,s}\big(\d\rho^{m+1}_t+\epsilon\d\Delta^2\rho^{m+1}_t\big)\Big|\leq
||\nabla G^m_{\mu,s}||_2||\d\nabla\rho^{m+1}_t||_2+\\
\displaystyle
+||\sqrt{\epsilon}\Delta\nabla G^m_{\mu,s}||_2||\sqrt{\epsilon}\d\Delta\nabla\rho^{m+1}_t||_2
\leq C\sqrt{\E^m}\sqrt{\D^m}\sqrt{\D^{m+1}}\leq C\sqrt{\E^m}(\D^m+\D^{m+1}).
\end{array}
\label{eq:z29}
\ee
We integrate by parts in the  fourth and the last term of $\int_{\T^{n-1}}T^m_{\mu,s}$
(given by the last term on RHS of~(\ref{eq:bigT}) and~(\ref{eq:RST})).
Recall
$
\Delta\U|_{\T^{n-1}}=\Delta\d u^{m+1}|_{\T^{n-1}}=\Delta\d\kappa^m.
$
\beas
&&\int_{\T^{n-1}}\sm^2h^m_{\mu,s}\Delta\d\kappa^m=\\
&&\quad\quad-\int_{\T^{n-1}}\nabla(\sm^2h_1)\cdot\d\nabla\kappa^m-
\epsilon\int_{\T^{n-1}}\nabla(\sm^2h_2)\cdot\d\nabla\kappa^m.
\eeas
By~(\ref{eq:ha12}), we set $h^m_{\mu,s}=h_1+\epsilon h_2$. By the trace inequality,
\[
||\nabla\d\kappa^m||_2=||\d\nabla_{x'}u^{m+1}||_{L^2(\T^{n-1})}\leq
C||\d\nabla u^{m+1}||_{H^1(\Omega)}\leq C\sqrt{\D^{m+1}}.
\]
By~(\ref{eq:bosna}), we obtain
\be\label{eq:hope1}
\Big|\int_{\T^{n-1}}\nabla(\sm^2h_1)\cdot\d\nabla\kappa^m\Big|\leq
||\nabla(\sm^2h_1)||_2||\d\nabla\kappa^m||_2\leq C\sqrt{\E^m}\D^{m+1}.
\ee
On the other hand, by Lemma~\ref{lm:aux} and $L^{\infty}-L^2$ type estimates,
$||\sqrt{\epsilon}\d\nabla\kappa^m||_2\leq C\sqrt{\D^m}$.
By~(\ref{eq:hercegovina}),
\be
\begin{array}{l}
\displaystyle
\Big|\epsilon\int_{\T^{n-1}}\nabla(\sm^2h_2)\cdot\d\nabla\kappa^m\Big|\leq
||\sqrt{\epsilon}\nabla(\sm^2h_2)||_2||\sqrt{\epsilon}\d\nabla\kappa^m||_2\\
\displaystyle
\quad\leq C\sqrt{\E^m}\sqrt{\D^{m+1}}\sqrt{\D^m}\leq C\sqrt{\E^m}(\D^m+\D^{m+1}),
\end{array}
\label{eq:hope2}
\ee
where we recall~(\ref{eq:ha12}) again.
From the estimates~(\ref{eq:hope1}) and~(\ref{eq:hope2}),
the last term in $\int_{\T^{n-1}}T^m_{\mu,s}$ is bounded by
\be\label{eq:z30}
\Big|\int_{\T^{n-1}}\sm^2h^m_{\mu,s}\Delta\d\kappa^m\Big|\leq C\sqrt{\E^m}(\D^m+\D^{m+1}).
\ee
Using the identity~(\ref{eq:energymm}) and
summing the estimates~(\ref{eq:z1}) -~(\ref{eq:z21}) and~(\ref{eq:z22}) -~(\ref{eq:z30}) to
get a bound on the right-hand side of the identity~(\ref{eq:energymm}), we arrive at
\be
\begin{array}{l}
\displaystyle
\frac{d}{dt}\E^{m+1}(t)+\D^{m+1}(t)\leq C\sqrt{\E^m}(\D^m+\D^{m+1})+C\sqrt{\D^m}\sqrt{\D^{m+1}}\sqrt{\E^{m+1}}\\
\displaystyle
+C(\eta+\lambda)\D^{m+1}+\frac{C}{\eta}(\E^m+\E^{m+1})+\frac{C}{\lambda}\E^{m+1}+\eta\D^m\\
\displaystyle
+\frac{C}{\eta\lambda\epsilon^4}\E^{m+1}+\frac{C\lambda\epsilon}{\eta}\D^{m+1}+
\frac{C}{\lambda\epsilon}(\E^m+\E^{m+1})+C\E^{m+1}\D^m.
\end{array}
\label{eq:energyest}
\ee
Now integrate in time over the interval $[0,t]$ to get
\be
\begin{array}{l}
\displaystyle
\E^{m+1}(t)+\int_0^t\D^{m+1}(\tau)\,d\tau\leq\E^{m+1}(0)+C\sup_{s\leq t}\sqrt{\E^m(s)}
\Big(\int_0^t\D^m(\tau)\,d\tau+\int_0^t\D^{m+1}(\tau)\,d\tau\Big)+\\
\displaystyle
+C\sup_{s\leq t}\sqrt{\E^{m+1}(s)}\int_0^t\sqrt{\D^m(s)}\sqrt{\D^{m+1}(s)}\,ds+
C(\eta+\lambda)\int_0^t\D^{m+1}(\tau)\,d\tau\\
\displaystyle
+Ct\big(\frac{1}{\eta}+\frac{1}{\lambda}\big)\sup_{s\leq t}\E^{m+1}(s)+\frac{Ct}{\eta}\sup_{s\leq t}\E^m(s)+
\frac{Ct}{\eta\lambda\epsilon^4}\sup_{s\leq t}\E^{m+1}(s)+\\
\displaystyle
+\frac{C\lambda\epsilon}{\eta}\int_0^t\D^{m+1}(\tau)\,d\tau+\eta\int_0^t\D^m(\tau)\,d\tau+
\frac{Ct}{\lambda\epsilon}\Big(\sup_{s\leq t}\E^m(s)+\sup_{s\leq t}\E^{m+1}(s)\Big)+\\
\displaystyle
+\sup_{s\leq t}\E^{m+1}(s)\int_0^t\D^m(\tau)\,d\tau.
\end{array}
\label{eq:energyest1}
\ee
Note that by Cauchy-Schwarz inequality we have
\be
\begin{array}{l}
\displaystyle
\sup_{s\leq t}\sqrt{\E^{m+1}(s)}\int_0^t\sqrt{\D^m(s)}\sqrt{\D^{m+1}(s)}\,ds\leq\\
\displaystyle
\leq\sup_{s\leq t}\sqrt{\E^{m+1}(s)}\sqrt{\int_0^t\D^m(s)\,ds}\sqrt{\int_0^t\D^{m+1}(s)\,ds}\\
\displaystyle
\leq\big(\int_0^t\D^m(\tau)\,d\tau\big)^{1/2}\big(\sup_{s\leq t}\E^{m+1}(s)+\int_0^t\D^{m+1}(\tau)\,d\tau\big).
\end{array}
\label{eq:vazno9}
\ee
By assumption $\sup_{s\leq t}\E^m(s)+\int_0^t\D^m(\tau)\,d\tau\leq L$, and thus from~(\ref{eq:energyest1})
and~(\ref{eq:vazno9}),
for any $t'\leq t$:
\beas
&&\E^{m+1}(t')+\int_0^{t'}\D^{m+1}(\tau)\,d\tau\leq\frac{L}{2}+L^{3/2}+\frac{Ct}{\eta}L+
\eta L+\frac{CtL}{\lambda\epsilon}+\\
&&\big(\sup_{s\leq t}\E^{m+1}(s)+\int_0^t\D^{m+1}(\tau)\,d\tau\big)\Big\{L^{1/2}+C(\eta+\lambda)+
Ct(\frac{1}{\eta}+\frac{1}{\lambda})+
\frac{Ct}{\eta\lambda\epsilon^4}+\frac{C\lambda\epsilon}{\eta}+\frac{Ct}{\lambda\epsilon}+L\Big\}.
\eeas
Since the above inequality holds for any $t'\leq t$, we obtain
\beas
&&\big(\sup_{s\leq t}\E^{m+1}(s)+\int_0^t\D^{m+1}(\tau)\,d\tau\big)\Big\{1-\Big(L^{1/2}+C(\eta+\lambda)+
Ct(\frac{1}{\eta}+\frac{1}{\lambda})+
\frac{Ct}{\eta\lambda\epsilon^4}+\frac{C\lambda\epsilon}{\eta}+\frac{Ct}{\lambda\epsilon}+L\Big)\Big\}\\
&&\leq\frac{L}{2}+CL^{3/2}+\frac{Ct}{\eta}L+
\eta L+\frac{CtL}{\lambda\epsilon}.
\eeas
We first choose $\eta$ small and then $\lambda$ small so that $C(\lambda+\eta)+\frac{C\lambda\epsilon}{\eta}$ is small.
Further, we choose $t$ ($t$ depends on $\epsilon$) and $L$ small so that
\[
L^{1/2}+C(\eta+\lambda)+
Ct(\frac{1}{\eta}+\frac{1}{\lambda})+
\frac{Ct}{\eta\lambda\epsilon^4}+\frac{C\lambda\epsilon}{\eta}+\frac{Ct}{\lambda\epsilon}+L<\frac{1}{3},
\]
and
\[
\frac{L}{2}+C(L)^{3/2}+\frac{Ct}{\eta}L+
\eta L+\frac{CtL}{\lambda\epsilon}\leq\frac{2}{3}L,
\quad L\leq\theta.
\]
With such a choice of $L$ and $t=:{\cal T}^{\epsilon}$, we obtain
$
\sup_{s\leq {\cal T}^{\epsilon}}\E^{m+1}(s)+\int_0^{{\cal T}^{\epsilon}}\D^{m+1}(\tau)\,d\tau\leq L
$
and this finishes the proof of Lemma~\ref{lm:local3}.
\prfe
\section{Regularized Stefan problem.}\label{ch:localwp}
The principal goal of this section is the following local existence theorem:
\begin{theorem}\label{th:local}
For any sufficiently small $L>0$ there exists $t^{\epsilon}>0$ depending on $L$ and $\epsilon$
such that
if for given initial data $(u^{\epsilon}_0,\rho^{\epsilon}_0)$
\[
\E_{\epsilon}(u^{\epsilon}_0,\rho^{\epsilon}_0;\rho^{\epsilon}_0)+\Big|\int_{\T^{n-1}}
\rho^{\epsilon}_0-\int_{\Omega}u^{\epsilon}_0(1+\phi'\rho^{\epsilon}_0)\Big|
\leq\frac{L}{2}
\]
then there exists a unique solution $(u^{\epsilon},\rho^{\epsilon})$ to the regularized Stefan
problem~(\ref{eq:temp})- (\ref{eq:inrho}) and~(\ref{eq:jumpr}) defined
on the time interval $[0,t^{\epsilon}]$. Moreover,
\[
\sup_{0\leq t\leq t^{\epsilon}}\E_{\epsilon}(u^{\epsilon},\rho^{\epsilon};\rho^{\epsilon})(t)+
\int_0^{t^{\epsilon}}\D_{\epsilon}(u^{\epsilon},\rho^{\epsilon};\rho^{\epsilon})(\tau)\,d\tau\leq L
\]
and $\E_{\epsilon}(u^{\epsilon},\rho^{\epsilon})(\cdot)$ is continuous on $[0,t^{\epsilon}[$.
\end{theorem}
\smallskip

\noindent
{\bf Remark.} Note that the constant $L$ is independent of $\epsilon$.\\
\prf
{\bf Convergence.}
Combining Lemmas~\ref{lm:tau} and~{\ref{lm:local3}}, we obtain a uniform-in-$m$ bound
on the sequence $\big\{(u^m,\rho^m)\big\}_m$.
Our goal is to show that $\big\{(u^m,\rho^m)\big\}_m$ is a Cauchy sequence in the energy space.
For any $l\in\N$ let $v^{l+1}:=u^{l+1}-u^l$ and $\sigma^{l+1}=\rho^{l+1}-\rho^l$.
By subtracting two consecutive equations in the iteration process,
we obtain
\be\label{eq:tempd}
v^{m+1}_t-\Delta_{x'}v^{m+1}-\am v^{m+1}_{nn}=f^{\circ}_m,\qquad\qquad\qquad\qquad\qquad\qquad
\qquad\qquad\qquad
\ee
\be
\begin{array}{l}
\displaystyle
v^{m+1}=\Delta\sigma^m\sm^{-1}+g^{\circ}_m\\
\displaystyle
\quad\quad\,\,\,\,=\Delta\sigma^m\sm^{-1}-\rho^m_i\rho^m_j\sigma^m_{ij}\sm^{-1}+G^{\circ}_m\quad\mbox{on}\quad\T^{n-1}\times\{x_n=0\},
\end{array}
\label{eq:gtd}
\ee
\be
\partial_nv^{m+1}=0\quad\quad\textrm{on}\quad\T^{n-1}\times\big\{x_n=\pm1\big\},
\qquad\qquad\qquad\qquad\qquad\qquad\qquad\quad
\ee
\be\label{eq:jumpd}
[v^{m+1}_n]^-_+=\big(\sigma^{m+1}_t+\epsilon\Delta^2\sigma^{m+1}_t\big)\sm^{-2}+h^{\circ}_m\quad\mbox{on}\quad\T^{n-1}\times\{x_n=0\}.
\qquad\qquad
\ee
Here
\be
\begin{array}{l}
\displaystyle
f^{\circ}_m=-\Bm\cdot\nabla_{x'}v^{m+1}_n-\cm v^{m+1}_n+u^m_{nn}(\am-a_{\rho^{m-1}})\qquad\qquad\qquad\quad\\
\displaystyle
-\nabla_{x'}u^m_n\cdot\big(\Bm-B_{\rho^{m-1}}\big)-
u^m_n(\cm-c_{\rho^{m-1}}),
\end{array}
\label{eq:efkrug}
\ee
\beas
G^{\circ}_m&=&\Delta\rho^{m-1}\big(\sm^{-1}-\left\langle \rho^{m-1}\right\rangle\big)+
\rho^{m-1}_{ij}\big(\sigma^m_i\rho^m_j\sm^{-1}+\rho^{m-1}_i\sigma^m_j\sm^{-1}+\\
&&\quad\rho^{m-1}_i\rho^{m-1}_j(\sm^{-1}-\left\langle \rho^{m-1}\right\rangle)\big),
\eeas
\[
g^{\circ}_m=-\rho^m_i\rho^m_j\sigma^m_{ij}\sm^{-1}+G^{\circ}_m,\qquad\qquad\qquad\qquad\qquad\qquad\qquad
\qquad\qquad\quad
\]
\be\label{eq:hakrug}
h^{\circ}_m=-\big(|\nabla\rho^m|^2-|\nabla\rho^{m-1}|^2\big)[u^m_n]^-_+\sm^{-2}.\quad
\qquad\qquad\qquad\qquad\qquad
\ee
After applying the differential operator $\d$ to the equations~(\ref{eq:tempd}), (\ref{eq:gtd})
and~(\ref{eq:jumpd}) and singling out
the leading-order terms we arrive at:
\be\label{eq:tempd1}
\d v^{m+1}_t-\Delta_{x'}v^{m+1}-\am v^{m+1}_{nn}=f'_m,\qquad\qquad\qquad\qquad
\qquad\qquad\quad
\ee
\be\label{eq:gtd1}
\d v^{m+1}=\d\Delta\sigma^m\sm^{-1}+g'_m
=\d\Delta\sigma^m\sm^{-1}-\rho^m_i\rho^m_j\d\sigma^m_{ij}\sm^{-1}+G'_m,
\ee
\be
\label{eq:jumpd1}
[\d v^{m+1}_n]^-_+=\big(\d\sigma^{m+1}_t+\epsilon\d\Delta^2\sigma^{m+1}_t\big)\sm^{-2}+h'_m,
\qquad\qquad\qquad\qquad\qquad
\ee
where
\be\label{eq:efem}
f'_m=\d f^{\circ}_m
+\big(\d(\am\d v^{m+1}_{nn})-\am\d v^{m+1}_{nn}\big),
\ee
\be\label{eq:gem}
g'_m=\d g^{\circ}_m+\sum_{|\mu'|+s'\atop <|\mu|+s}\partial^{\mu'}_{s'}\Delta\sigma^m
\partial^{\mu-\mu'}_{s-s'}(\sm^{-1}),
\ee
\be
\begin{array}{l}
\displaystyle
G'_m=\d G^{\circ}_m+\sum_{|\mu'|+s'\atop <|\mu|+s}\partial^{\mu'}_{s'}\Delta\sigma^m
\partial^{\mu-\mu'}_{s-s'}(\sm^{-1})-\\
\displaystyle
\Big(\d\big(\rho^m_i\rho^m_j\sigma^m_{ij}\sm^{-1}\big)-
\rho^m_i\rho^m_j\d\sigma^m_{ij}\sm^{-1}\Big),
\end{array}
\label{eq:Gem}
\ee
\be\label{eq:haem}
h'_m=\d h^{\circ}_m+
\sum_{|\mu'|+s'\atop <|\mu|+s}\partial^{\mu'}_{s'}\big(\sigma^{m+1}_t+\epsilon\Delta^2\sigma^{m+1}_t\big)
\partial^{\mu-\mu'}_{s-s'}(\sm^{-2}).
\ee
As a next step, we use the identities from Chapter~\ref{ch:basic} to obtain the energy identities for
the problem~(\ref{eq:tempd1}) -~(\ref{eq:jumpd1}). Respecting the notations of Chapter~\ref{ch:basic} we
set for any $l\in\N$
\be\label{eq:subst}
f=f'_l,\,g=g'_l,\,\, G=G'_l,\,\, h=h'_l,\,\,
\U=\d v^{l+1},\,\, \omega=\d\sigma^{l+1},\,\, \chi=\d\sigma^l\,\, \mbox{and}\,\, \psi=\rho^l.
\ee
Additionally, we introduce the notations
\[
e^l:=\E_{\epsilon}(v^l,\sigma^l;\rho^{l-1}),\quad
d^l:=\D_{\epsilon}(v^l,\sigma^l;\rho^{l-1}),
\]
where $\E_{\epsilon}$ and $\D_{\epsilon}$ are defined by~(\ref{eq:energyeps}) and~(\ref{eq:dissipationeps}) respectively.
Using~(\ref{eq:mainenergy}) and~(\ref{eq:subst}), we arrive at
\be\label{eq:energylittle}
\frac{d}{dt}e^{m+1}+d^{m+1}=\int_{\Omega}\big\{p^m+r^m\big\}-\int_{\T^{n-1}}\big\{q^m+s^m+t^m\big\},
\ee
where
\be
\begin{array}{l}
\displaystyle
p^m:=\sum_{|\mu|+2s\leq 2k}P(\d v^{m+1},\rho^m,f'_m),\quad r^m:=\sum_{|\mu|+2s\leq 2k}R(\d v^{m+1},\rho^m,f'_m)\\
\displaystyle
q^m:=\sum_{|\mu|+2s\leq 2k}Q(\d\sigma^m,\d\sigma^{m+1},\rho^m,g'_m,h'_m),\\
\displaystyle
s^m:=\sum_{|\mu|+2s\leq 2k}S(\d\sigma^m,\d\sigma^{m+1},\rho^m,g'_m,h'_m),\\
\displaystyle
t^m:=\sum_{|\mu|+2s\leq 2k}T(\d\sigma^m,\d\sigma^{m+1},\rho^m,g'_m,h'_m).
\end{array}
\label{eq:defns}
\ee
 Here
$P$, $Q$, $R$, $S$ and $T$ are defined by~(\ref{eq:P}), ~(\ref{eq:Q}), ~(\ref{eq:R}), ~(\ref{eq:S})
and~(\ref{eq:bigT}) respectively. Our aim is to prove that for suitably small $t\leq t(\epsilon)$
there exists a $\Lambda<1$ such that
\[
e^{m+1}(t)+\int_0^td^{m+1}(\tau)\,d\tau\leq \Lambda\big(e^m(t)+\int_0^td^m(\tau)\,d\tau\big).
\]
We shall accomplish this by estimating the terms $p^m$, $r^m$,
$q^m$, $s^m$ and $t^m$ on RHS of the
identity~(\ref{eq:energylittle}). These estimates will be largely
analogous to the estimates from Chapter~\ref{ch:energye}. However,
due to the formally new terms $\d f^{\circ}_m$, $\d g^{\circ}_m$,
$\d G^{\circ}_m$, $\d h^{\circ}_m$ appearing in the
definitions~(\ref{eq:efem}), ~(\ref{eq:gem}), ~(\ref{eq:Gem})
and~(\ref{eq:haem}) of $f'_m$, $g'_m$, $G'_m$ and $h'_m$
respectively, we need to make several preparatory steps. First, for
any $l\in\N$ \be\label{eq:w1}
||\am-a_{\rho^{m-1}}||_{H^l(\Omega)}+||\Bm-B_{\rho^{m-1}}||_{H^l(\Omega)}
\leq C||\sigma^m||_{H^{l+1}(\Omega)}, \ee and \be\label{eq:w2}
||\cm-c_{\rho^{m-1}}||_{H^l(\Omega)}\leq
C||\sigma^m_t||_{H^l}+C||\sigma^m||_{H^{l+2}(\Omega)}. \ee Note that
$ \sup_{0\leq t\leq
t^{\epsilon}}\E^m(t)+\int_0^{t^{\epsilon}}\D^m(\tau)\,d\tau\leq L. $
In particular, for $|\mu|+2s\leq 2k$ \be\label{eq:uniformb} ||\d
u^m||_{H^1(\Omega)}\,\,,||\d \am||_{L^2(\Omega)}\,\,, ||\d
\Bm||_{L^2(\Omega)}\leq C\sqrt{L}. \ee Furthermore, $ ||\d
\cm||_{L^2(\Omega)}\leq C\sqrt{L} $ for $|\mu|+2s\leq 2k-1$. We can
now use the Sobolev inequality to bound the lower order derivatives
of $u^m$, $\am$, $\Bm$ and $\cm$ in $L^{\infty}$-norm by
$C\sqrt{L}$. The major step is to provide the analogues of part (c)
of Lemma~\ref{lm:aux} for the function $\sigma^{m+1}$ instead of
$\rho^{m+1}$ and Lemma~\ref{lm:auximp} for the function $v^{m+1}$
instead of $u^{m+1}$. By the boundary condition~(\ref{eq:jumpd1})
and the proof of Lemma~\ref{lm:aux}, part (c), we deduce
\be\label{eq:d1}
||\d\sigma^{m+1}_t||_2^2+2\epsilon||\d\Delta\sigma^{m+1}_t||_2^2+\epsilon^2||\d\Delta^2\sigma^{m+1}_t||_2^2\leq
\lambda d^{m+1}+\frac{C}{\lambda}e^{m+1}+e^m(\E^m+\D^m) \ee and
\be\label{eq:d2}
||\d\sigma^{m+1}_t||_2^2+2\epsilon||\d\Delta\sigma^{m+1}_t||_2^2+\epsilon^2||\d\Delta^2\sigma^{m+1}_t||_2^2\leq
Cd^{m+1}+e^m\D^m. \ee As in Lemma~\ref{lm:auximp}, \be\label{eq:d3}
||\d\partial_{n^r}v^{m+1}||_{L^2(\Omega)}^2\leq Ce^{m+1}+Ce^m. \ee
The proof of~(\ref{eq:d3}) is completely analogous to the proof of
Lemma~\ref{lm:auximp} whereby, due to the addition of the formally
new term $\d f^{\circ}_m$ in the definition of $f'_m$ we need to
exploit the relations~(\ref{eq:w1}) and~(\ref{eq:w2}), which are
responsible for the occurrence of the term $Ce^m$ on RHS
of~(\ref{eq:d3}). We proceed fully analogously to the energy
estimates in Chapter~\ref{ch:energye} to estimate the right-hand
side of the energy identity~(\ref{eq:energylittle}). The terms
involving $\d f^{\circ}_m$ and $\d h^{\circ}_m$ require an
additional care. In the estimate for the term
$\int_{\Omega}p^m=\int_{\Omega}f'_m\d v^{m+1}$, where $f'_m$ is
given by~(\ref{eq:efem}), we single out the term $\int_{\Omega}\d
f^{\circ}_m\d v^{m+1}$. Here $f^{\circ}_m$ is given
by~(\ref{eq:efkrug}). Writing \be\label{eq:splitting} \d(\cm
v^{m+1}_n)=\sum_{\mu',s'}C^{\mu'}_{s'}\partial^{\mu'}_{s'}\cm\partial^{\mu-\mu'}_{s-s'}v^{m+1}_n=
\sum_{\mu'+2s'\leq k}+\sum_{\mu'+2s'> k} \ee we estimate the
lower-order terms in both sums in $L^{\infty}$-norm and the
higher-order terms in $L^2$-norm. By the Cauchy-Schwarz inequality
and~(\ref{eq:d3}), \be\label{eq:vazno7} |\int_{\Omega}\d(\cm
v^{m+1}_n)\d v^{m+1}|\leq
C\sqrt{\D^m}\sqrt{e^{m+1}}(\sqrt{e^m}+\sqrt{e^{m+1}}). \ee
Similarly, $ |\int_{\Omega}\d(\Bm\cdot\nabla_{x'}v^{m+1}_n)\d
v^{m+1}|\leq C\sqrt{\E^m}(e^{m+1}+d^{m+1}). $ In analogous fashion,
we find \be
\begin{array}{l}
\displaystyle
\Big|\int_{\Omega}\d\big(u^m_{nn}(\am-a_{\rho^{m-1}})-\nabla_{x'}u^m_n\cdot\big(\Bm-B_{\rho^{m-1}}\big)-
u^m_n(\cm-c_{\rho^{m-1}})\big)\d v^{m+1}\Big|\\
\displaystyle
\quad
\leq C\sqrt{\D^m}\sqrt{e^m}\sqrt{e^{m+1}}+C\sqrt{\E^m}\sqrt{d^m}\sqrt{e^{m+1}},
\end{array}
\label{eq:vazno8}
\ee
where we rely on~(\ref{eq:w1}) and~(\ref{eq:w2}). The term $r^m$ is of the
form $r^m=(f'_m)^2+(\textrm{rest})$, ($r^m$ is defined
in~(\ref{eq:defns})). The formally new term to estimate has the form $\int_{\Omega}(\d f^{\circ}_m)^2$, where
$f^{\circ}_m$ is given by~(\ref{eq:efkrug}). By the same splitting idea
as in~(\ref{eq:splitting}) and the estimate~(\ref{eq:d3}),
\[
|\int_{\Omega}(\d f^{\circ}_m)^2|\leq C\E^m(e^m+e^{m+1}+d^{m+1})+C\D^me^{m+1}+C\D^me^m+C\E^md^m.
\]
Recall now that $s^m$ and $t^m$ are defined by~(\ref{eq:defns}). The
last term in each of the expressions~(\ref{eq:S}) and~(\ref{eq:bigT})
has the form $\sm^2h'_mv^{m+1}_t|_{\T^{n-1}}$ and
$\sm^2h'_m\Delta_{x'}v^{m+1}|_{\T^{n-1}}$, respectively.
By~(\ref{eq:haem}), the formally new term (with respect to $h^m_{\mu,s}$
defined by~(\ref{eq:ha})) is $\d h^{\circ}_m$.
By~(\ref{eq:hakrug}), we conclude
\[
\d h^{\circ}_m=-\sum_{\mu',s'}C^{\mu'}_{s'}\partial^{\mu'}_{s'}([u^m_n]^-_+)
\partial^{\mu-\mu'}_{s-s'}\big((|\nabla\rho^m|^2-|\nabla\rho^{m-1}|^2)\sm^{-2}\big)=
\sum_{\mu'+2s'\leq k}+\sum_{\mu'+2s'> k}
\]
Hitting the lower order terms with $L^{\infty}$-norm and the higher order terms
with $L^2$-norm and using the trace inequality to estimate $||\partial^{\mu'}_{s'}[u^m_n]^-_+||_2$,
we obtain
\be\label{eq:haest}
||\d h^{\circ}_m||_2\leq(\lambda\sqrt{\D^m}+\frac{C}{\lambda}\sqrt{\E^m})\sqrt{e^{m}},\quad\quad
||\d h^{\circ}_m||_2\leq C\sqrt{\D^m}\sqrt{e^{m}}.
\ee
Note that $v^{m+1}_t|_{\T^{n-1}}=(\kappa_{\rho^m}-\kappa_{\rho^{m-1}})_t$. It is easy to check
that
\[
||\d v^{m+1}_t|_{\T^{n-1}}||_2\leq C(||\nabla^2\sigma^m_t||_2+||\nabla\sigma^m_t||_2+
||\nabla^2\sigma^m||_2+||\nabla\sigma^m||_2).
\]
By the definitions of $d^{m+1}$ and $e^{m+1}$,
$
||\d v^{m+1}_t|_{\T^{n-1}}||_2\leq\frac{C}{\sqrt{\epsilon}}\sqrt{d^m}+C\sqrt{e^m}.
$
Combining this with~(\ref{eq:haest}) and the Cauchy-Schwarz inequality,
we obtain
\[
\Big|\int_{\T^{n-1}}\sm^2\d h^{\circ}_mv^{m+1}_t|_{\T^{n-1}}\Big|\leq
(\lambda\sqrt{\D^m}+\frac{C}{\lambda}\sqrt{\E^m})\sqrt{e^{m}}(\frac{C}{\sqrt{\epsilon}}\sqrt{d^m}+C\sqrt{e^m}).
\]
Choosing $\lambda=\sqrt{\epsilon}$ we arrive at
\be
\begin{array}{l}
\displaystyle
\Big|\int_{\T^{n-1}}\sm^2\d h^{\circ}_mv^{m+1}_t|_{\T^{n-1}}\Big|\leq
C\sqrt{\D^m}\sqrt{e^m}\sqrt{d^m}+C\sqrt{\epsilon}\sqrt{\D^m}e^{m}+\\
\displaystyle
\quad\quad\frac{C}{\epsilon}\sqrt{\E^m}\sqrt{e^{m}}\sqrt{d^m}+
\frac{C}{\sqrt{\epsilon}}\sqrt{\E^m}e^m.
\end{array}
\label{eq:vazno1}
\ee
We want to estimate the time integral of the right-hand side of
the inequality~(\ref{eq:vazno1}). For any $t\leq t^{\epsilon},\lambda>0$, we have
\be\label{eq:vazno3}
\int_0^t\sqrt{\D^m}\sqrt{e^m}\sqrt{d^m}\,d\tau\leq\lambda\int_0^td^m(\tau)\,d\tau+
\frac{C}{\lambda}\sup_{0\leq s\leq t}e^m(s)\int_0^t\D^m(\tau)\,d\tau.
\ee
Similarly,
\be\label{eq:vazno4}
\int_0^t\sqrt{\D^m}e^m\,d\tau\leq\int_0^te^m+\sup_{0\leq s\leq t}e^m(s)\int_0^t\D^m(\tau)\,d\tau\leq
t\sup_{0\leq s\leq t}e^m(s)+L\sup_{0\leq s\leq t}e^m(s).
\ee
Since $\E^m\leq L$, for any $\lambda>0$, we obtain
\be\label{eq:vazno5}
\int_0^t\frac{C}{\epsilon}\sqrt{\E^m}\sqrt{e^{m}}\sqrt{d^m}\,d\tau\leq
\lambda\int_0^td^m(\tau)\,d\tau+\frac{C}{\lambda\epsilon^2}Lt\sup_{0\leq s\leq t}e^m(s),
\ee
and similarly,
\be\label{eq:vazno6}
\int_0^t\frac{C}{\sqrt{\epsilon}}\sqrt{\E^m}e^m\,d\tau\leq\frac{C}{\sqrt{\epsilon}}t\sqrt{L}
\sup_{0\leq s\leq t}e^m(s).
\ee
Noting that $||\Delta_{x'}v^{m+1}|_{\T^{n-1}}||_2\leq\frac{C}{\sqrt{\epsilon}}\sqrt{e^m}$, we can use
the Cauchy-Schwarz inequality and the estimate~(\ref{eq:haest}) to get
\be\label{eq:vazno2}
\Big|\int_{\T^{n-1}}\sm^2\d h^{\circ}_m\Delta_{x'}v^{m+1}|_{\T^{n-1}}\Big|\leq
\frac{C}{\sqrt{\epsilon}}\sqrt{\D^m}e^m.
\ee
We get
\be\label{eq:vazno10}
\int_0^t\frac{C}{\sqrt{\epsilon}}\sqrt{\D^m}e^m\,d\tau\leq\frac{C}{\epsilon}t\sup_{0\leq s\leq t}e^m(s)+
L\sup_{0\leq s\leq t}e^m(s).
\ee
In analogy to the estimate~(\ref{eq:energyest}) together with~(\ref{eq:vazno7}), ~(\ref{eq:vazno8}),
~(\ref{eq:vazno1}) and ~(\ref{eq:vazno2}), we obtain
\[
\begin{array}{l}
\displaystyle
\frac{d}{dt}e^{m+1}+d^{m+1}\leq Ce^m(\frac{1}{\eta}+\frac{C}{\lambda\epsilon})+
\eta d^m+Ce^{m+1}(\frac{1}{\eta}+\frac{1}{\lambda}+\frac{C}{\eta\lambda\epsilon^4}+\frac{C}{\lambda\epsilon})+\\
\displaystyle
Cd^{m+1}(\sqrt{L}+\eta+\lambda+\frac{C\lambda\epsilon}{\eta})+C\sqrt{\D^m}\sqrt{d^{m+1}}\sqrt{e^{m+1}}+
C\D^me^{m+1}+\\
\displaystyle
C\sqrt{\D^{m+1}}\sqrt{e^{m+1}}(\sqrt{e^m}+\sqrt{e^{m+1}})+
C\sqrt{\D^m}\sqrt{e^m}\sqrt{e^{m+1}}+C\sqrt{\E^m}\sqrt{d^m}\sqrt{e^{m+1}}+\\
\displaystyle
C\sqrt{\D^m}\sqrt{e^m}\sqrt{d^m}+C\sqrt{\epsilon}\sqrt{\D^m}e^{m}+
\frac{C}{\epsilon}\sqrt{\E^m}\sqrt{e^{m}}\sqrt{d^m}+
\frac{C}{\sqrt{\epsilon}}\sqrt{\E^m}e^m+\frac{C}{\sqrt{\epsilon}}\sqrt{\D^m}e^m.
\end{array}
\]
Integrating the above inequality over $[0,t]$ and proceeding as in~(\ref{eq:energyest1}),
using~(\ref{eq:vazno3}), ~(\ref{eq:vazno4}), ~(\ref{eq:vazno5}), ~(\ref{eq:vazno6}), ~(\ref{eq:vazno10}),  we conclude
\be
\begin{array}{l}
\displaystyle
e^{m+1}(t)+\int_0^td^{m+1}(\tau)\,d\tau\leq
Ct(\frac{1}{\eta}+\frac{C}{\lambda\epsilon})\sup_{0\leq s\leq t}e^m(s)+
\eta\int_0^td^m(\tau)\,d\tau+\\
\displaystyle
Ct(\frac{1}{\eta}+\frac{1}{\lambda}+\frac{C}{\eta\lambda\epsilon^4}+\frac{C}{\lambda\epsilon})
\sup_{0\leq s\leq t}e^{m+1}(s)
+C(\sqrt{L}+\eta+\lambda+\frac{C\lambda\epsilon}{\eta})\int_0^td^{m+1}(\tau)\,d\tau\\
\displaystyle
+C\lambda\big(\int_0^td^m(\tau)\,d\tau+\int_0^td^{m+1}(\tau)\,d\tau\big)+\\
\displaystyle
(\frac{CL}{\lambda}+L+t+\frac{CL}{\lambda\epsilon^2}t+\frac{C\sqrt{L}}{\sqrt{\epsilon}}t+\frac{C}{\epsilon}t)
\sup_{0\leq s\leq t}e^m(s)+
C(t+L+\sqrt{L})
\sup_{0\leq s\leq t}e^{m+1}(s)
\end{array}
\ee
We choose $\eta$, $\lambda$, $L$ and  $t=:t^{\epsilon}\leq{\cal T}^{\epsilon}$ small, so that $\eta+C\lambda<\chi<\frac{1}{5}$
\[
Ct(\frac{1}{\eta}+\frac{C}{\lambda\epsilon})
+\frac{CL}{\lambda}+L+t^{\epsilon}+\frac{CL}{\lambda\epsilon^2}t^{\epsilon}+\frac{C\sqrt{L}}{\sqrt{\epsilon}}t^{\epsilon}+
\frac{C}{\epsilon}t^{\epsilon}
<\chi<\frac{1}{5},
\]
\[
Ct^{\epsilon}(\frac{1}{\eta}+\frac{1}{\lambda}+\frac{C}{\eta\lambda\epsilon^4}+\frac{C}{\lambda\epsilon})+
C(t^{\epsilon}+L+\sqrt{L})<\chi<\frac{1}{5}
\]
and $C(\sqrt{L}+\eta+\lambda+\frac{C\lambda\epsilon}{\eta})<\chi<\frac{1}{5}$.
Taking supremum over $[0,t^{\epsilon}]$, we arrive at
\beas
&&\sup_{0\leq \tau\leq t^{\epsilon}}e^{m+1}(\tau)+\int_0^{t^{\epsilon}}d^{m+1}(s)\,ds\leq\chi
\Big\{\sup_{0\leq \tau\leq t^{\epsilon}}e^m(\tau)+\int_0^{t^{\epsilon}}d^m(s)\,ds\Big\}+\\
&&+\chi\Big\{\sup_{0\leq\tau\leq t^{\epsilon}}e^{m+1}(\tau)+\int_0^{t^{\epsilon}}d^{m+1}(s)\Big\}\,ds.
\eeas
Therefore,
\be\label{eq:almost}
\sup_{0\leq \tau\leq t^{\epsilon}}e^{m+1}(\tau)+\int_0^{t^{\epsilon}}d^{m+1}(s)\,ds\leq
\frac{\chi}{1-\chi}\Big\{\sup_{0\leq \tau\leq t^{\epsilon}}e^m(\tau)+\int_0^{t^{\epsilon}}d^m(s)\,ds\Big\}.
\ee
We observe now that the conservation law~(\ref{eq:consm}) and the fact that
$\sigma^{m+1}(0)=v^{m+1}(0)=0$ imply
\[
\int_{\T^{n-1}}\sigma^{m+1}=\int_{\Omega}v^{m+1}(1+\phi'\rho^m)+\int_{\Omega}\phi'u^m\sigma^m.
\]
By the Poincar\'e inequality, previous identity and the uniform bounds on $\rho^m$ and $u^m$ we get
\be
\begin{array}{l}
\displaystyle
||\sigma^{m+1}||_2^2\leq C||\nabla\sigma^{m+1}||_2^2+C||\int_{\T^{n-1}}\sigma^{m+1}||_2^2\\
\displaystyle
\leq C||\nabla\sigma^{m+1}||_2^2+C||v^{m+1}||_2^2+C||u^m||_{L^2(\Omega)}^2||\sigma^m||_2^2\\
\displaystyle
\leq Ce^{m+1}+C\E^m||\sigma^m||_2^2\leq Ce^{m+1}+CL||\sigma^m||_2^2.
\end{array}
\label{eq:controlm}
\ee
With $L$ and $\chi$ so small that $CL+2C\frac{\chi}{1-\chi}=:\Lambda<1$, we obtain by~(\ref{eq:almost}),
\beas
&&\sup_{0\leq\tau\leq t^{\epsilon}}||\sigma^{m+1}(\tau)||_2^2\leq
C\frac{\chi}{1-\chi}\Big\{\sup_{0\leq \tau\leq t^{\epsilon}}e^m(\tau)+\int_0^{t^{\epsilon}}d^m(s)\,ds\Big\}+\\
&&+CL\sup_{0\leq \tau\leq t^{\epsilon}}||\sigma^m(\tau)||_2^2
\leq\frac{\Lambda}{2}\Big\{\sup_{0\leq \tau\leq t^{\epsilon}}e^m(\tau)+\int_0^{t^{\epsilon}}d^m(s)\,ds\Big\}+
\Lambda\sup_{0\leq\tau\leq t^{\epsilon}}||\sigma^m(\tau)||_2^2.
\eeas
Adding  $\big\{\sup_{0\leq \tau\leq t^{\epsilon}}e^{m+1}(\tau)+\int_0^{t^{\epsilon}}d^{m+1}(s)\,ds\big\}$ to the both sides of the above inequality and
using~(\ref{eq:almost}) again we get
\beas
&&\sup_{0\leq\tau\leq t^{\epsilon}}\Big\{e^{m+1}(\tau)+\int_0^{t^{\epsilon}}d^{m+1}(s)\,ds+||\sigma^{m+1}(\tau)||_2^2\Big\}\leq\\
&&C\frac{\chi}{1-\chi}\Big\{\sup_{0\leq \tau\leq t^{\epsilon}}e^m(\tau)+\int_0^{t^{\epsilon}}d^m(s)\,ds\Big\}+
\frac{\Lambda}{2}\Big\{\sup_{0\leq \tau\leq t^{\epsilon}}e^m(\tau)+\int_0^{t^{\epsilon}}d^m(s)\,ds\Big\}+\\
&&\sup_{0\leq\tau\leq t^{\epsilon}}\Lambda||\sigma^m(\tau)||_2^2
\leq\Lambda\sup_{0\leq \tau\leq t^{\epsilon}}\Big\{e^m(\tau)+\int_0^{t^{\epsilon}}d^m(s)\,ds
+||\sigma^m(\tau)||_2^2\Big\}.
\eeas
Define the Banach spaces
\[
X:=\Big\{(v,\sigma)\Big|\,\,|(v,\sigma)||_{\E_{\epsilon}}+||\sigma||_2^2<\infty\Big\},\quad
Y:=\Big\{(v,\sigma)\Big|\,\,||(v,\sigma)||_{\D_{\epsilon}}<\infty\Big\},
\]
where $||(\cdot,\cdot)||_{\E_{\epsilon}}$ and $||(\cdot,\cdot)||_{\D_{\epsilon}}$ are
defined by~(\ref{eq:equivE}) and~(\ref{eq:equivD}) respectively.
Let
$
Z:=L^{\infty}\big(X;[0,t^{\epsilon}]\big)\oplus L^2\big(Y;[0,t^{\epsilon}]\big).
$
Since $\Lambda$ can be chosen arbitrarily small,
we have proven that the sequence $\big((v^l,\sigma^l)\big)_{l\in\N}$ satisfies
$
||(v^{m+1},\sigma^{m+1})||_Z\leq\Lambda'||(v^m,\sigma^m)||_Z,
$
for some $\Lambda'<1$. This implies that $\big((u^l,\rho^l)\big)_{l\in\N}$ is a Cauchy sequence
and converges strongly in $Z$. Thus the {\em whole} sequence $\big((u^l,\rho^l)\big)_{l\in\N}$
converges to a solution $(u^{\epsilon},\rho^{\epsilon})$ of the regularized Stefan problem in
the original energy space.
In addition to this, passing to a limit in~(\ref{eq:consm}),
we obtain the conservation law
\be\label{eq:cons}
\partial_t\Big\{\int_{\Omega} u(1+\phi'\rho)\Big\}=\partial_t\Big\{\int_{\T^{n-1}}\rho\Big\}.
\ee
\par
{\bf Uniqueness.}
 We want to prove uniqueness in the class of functions $(u,\rho)$ satisfying
$\sup_{0\leq t<t^{\epsilon}}\E_{\epsilon}(u,\rho)(t)+\int_0^{t^{\epsilon}}\D_{\epsilon}(\tau)\,d\tau\leq L$,
where $L$ may be chosen smaller if necessary.
Let us assume that there exists another solution $(v,\sigma)$ satisfying the same
initial conditions $(v(x,0),\sigma(x',0))=(u_0(x),\rho_0(x'))$ and the bound
$\sup_{0\leq t<t^{\epsilon}}\E_{\epsilon}(v,\sigma)(t)+
\int_0^{t^{\epsilon}}\D_{\epsilon}(v,\sigma)(\tau)\,d\tau\leq L$.
After subtracting them
and setting $w:=u-v$, $\tau:=\rho-\sigma$, we obtain
\be
w_t-\Delta_{x'}w-\a w_{nn}=f^{\ast},\qquad\qquad\qquad\qquad\qquad
\qquad\qquad\qquad\qquad\qquad\quad
\ee
\be\label{eq:uniq1e'}
w=\Delta\tau\s^{-1}+g^{\ast}
=\Delta\tau\s^{-1}-\rho_i\rho_j\tau_{ij}\s^{-1}+G^{\ast}\quad\mbox{on}\quad\T^{n-1}\times\{x_n=0\},
\ee
\be\label{eq:uniq2e'}
[w_n]^-_+=\big(\tau_t+\epsilon\Delta^2\tau_t\big)\s^{-2}+h^{\ast}\quad\mbox{on}\quad\T^{n-1}\times\{x_n=0\},
\qquad\qquad\quad
\ee
where
\beas
&&f^{\ast}=-\B\cdot\nabla_{x'}w_n-\c w_n+(\a-\as)v_{nn}+(\B-\Bs)\nabla_{x'}v_n+(\cs-\c)v_n,\\
&&G^{\ast}=\Delta\sigma\big(\s^{-1}-\si^{-1}\big)+
\sigma_{ij}\big(\tau_i\rho_j\s^{-1}+\sigma_i\tau_j\s^{-1}+
\sigma_i\sigma_j\big(\s^{-1}-\si^{-1}\big)\big),\\
&&g^{\ast}=-\rho_i\rho_j\tau_{ij}\s^{-1}+G^{\ast}\\
&&h^{\ast}=-(|\nabla\rho|^2-|\nabla\sigma|^2)[v_n]^-_+\s^{-2}.
\eeas
We use
Chapter~\ref{ch:basic} to derive the accompanying energy identities.
To this end
we set $\U=w$, $\psi=\rho$, $\omega=\chi=\tau$,
$f=f^{\ast}$, $g=g^{\ast}$, $G=G^{\ast}$ and $h=h^{\ast}$. Here $\rho$ takes
the role of $\rho^{m+1}$ and $\sigma$ the role of $\rho^m$ and additionally, the cross-terms vanish
since $\omega=\chi=\tau$. With $k\geq n$ sufficiently large, the regularity assumptions of Lemma~\ref{lm:model}
are fulfilled.
We are thus naturally led to the following energy quantities:
\[
\E^{\ast}:=\E_{\epsilon}(w,\tau;\rho),\quad\D^{\ast}:=\D(w,\tau;\rho).
\]
In addition to this we define $P^{\ast}=P(w,\rho,f^{\ast})$ and
analogously $Q^{\ast}$, $R^{\ast}$, $S^{\ast}$ and $T^{\ast}$. Using
the identity~(\ref{eq:mainenergy}), we obtain \be\label{eq:energyw}
\frac{d}{dt}\E^{\ast}+\D^{\ast}=\int_{\Omega}\Big\{P^{\ast}+R^{\ast}\Big\}
-\int_{\T^{n-1}}\Big\{Q^{\ast}+S^{\ast}+T^{\ast}\Big\}. \ee Our goal
at this stage is to prove the inequality of the form
\be\label{eq:gronwallgentle}
\frac{d}{dt}\E^{\ast}(t)+\D^{\ast}(t)\leq C\E^{\ast}(t)+C\sqrt{L}
\D^{\ast}(t), \ee which would enable us to absorb the multiple of
$\D^{\ast}$ on the right-hand side into the left-hand side and then
use the Gronwall's inequality to conclude that $\E^{\ast}(t)=0$ for
any $t\geq0$. It is essential that the constant $C$ in the above
estimate does not depend on $\epsilon$ so that the smallness bound
on $L$ remains independent of $\epsilon$. That the
identity~(\ref{eq:gronwallgentle}) indeed holds, follows analogously
to the energy estimates from the Chapter~\ref{ch:energye} applied to
the right-hand side of~(\ref{eq:energyw}). Here we strongly exploit
the uniform bounds on $\E_{\epsilon}(u,\rho)$ and
$\E_{\epsilon}(v,\sigma)$. In particular we know that
\[
||(\a)_t||_{\infty},\,\,||\nabla_{x'}(\a)||_{\infty},\,\,||(\a)_n||_{\infty},\,\,
||\B||_{\infty},\,\,||\c||_{\infty},\,\,||\Bs||_{\infty},\,\,||\cs||_{\infty}\,\,\leq C\sqrt{L},
\]
\[
||v_n||_{L^{\infty}(\Omega)},\,\,||\nabla_{x'}v_n||_{L^{\infty}(\Omega)},\,\,
||v_{nn}||_{L^{\infty}(\Omega)}\,\,\leq C\sqrt{L}
\]
A major difference from the existence part of the proof is the absence
of cross-terms in the energy identities (since $\omega=\chi$ in the notation of Chapter~\ref{ch:basic}).
In addition to that, we work in a lower order energy space and we can thus use the above uniform
estimates to bound the term $[v_n]^-_+$ by $C\sqrt{L}$ in $L^{\infty}$-norm. This
observation is crucial when estimating $h^{\ast}$.
Knowing that the $\epsilon$-dependence comes only from the estimates of the cross-terms
(cf.~(\ref{eq:z16}), ~(\ref{eq:z17}), ~(\ref{eq:z23}), ~(\ref{eq:z24}), ~(\ref{eq:z26}) and~(\ref{eq:z27})),
we conclude that the constants on the right-hand side
of~(\ref{eq:gronwallgentle}) {\em do not} depend on $\epsilon$.
Choosing $L$ suitably small~(\ref{eq:gronwallgentle}) implies
$
\frac{d}{dt}\E^{\ast}\leq C\E^{\ast}
$
implying
$
\E^{\ast}(t)\leq C\int_0^t\E^{\ast}(s)\,ds,
$
since $\E^{\ast}(0)=0$.
By Gronwall's inequality, we conclude $\E^{\ast}(t)=0$.
In addition to this
the conservation law~(\ref{eq:cons}) gives
$
||\tau||_2^2\leq C\E^{\ast}.
$
This estimate follows in the same way as~(\ref{eq:controlm}).
Thus $(u,\rho)=(v,\sigma)$. This finishes the proof of the uniqueness claim.
\par
{\bf Continuity.}
Integrating the identity~(\ref{eq:energymm}) over the time interval $[s,t]$, we obtain
\be\label{eq:cont}
\E^{m+1}(t)-\E^{m+1}(s)+\int_s^t\D^{m+1}(s)=\int_s^t\int_{\Omega}P^m+R^m-
\int_s^t\int_{\T^{n-1}}Q^m+S^m+T^m.
\ee
However, since $(u^l,\rho^l)\to(u,\rho)$ strongly in the energy space, we may pass to the limit in~(\ref{eq:cont}) to
conclude
\be\label{eq:cont2}
\E_{\epsilon}(t)-\E_{\epsilon}(s)+
\int_s^t\D_{\epsilon}(\tau)\,d\tau
=\int_s^t\int_{\Omega}\bar{P}+\bar{R}-\int_s^t\int_{\T^{n-1}}
\bar{Q}+\bar{S}+\bar{T}.
\ee
Here $\bar{P}=\sum_{|\mu|+2s\leq2k}P(\d u,\rho, f_{\mu,s})$,
where $f_{\mu,s}$ is defined by dropping the index $m$ in the
definition~(\ref{eq:ef}) of $f^m_{\mu,s}$. The terms $\bar{Q}$, $\bar{R}$,
$\bar{S}$ and $\bar{T}$ are defined analogously. We claim that
\be\label{eq:cont3}
\Big|\int_s^t\int_{\Omega}\bar{P}+\bar{R}-\int_s^t\int_{\T^{n-1}}
\bar{Q}+\bar{S}+\bar{T}\Big|\leq C\int_s^t\sqrt{\E_{\epsilon}(\tau)}
\D_{\epsilon}(\tau)\,d\tau.
\ee
The inequality follows easily from the energy estimates in Chapter~\ref{ch:energye}.
We observe that the estimates involving $\epsilon$ on the right-hand side, are used
only when estimating the cross-terms
(cf.~(\ref{eq:z16}), ~(\ref{eq:z17}), ~(\ref{eq:z23}), ~(\ref{eq:z24}), ~(\ref{eq:z26}) and~(\ref{eq:z27})). However, the cross-terms vanish
as $m$ goes to $\infty$ (since $\chi=\omega=\d\rho$). As a result, we obtain
the estimate~(\ref{eq:cont3}) with the constant $C$ on the right-hand side which does not
depend on $\epsilon$. Using~(\ref{eq:cont2}) and~(\ref{eq:cont3}), we obtain
\be\label{eq:prime2}
\Big|\E_{\epsilon}(t)-\E_{\epsilon}(s)+\int_s^t\D_{\epsilon}(\tau)\,d\tau\Big|\leq
C\int_s^t\sqrt{\E_{\epsilon}(\tau)}\D_{\epsilon}(\tau)\,d\tau.
\ee
In addition to that, for any $0\leq s<t\leq t^{\epsilon}$ we have
\[
\big|\E_{\epsilon}(t)-\E_{\epsilon}(s)\big|\leq C\Big|\int_s^t\D_{\epsilon}(\tau)\,d\tau\Big|
\Big(1+\sup_{0\leq s\leq t^{\epsilon}}\sqrt{\E_{\epsilon}(s)}\Big)\longrightarrow0\quad\textrm{as}\quad s\to t,
\]
since $\sup_{0\leq s\leq t^{\epsilon}}\sqrt{\E_{\epsilon}(s)}\leq\sqrt{L}$. This finishes the proof of
Theorem~{\ref{th:local}}.
\prfe
\section{Global stability}\label{ch:main}
{\em Proof of Theorem~\ref{th:globalr}}:
We exploit the estimate~(\ref{eq:prime2}) to prove the theorem.
We shall abbreviate
$(\E,\D)(u,\rho;\rho)(t)=:(\E,\D)(t)$ and $(\E_{\epsilon},\D_{\epsilon})(u^{\epsilon},\rho^{\epsilon};\rho^{\epsilon})(t)=:(\E_{\epsilon},\D_{\epsilon})(t)$.
\par
{\bf Existence.}
Let $M\leq L/2$ where $L$ is given in Lemma~\ref{lm:local3}.
Let $(u^{\epsilon},\rho^{\epsilon})$ be the associated solution to the regularized Stefan problem
on the time interval $[0,t^{\epsilon}]$ given by Theorem~\ref{th:local}.
Set
\[
{\cal T}:=\sup_t\Big\{t:\sup_{0\leq s\leq t}\E_{\epsilon}(u^{\epsilon},\rho^{\epsilon})(s)+
\int_0^t\D_{\epsilon}(u^{\epsilon},\rho^{\epsilon})(\tau)\,d\tau\leq 2M\Big\}.
\]
Theorem~\ref{th:local} guarantees ${\cal T}\geq t^{\epsilon}>0$. For
any $t<{\cal T}$, the estimate~(\ref{eq:prime2}) with $s=0$ implies
\be\label{eq:star2}
\E_{\epsilon}(t)+\int_0^t\D_{\epsilon}(\tau)\,d\tau\leq
\E_{\epsilon}(u^{\epsilon}_0,\rho^{\epsilon}_0)+ C\sup_{0\leq
s\leq{\cal
T}}\E_{\epsilon}^{\frac{1}{2}}(s)\int_0^t\D_{\epsilon}(\tau)\,d\tau,
\ee and thus \be\label{eq:imp2} \sup_{0\leq t\leq {\cal
T}}\E_{\epsilon}(t)+\int_0^{\cal T}\D_{\epsilon}(\tau)\,d\tau\leq
M+C\sqrt{2M}\int_0^{\cal T}\D_{\epsilon}(\tau)\,d\tau. \ee Choose
$M<\min\{\frac{1}{32C^2},L/2\}$. Inequality~(\ref{eq:imp2}) implies
\[
\sup_{0\leq t\leq{\cal T}}\E_{\epsilon}(t)+\int_0^{\cal T}\D_{\epsilon}(\tau)\,d\tau\leq\frac{4}{3}M<2M,
\]
which would contradict the choice of ${\cal T}$ in case ${\cal T}$ were finite. Thus ${\cal T}=\infty$ and
the estimate~(\ref{eq:star}) follows easily from~(\ref{eq:star2}) and the above choice of $M$.
This proves the theorem.
\par
{\bf Uniqueness.} We want to prove uniqueness in the class of functions $(u,\rho)$ satisfying
$\sup_{0\leq t<\infty}\E_{\epsilon}(u,\rho)(t)+\int_0^t\D_{\epsilon}(u,\rho)(\tau)\,d\tau\leq2M$,
where $M$ may be chosen smaller if necessary.
It is done in exactly the same way as the uniqueness proof in Theorem~\ref{th:local}.
\prfe
{\em Proof of Theorem~\ref{th:global}}:
\emph{Claim 1}:
Let $K\leq M$ be any positive number, where $M$ is given by Theorem~\ref{th:globalr}.
If the initial data
$(u_0,\rho_0)$ satisfy
\[
\E(u_0,\rho_0)+
\Big|\int_{\T^{n-1}}\rho_0-\int_{\Omega}u_0(1+\phi'\rho_0)\Big|
\leq\frac{K}{3},
\]
then there exists a unique global solution to the Stefan
problem~(\ref{eq:temp}) -~(\ref{eq:jump}). Moreover, we obtain the global bound
$
\sup_{0\leq t<\infty}\E(t)+\int_0^t\D(\tau)\,d\tau\leq K.
$\\
\emph{Proof of Claim 1}.
Let $\{(u^{\epsilon},\rho^{\epsilon})\}_{\epsilon}$ be a family of solutions of the
regularized Stefan problem satisfying the given
initial condition
$
(u^{\epsilon}(x,0),\rho^{\epsilon}(x',0))=(u_0^{\epsilon},\rho_0^{\epsilon}),
$
where we choose $(u_0^{\epsilon},\rho_0^{\epsilon})$ so that
\[
(u^{\epsilon}_0,\rho^{\epsilon}_0)\to(u_0,\rho_0),\quad \E(u_0^{\epsilon},\rho_0^{\epsilon};\rho^{\epsilon}_0)\to \E(u_0,\rho_0;\rho_0)
\quad\textrm{as}\quad\epsilon\to0
\]
and
\[
\sum_{|\mu|+2s\leq2k}\epsilon\int_{\T^{n-1}}|\nabla^2\d\Delta\rho_0^{\epsilon}|\leq\sqrt{\epsilon}.
\]
Thus for $\epsilon$ small, we have
$\E_{\epsilon}(u_0^{\epsilon},\rho_0^{\epsilon})
+\Big|\int_{\T^{n-1}}\rho^{\epsilon}_0-\int_{\Omega}u^{\epsilon}_0(1+\phi'\rho^{\epsilon}_0)\Big|
\leq\frac{K}{2}$. Theorem~\ref{th:globalr} guarantees global
existence of the solution $(u^{\epsilon},\rho^{\epsilon})$ and also
gives the estimate $ \sup_{0\leq t<\infty}\E_{\epsilon}(t)+
\int_0^{\infty}\D_{\epsilon}(\tau)\,d\tau\leq K. $ Since
$\E\leq\E_{\epsilon}$ and $\D\leq\D_{\epsilon}$, we obtain $
\sup_{0\leq t<\infty}\E(u^{\epsilon},\rho^{\epsilon})(t)+
\int_0^{\infty}\D(u^{\epsilon},\rho^{\epsilon})(\tau)\,d\tau\leq K.
$ Passing to the limit as $\epsilon\to0$, we obtain the solution
$(u,\rho)$ to the original Stefan problem~(\ref{eq:temp})
-~(\ref{eq:jump}). The uniqueness claim follows by setting
$\epsilon=0$ in the proof of the uniqueness statement of
Theorem~\ref{th:globalr}. This finishes the proof of {\em Claim 1}.
In the same way as we derived the inequality~(\ref{eq:star2}), we
deduce for any $t>0$: \be \label{eq:beforem1} \sup_{0\leq\tau\leq t}
\E_{\epsilon}(\tau)+ \int_0^t\D_{\epsilon}(\tau)\,d\tau \leq
\E_{\epsilon}(u^{\epsilon}_0,\rho^{\epsilon}_0)+
C\sqrt{K}\int_0^t\D_{\epsilon}(\tau)\,d\tau. \ee If we choose
$K<\min\{\frac{1}{(2C)^2},M\}=:M_1$, absorb the right-most term into
the left-hand side and drop the supremum sign, we obtain for any
$t>0$
\[
\E_{\epsilon}(t)+\frac{1}{2}\int_0^t\D_{\epsilon}(s)\,ds\leq \E_{\epsilon}(u^{\epsilon}_0,\rho^{\epsilon}_0).
\]
We let $\epsilon\to0$ and by lower semicontinuity and the assumptions on initial data, we obtain
\be\label{eq:decay1}
\E(t)+\frac{1}{2}\int_0^t\D(s)\,ds\leq \E(u_0,\rho_0).
\ee
In addition to this, we obtain the conservation law
\be\label{eq:conslaw2}
\partial_t\big\{\int_{\T^{n-1}}\rho\big\}=\partial_t\big\{\int_{\Omega}u(1+\phi'\rho)\big\}.
\ee
\par
Let us set $M^{\ast}:=\frac{M_1}{12}$, where $M_1$ is defined in the line after~(\ref{eq:beforem1}), and assume
$
\E(u_0,\rho_0)+
\Big|\int_{\T^{n-1}}\rho_0-\int_{\Omega}u_0(1+\phi'\rho_0)\Big|
\leq M^{\ast}.
$
{\em Claim 1} guarantees the global existence of the solution $(u,\rho)$ and also gives
the global bound
$
\sup_{0\leq t<\infty}\E(t)+\int_0^{\infty}\D(\tau)\,d\tau\leq \frac{M_1}{4}.
$
In order to prove~(\ref{eq:star1}), we
first fix any $s>0$. The idea is to solve the Stefan problem with the new initial data
$(u^1(x,0),\rho^1(x',0))=(u(x,s),\rho(x',s))$. The problem allows for unique solutions by~{\em Claim 1},
since
\beas
&&\E(u^1_0,\rho^1_0)+
\Big|\int_{\T^{n-1}}\rho^1_0-\int_{\Omega}u^1_0(1+\phi'\rho^1_0)\Big|\\
&&=\E(u,\rho)(s)+\Big|\int_{\T^{n-1}}\rho_0-\int_{\Omega}u_0(1+\phi'\rho_0)\Big|\leq\frac{M_1}{4}+\frac{M_1}{12}=
\frac{M_1}{3}.
\eeas
In addition to this we have the global bound
$
\sup_{0\leq t<\infty}\E(u^1,\rho^1)(t)
+\int_0^{\infty}\D(u^1,\rho^1)(\tau)\,d\tau\leq M_1
$
(again by \emph{Claim 1}).
We are thus in the uniqueness regime and we conclude $(u^1,\rho^1)(t)=(u,\rho)(t+s)$ for any $t\geq0$.
We may now use the estimate~(\ref{eq:decay1}) to obtain~(\ref{eq:star1}).
\par
The second main ingredient in proving the decay is to control the instant energy
in terms of the dissipation, i.e. to prove that there exists a constant
$C>0$ such that $\E(t)\leq C\D(t)$. We know that for $|\mu|+2s\leq2k$
\be\label{eq:auxaux}
||\partial^{\mu}_s\nabla\rho||_{H^1}\leq C\sqrt{\D}.
\ee
Thus, the only non-trivial term left to estimate
is $||u||_{L^2(\Omega)}$.
\par
\emph{Claim 2}: There exists a constant $C>0$ such that $||u||_{L^2(\Omega)}\leq C\sqrt{\D}$.\\
\emph{Proof of Claim 2}. Let $x\in\Omega$ and $x'\in\T^{n-1}$ be arbitrarily chosen. By
the mean value theorem
\[
u(x)=u(x')+\int_0^1\nabla u(tx+(1-t)x')\cdot(x-x')\,dt.
\]
Note that
$\int_{\T^{n-1}}u(x')\,dx'=0$ because $u=\nabla\cdot\Big(\frac{\nabla\rho}{\s}\Big)$ on $\T^{n-1}$.
We thus integrate with respect to $x'$ over $\T^{n-1}$ and then with respect to $x$
over $\Omega$ to obtain
\[
|\T^{n-1}|\int_{\Omega}u(x)\,dx=\int_{\Omega}\int_{\T^{n-1}}\int_0^1\nabla u(tx+(1-t)x')(x-x')\,dtdx'dx.
\]
Therefore
\beas
\Big|\int_{\Omega}u(x)\,dx\Big|&\leq&\frac{1}{|\T^{n-1}|}\max_{x\in\Omega\atop x'\in\T^{n-1}}|x-x'||\Omega|
|\T^{n-1}||\nabla u||_{L^{\infty}(\Omega)}\\
&\leq& C||\nabla u||_{L^{\infty}(\Omega)}\leq
C||\nabla u||_{H^{\frac{n}{2}+1}(\Omega)}\leq C\sqrt{\D}.
\eeas
By the Poincar\'e
inequality,
\[
||u||_{L^2(\Omega)}\leq C||\int_{\Omega}u||_{L^2(\Omega)}+C||\nabla u||_{L^2(\Omega)}\leq
C\sqrt{\D}.
\]
This finishes the proof of \emph{Claim 2}.
As explained above, \emph{Claim 2} and the estimate~(\ref{eq:auxaux}) together, imply that there
exists $C>0$ such that
\be\label{eq:decay}
\E\leq C\D.
\ee
Plugging~(\ref{eq:decay}) into~(\ref{eq:star1}) yields for any $s>0$ and some constant $\alpha>0$:
\be\label{eq:decay3}
\E(t)+\alpha\int_s^t\E(\tau)\,d\tau\leq\E(s).
\ee
As in~\cite{Ma}, p. 135, define a function $V(s):=\int_s^{\infty}\E(\tau)\,d\tau$.
From~(\ref{eq:decay3}),
$\alpha V(s)\leq\E(s)$,
\[
V'(s)=-\E(s)\leq-\alpha V(s)
\]
and thus $V(s)\leq V(0)e^{-t\alpha}$. We integrate~(\ref{eq:decay3}) with respect to $s$ over
the time interval $[t/2,t]$ to get
\[
\E(t)\frac{t}{2}\leq V(\frac{t}{2}).
\]
Thus
$
\E(t)\leq\frac{2C}{t}e^{-\frac{\alpha}{2}t}.
$
There exist $k_1,K_2>0$ such that for any $t>0$
\be\label{eq:decayx}
\E(t)\leq k_1e^{-K_2t}.
\ee
Integrating the conservation law~(\ref{eq:conslaw2}) implies $\int_{\T^{n-1}}(\rho(x',t)-\bar{\rho})\,dx'
=\int_{\Omega}u(1+\phi'\rho)$. By an argument analogous to~(\ref{eq:controlm}),
$||\rho(t)-\bar{\rho}||_2^2\leq C\E(u,\rho)(t)$. Combining this
inequality with~(\ref{eq:decayx}), we conclude
\[
\E(u,\rho)(t)+||\rho(t)-\bar{\rho}||_2^2\leq K_1e^{-K_2t}
\]
for some new constant $K_1>0$, $K_2$ as in~(\ref{eq:decayx}) and for all $t\geq 0$.
This finishes the proof of the theorem.
\prfe


\begin{thebibliography}{99}
\bibitem{Lu}
Luckhaus, S.:
Solutions for the two-phase Stefan problem with the Gibbs-Thomson law for the
melting temperature.
{\em Europ. J. Appl. Math.}
{\bf 1}, 101-111 (1990)
\end{thebibliography}
\end{document}